\documentclass[11pt,leqno,final]{amsart}
\usepackage{float}
\usepackage{source_1}

\begin{document}
\title{Fluid structure system with boundary conditions involving the pressure}
\author{Jean-J\'er\^ome Casanova}
\thanks{The author is partially supported by the ANR-Project IFSMACS (ANR 15-CE40.0010).}
\maketitle
\begin{abstract}
We study a coupled fluid--structure system involving boundary conditions on the pressure. The fluid is described by the incompressible Navier--Stokes equations in a 2D rectangular type domain where the upper part of the domain is described by a damped Euler--Bernoulli beam equation. Existence and uniqueness of local strong solutions without assumptions of smallness on the initial data is proved.
\end{abstract}


\blfootnote{Institut de Math\'ematiques de Toulouse, 
Universit\'e Paul Sabatier, 118, route de Narbonne F-31062 Toulouse Cedex 9 (France)
and School of Mathematical Sciences, Monash University, Melbourne (Australia). {\em E-mail address:} Jean-Jerome.Casanova@math.univ-toulouse.fr, jean-jerome.casanova@monash.edu}

\section{Introduction}
\subsection{Setting of the problem}
We study the coupling between the 2D Navier--Stokes equations and a damped Euler--Bernoulli beam equation in a rectangular type domain, where the beam is a part of the boundary. Let $T>0$, $L>0$ and consider the spatial domain $\Omega$ in $\mathbb{R}^2$ defined by $\Omega=(0,L)\times(0,1)$. Let us set $\Gamma_i=\{0\}\times (0,1)$ and $\Gamma_o=\{L\}\times (0,1)$ the left and right boundaries, $\Gamma_{s}=(0,L)\times\{1\}$, $\Gamma_{b}=(0,L)\times\{0\}$ and $\Gamma=\partial\Omega$ the boundary of $\Omega$. Let $\eta$ be the displacement of the beam. The function $\eta$ is defined on $\Gamma_s\times(0,T)$ with values in $(-1,+\infty)$. Let $\Omega_{\eta(t)}$ and $\Gamma_{\eta(t)}$ be the sets defined by
\[
\begin{array}{rr}
\begin{aligned}
\Omega_{\eta(t)}&=\{(x,y)\in \mathbb{R}^2\mid x\in (0,L),\, 0<y<1+\eta(x,1,t)\},\\
\Gamma_{\eta(t)}&=\{(x,y)\in \mathbb{R}^2\mid x\in (0,L),\, y=1+\eta(x,1,t)\}.\\
\end{aligned}
\end{array}
\]
\begin{figure}[H]
\centering
\begin{tikzpicture}[scale=0.65]
\draw (5,0)node [below] {$\Gamma_b$};
\draw (0,0)node [below left] {$0$};
\draw (10,0)node [below right] {$L$};
\draw (0,5)node [above] {$1$};
\draw (10,2.5)node [above right] {$\Gamma_{o}$};
\draw (0,2.5) node [above left] {$\Gamma_{i}$};
\draw (5,5) node [below] {$\Gamma_s$};
\draw (7,6.4) node [above] {$\Gamma_{\eta(t)}$};
\draw (0,0) -- (10,0);
\draw (10,0) -- (10,5);
\draw (0,5) -- (0,0);
\draw (5,5) -- (5,7);
\draw (0,5) -- (10,5);
\draw (5,5.5)node [above left] {$\eta(x,t)$};
\draw (0,5) .. controls (2,5) and (3,7) .. (5,7);
\draw (5,7) .. controls (7,7) and (8,5) .. (10,5);
\end{tikzpicture}
\caption{fluid--structure system.} 
\end{figure}
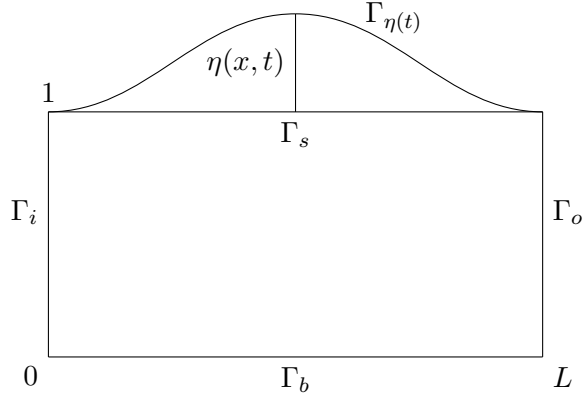
We also set  $\Gamma_{i,o}=\Gamma_{i}\cup \Gamma_{o}$. The space-time domains are denoted by
\begin{align*}
&\Sigma^{s}_{T}=\Gamma_{s}\times(0,T),\,\Sigma^{i,o}_{T}=\Gamma_{i,o}\times(0,T),\,\Sigma^{b}_{T}=\Gamma_{b}\times(0,T),\\
&\widetilde{\Sigma}_{T}^s=\underset{t\in(0,T)}{\bigcup}\Gamma_{\eta(t)}\times\{t\},\,\,\widetilde{\mathcal{Q}}_{T}=\underset{t\in(0,T)}{\bigcup}\Omega_{\eta(t)}\times\{t\}.\\
\end{align*}
We study the following fluid structure system coupling the Navier--Stokes equations and the damped Euler--Bernoulli beam equation
\begin{equation}\label{s1main}
\begin{aligned}
&\ug_{t} + (\ug\cdot\nabla)\ug-\text{div }\sigma(\ug,p)=0,\,\,\,\text{div }\ug=0\,\text{ in }\widetilde{\mathcal{Q}}_{T},\\
&\ug(x,y,t)=\eta_{t}(x,1,t)\textbf{e}_{2}\,\text{ for }(x,y,t)\in\widetilde{\Sigma}_{T}^s,\\
&u_{2}=0\,\text{ and }\,p+(1/2)\vert\ug\vert^{2}=0\,\text{ on }\Sigma^{i,o}_{T},\\
&\ug=0\,\text{ on }\Sigma^{b}_{T},\,\,\,\ug(0)=\ug^{0}\,\text{ in }\Omega_{\eta_1^0},\\
&\big[\eta_{tt}-\beta\eta_{xx}-\gamma\eta_{txx}+\alpha\eta_{xxxx}\big](x,1,t)
\\
& \hspace{1.cm}
=\psi[\ug,p,\eta](x,1+\eta(x,1,t),t)\,\text{ for }(x,t)\in(0,L)\times(0,T),\\
&\eta(\cdot,1,\cdot)=0\,\text{ and }\,\eta_{x}(\cdot,1,\cdot)=0\,\text{ on }\{0,L\}\times(0,T),\\
&\eta(\cdot,0)=\eta_1^{0}\,\text{ and }\,\eta_{t}(\cdot,0)=\eta_2^{0}\,\text{ in }\Gamma_{s},
\end{aligned}
\end{equation}
where $\ug$ is the velocity, $p$ the pressure, $\eta$ the displacement of the beam and
\begin{align*}
\sigma(\ug,p)&=-pI+\nu(\nabla\ug + (\nabla\ug)^{T}),\\
\psi[\ug,p,\eta](x,y,t)&=-\sigma(\ug,p)(x,y,t)(-\eta_{x}(x,1,t)\textbf{e}_{1}+\textbf{e}_{2})\cdot\textbf{e}_2,
\end{align*}
for all $(x,y,t)\in \widetilde{\Sigma}_{T}^s$. For a function $f$ defined on the flat domain $\Gamma_{s}$ or on $(0,L)$ we use the following abuse of notation : $f(x)=f(x,1)=f(x,y)$ for $(x,y)\in(0,L)\times\mathbb{R}$. This notation will typically be used for $f=\eta$ or $f=\psi[\ug,p,\eta]$. Hence the beam equation can be written 
\[\eta_{tt}-\beta\eta_{xx}-\gamma\eta_{txx}+\alpha\eta_{xxxx}=\psi[\ug,p,\eta]\,\text{ on }\widetilde{\Sigma}_{T}^s.\]
In the previous statement $\textbf{e}_1=(1,0)^{T}$, $\textbf{e}_2=(0,1)^{T}$, $\ug=u_1\textbf{e}_1 +u_2\textbf{e}_2$, $\nu>0$ is the constant viscosity of the fluid and $\psi$ is a force term modelling the interaction between the fluid and the beam (see  \cite{MR2745779}, \cite{MR3017292}). The constant $\beta\geq 0$, $\gamma> 0$ and $\alpha>0$ are parameters relative to the structure. This system can be used to model the blood flow through human arteries, provided that the arteries are large enough (see \cite{MR3017292}). The homogeneous Dirichlet boundary condition on $\Gamma_b$ is used to simplify the presentation; the same system with two beams can be studied in the same way.

The existence of weak solutions to system \eqref{s1main} is proved in \cite{MR3017292}. Here we would like to prove the existence of strong solutions for the same system.
A similar system is studied in \cite{MR2765696} with Dirichlet inflow and outflow boundary conditions, and in \cite{MR3466847} with periodic inflow and outflow boundary conditions.
In \cite{MR2765696} a local in time existence of strong solutions is proved without smallness assumptions on $\ug^0$ and $\eta_2^0$.
The initial condition $\eta_1^0$ is not zero, but, as far as we understand, the proof in \cite{MR2765696} is valid only if $\eta_1^0$ is small enough (see below).
Since some results in \cite{MR3466847} rely on the techniques of \cite{MR2765696}, it seems that the global existence result of \cite{MR3466847} is also only valid if $\eta^{0}_{1}$ is small. The existence of strong solutions to the fluid--structure system with a non-small $\eta^{0}_{1}$ therefore still seems to be an open question. The present paper brings an answer to this question, by establishing the local-in-time existence of strong solutions without smallness conditions on $\ug^0$, $\eta_1^0$ and $\eta_2^0$.

We prove this result for \eqref{s1main}, that is to say with boundary conditions involving the pressure. However the issue raised by a non-small $\eta^{0}_{1}$ is purely a nonlinear one, whose treatment is independent of the boundary conditions (once the proper regularity results for the linearized system have been established). 
The technique developped here, based on a novel change of variables, therefore fills the gap in \cite{MR2765696}. The existence of strong solutions for \eqref{s1main} relies on regularity results of the underlying
Stokes system and Leray projector. Three elements challenge this regularity here: the
change of variable used to deal with a generic $\eta_1^0$, the corners of the domain,
and the junctions between Dirichlet and pressure boundary conditions. To overcome
these challenges, we use symmetry techniques and a minimal-regularity transport of $H^3$ functions.
We note that, for smooth domains (no corner, no minimal-regularity change of variables),
the regularity result for the Stokes system was established in \cite{MR1900648,MR1978563}.
As a by-product of our analysis, we also obtain the existence over an arbitrary time interval $[0,T]$ of strong solutions to system \eqref{s1main}, provided that the initial data
$\ug^0$, $\eta_1^0$ and $\eta_2^0$ are small enough.

Let us detail the gap mentioned above. In \cite{MR2765696}, a key estimate, obtained through interpolation techniques, is
\begin{equation}\label{equ-lequeu}
\norme{\eta}_{L^{\infty}(\Sigma^{s}_{T})}+\norme{\eta_x}_{L^{\infty}(\Sigma^{s}_{T})}+\norme{\eta_{xx}}_{L^{\infty}(\Sigma^{s}_{T})}\leq C T^{\chi}\norme{\eta}_{H^{4,2}(\Sigma^{s}_{T})},
\end{equation}
for some $\chi>0$ and $C>0$ (see Section 3.1 for the functional spaces). If $T$ goes to $0$ the previous estimate implies that $\norme{\eta^{1}_{0}}_{L^{\infty}(\Gamma_{s})}=0$ and thus $\eta^{1}_{0}=0$. A careful study of the interpolation techniques and the Sobolev embeddings used to prove \eqref{equ-lequeu} shows that the time dependency of the constants was omitted. The fundamental reason behind this issue is related to the change of variables, used to fix the domain to $\Omega$, that introduces additional `geometrical' nonlinearities.
These nonlinearities are not small for small $T$ if the change of variables is not the identity at $T=0$. To solve this issue we rewrite the system (\ref{s1main}) in the fixed domain $\Omega_{\eta^{0}_{1}}$ instead of $\Omega$. The geometrical nonlinear terms now involve the difference $\eta-\eta^{0}_{1}$ which is small when $T$ is small.

Since our technique fills the gap in \cite{MR2765696}, this means that the global-in-time existence
result of \cite{MR3466847} for periodic boundary conditions is now genuinely established without
smallness assumption on $\eta_1^0$. An interesting question is to consider if the result
in \cite{MR3466847} can be adapted, starting from our local-in-time existence result, to obtain a global-in-time
existence of solutions with non-standard boundary conditions involving the pressure.
To do so, additional estimates should be proved to ensure that a collision between the beam and the bottom of the fluid cavity does not occur in finite time.

Finally we would like to mention some references related to our work. The Stokes and Navier--Stokes system with pressure boundary conditions was initially study in \cite{MR1308419}, using weak variational solutions. A first rigorous existence result for \eqref{s1main} with periodic boundary conditions goes back to \cite{MR2027753} where an iterative method was used to handle the coupled system. The feedback stabilization of \eqref{s1main} with Dirichlet inflow and outflow boundary conditions is studied in \cite{MR2745779} and provides a semigroup approach for the linearized system, based on a splitting of the pressure, that is used in the present article.  This semigroup framework was already used in \cite{MR2359448,MR2379673} for a linear model.
\subsection{Main results}
The main result of this paper is Theorem \ref{s42thm2} which proves the existence of a unique local strong solution for the fluid--structure system (\ref{s1main}) without smallness assumptions on the initial data. We also state in Theorem \ref{s42thm3} the existence of a unique strong solution on the time interval $[0,T]$ with $T>0$ an arbitrary fixed time, for small enough initial data. Several changes of variables are done on (\ref{s1main}) and these results are given for equivalent system (see (\ref{s32e2})).

The structure of the article is as follows. In Section 2, we rewrite (\ref{s1main}) in a fixed domain and we explain ideas of the proof which consist in studying a linear system associated with (\ref{s1main}) and in using a fixed point argument. In Section 3 we eliminate the pressure in the beam equation by expressing it in terms the velocity. We then rewrite the system as an abstract evolution equation and we prove that the underlying operator generates an analytic semigroup. Finally we prove the nonlinear estimates, with explicit time dependency, in Section 4 and we conclude with a fixed point procedure. All this process is based on the extension to non-standard boundary conditions of known result on the Stokes equations. This is detailed in the appendix.
\section{Plan of the paper}
\subsection{Equivalent system in a reference configuration} In order to study system \eqref{s1main}, we are going to rewrite the system in a reference configuration which can be chosen arbitrarily. For that,  throughout what follows, we choose a function $\eta^0$  belonging to $H^3(\Gamma_s)\cap H_0^2(\Gamma_s)$, 
and satisfying $1+\eta^0(x) >0$ for all $x\in (0,L)$. 
Set 
\begin{equation}\label{domain-ref}
\begin{aligned}
\Omega_0 &{}= \{(x,y)\in {\mathbb R}^2 \mid x\in (0,L),\ 0<y<1+\eta^0(x)\},\quad &Q_T=\Omega_0\times (0,T), \\
\Gamma_0 &{}= \{(x,y)\in {\mathbb R}^2 \mid x\in (0,L),\ y=1+\eta^0(x)\}, &\Sigma_T^0=\Gamma_0\times (0,T),
\end{aligned}
\end{equation}
$\Gamma_{d}=\Gamma_{0}\cup\Gamma_{b}$ and $\Sigma^{d}_{T}=\Gamma_{d}\times(0,T)$. 
In order to rewrite  system (\ref{s1main})  in the cylindrical  domain $Q_{T}$  for all $t\in (0,T)$ consider the following diffeomorphism 
\begin{equation}\label{s31changvar}
\mathcal{T}_{\eta(t)} :
\begin{cases}
\Omega_{\eta(t)}&\longrightarrow \Omega_{0},\\
(x,y)&\longmapsto (x,z)=\left(x,\frac{1+\eta^{0}(x)}{1+\eta(x,t)}y\right).\\
\end{cases}
\end{equation}
The variable $z$ can be written under the form $z=\displaystyle\frac{y}{1+\widetilde{\eta}}$ with $\widetilde{\eta}=\frac{\eta-\eta^{0}}{1+\eta^{0}}$. We introduce the new unknowns
$$
\widehat{\textbf{u}}(x,z,t)=\ug(\mathcal{T}_{\eta(t)}^{-1}(x,z) ,t) 
\quad\mbox{and}\quad 
\widehat{p}(x,z,t)=p(\mathcal{T}_{\eta(t)}^{-1}(x,z),t),
$$
and we set $\widehat{\textbf{u}}^0(x,z) = \textbf{u}^0(\mathcal{T}_{\eta_1^0}^{-1}(x,z) ) $. 
With this change of variables,
$$
p(x,1+\eta(x,t),t) = \widehat{p}(x,1+\eta^0(x,t),t)\quad \mbox{and}\quad \widehat{\textbf{u}}(x,1+\eta(x,t),t) = {\textbf{u}}(x,1+\eta^0(x,t),t),
$$
for all $(x,t)\in (0,L)\times (0,T)$. 
The system satisfied by $(\widehat{\textbf{u}},\widehat{p}, \eta)$ is
\begin{equation}\label{s31e1}
\begin{aligned}
&\widehat{\textbf{u}}_{t}-\nu\Delta\widehat{\textbf{u}}+\nabla\widehat{p}=\textbf{G}(\widehat{\textbf{u}},\widehat{p},\eta),\,\,\,\text{div }\widehat{\textbf{u}}=\text{div }\textbf{w}(\widehat{\textbf{u}},\eta)\,\text{ in }Q_{T},\\
&\widehat{\textbf{u}}=\eta_{t}\textbf{e}_{2}\,\text{ on }\Sigma^{0}_{T},\\
&\widehat{u}_2=0\,\text{ and }\,\widehat{p}+(1/2)\vert\widehat{\textbf{u}}\vert^{2}=0\,\text{ on }\Sigma^{i,o}_{T},\\
&\widehat{\textbf{u}}=0\,\text{ on }\Sigma^{b}_{T},\,\,\,\widehat{\textbf{u}}(0)=\widehat{\textbf{u}}^{0}\,\text{ on }\Omega_{0},\\
&\eta_{tt}-\beta\eta_{xx}-\gamma\eta_{txx}+\alpha\eta_{xxxx}=\widehat{p}+{\Psi}(\widehat{\textbf{u}},\eta)\,\text{ on }\Sigma^{0}_{T},\\
&\eta=0\,\text{ and }\,\eta_{x}=0\,\text{ on }\{0,L\}\times(0,T),\\
&\eta(0)=\eta_1^{0}\,\text{ and }\,\eta_{t}(0)=\eta_2^{0}\,\text{ in }\Gamma_{s},
\end{aligned}
\end{equation}
with 
\begin{align*}
\textbf{G}(\widehat{\textbf{u}},\widehat{p},\eta)&=-\widetilde{\eta}\widehat{\textbf{u}}_t+\left[z\widetilde{\eta}_t+\nu z\left(\frac{\widetilde{\eta}_x^2}{1+\widetilde{\eta}}-\widetilde{\eta}_{xx}\right)\right]\widehat{\textbf{u}}_z\\
&+\nu\left[-2z\widetilde{\eta}_x\widehat{\textbf{u}}_{xz}+\widetilde{\eta}\widehat{\textbf{u}}_{xx}+\frac{z^2\widetilde{\eta}_x^2-\widetilde{\eta}}{1+\widetilde{\eta}}\widehat{\textbf{u}}_{zz}\right]\\
&+z(\widetilde{\eta}_x\widehat{p}_z-\widetilde{\eta}\widehat{p}_x)\textbf{e}_1-(1+\widetilde{\eta})\widehat{u}_1
\widehat{\textbf{u}}_x+(z\widetilde{\eta}_x
\widehat{u}_1-\widehat{u}_2)\widehat{\textbf{u}}_z,
\end{align*}
\[\textbf{w}[\widehat{\textbf{u}},\eta]=-\widetilde{\eta}\widehat{u}_1\textbf{e}_1+z\widetilde{\eta}_x\widehat{u}_1\textbf{e}_2,\]
and
\[
{\Psi}(\widehat{\textbf{u}},\eta)=\nu\left(\frac{\eta_x}{1+\widetilde{\eta}}\widehat{u}_{1,z}+\eta_x\widehat{u}_{2,x}-\frac{z\widetilde{\eta}_x\eta_{x}-2}{1+\widetilde{\eta}}\widehat{u}_{2,z}\right).
\]
In Section 5, in order to prove the existence of solution to system \eqref{s32e2}, derived from system \eqref{s31e1} by a change of variables, we assume that $\eta_1^0$ is equal to $\eta^0$. In that special case, the function $\widetilde{\eta}$ is equal to $0$ at time $t=0$ which implies that $\mathcal{T}_{\eta(0)}$ is the identity. We also obtain that $\textbf{w}(\widehat{\ug},\eta)$ is equal to $0$ at time $t=0$. But up to Section 5 and in the appendix, $\eta^0$ is chosen a priori, and not necessarily equal to $\eta_1^0$. 
\subsection{Final system and linearisation}
In order to come back to a divergence free system consider the function $\overline{\ug}$ defined by $\overline{\ug} = \widehat{\textbf{u}} -\textbf{w}(\widehat{\textbf{u}},\eta)$. Set
\begin{equation}\label{defA}
M(\overline{\ug},\eta)=
\begin{pmatrix}
\displaystyle\frac{\overline{u}_{1}}{1+\widetilde{\eta}}\\
\displaystyle\frac{z\widetilde{\eta}_{x}}{1+\widetilde{\eta}}\overline{u}_{1}+\overline{u}_{2}\\
\end{pmatrix}\quad
\mbox{and}\quad 
N(\overline{\ug},\eta)=\begin{pmatrix}\displaystyle\frac{-\widetilde{\eta}\overline{u}_{1}}{1+\widetilde{\eta}}\\
\displaystyle\frac{z\widetilde{\eta}_{x}\overline{u}_{1}}{1+\widetilde{\eta}}
\end{pmatrix}. 
\end{equation}
The function $\widehat{\textbf{u}}$ can be expressed in terms of 
$\overline{\ug}$ as follows 
$$
\widehat{\textbf{u}} = M(\overline{\ug},\eta) = \overline{\ug} + N(\overline{\ug},\eta). 
$$
To simplify the notation, we drop out the hat over $p$. 
Thus the system satisfied by  $(\overline{\ug},\widehat{p},\eta)=(\overline{\ug},{p},\eta)$ is 
\begin{equation}\label{s32e1}
\begin{aligned}
&\overline{\ug}_{t}-\text{div }\sigma(\overline{\ug},p)=\textbf{F}(\overline{\ug},p,\eta),\,\,\,\text{div }\overline{\ug}=0\,\text{ in }Q_{T},\\
&\overline{\ug}=\eta_{t}\textbf{e}_{2}-\textbf{w}(M(\overline{\ug},\eta),\eta)\,\text{ on }\Sigma^{0}_{T},\\
&\overline{u}_2=-w_2(M(\overline{\ug},\eta),\eta)\,\text{ and }\, p+(1/2)\vert\overline{\ug} +\textbf{w}(M(\overline{\ug},\eta),\eta)\vert^2=0\,\text{ on }\Sigma^{i,o}_{T},\\
&\overline{\ug} =-\textbf{w}(M(\overline{\ug},\eta),\eta)\,\text{ on }\Sigma^{b}_{T},\,\,\,\overline{\ug}(0)=\widehat{\textbf{u}}^{0}-
\textbf{w}(M(\overline{\ug},\eta),\eta)(0)\,\text{ in }\Omega_{0},\\
&\eta_{tt}-\beta\eta_{xx}-\gamma\eta_{txx}+\alpha\eta_{xxxx}=p+\Psi(M(\overline{\ug},\eta),\eta)\,\text{ on }\Sigma^{0}_{T},\\
&\eta=0\,\text{ and }\,\eta_{x}=0\,\text{ on }\{0,L\}\times(0,T),\\
&\eta(0)=\eta_1^{0}\,\text{ and }\,\eta_{t}(0)=\eta_2^{0}\,\text{ in }\Gamma_{s},
\end{aligned}
\end{equation}
with $\textbf{F}(\overline{\ug},p,\eta)=\textbf{G}(M(\overline{\ug},\eta),p,\eta)-\partial_tN(\overline{\ug},\eta) +\nu\Delta 
N(\overline{\ug},\eta)$.  

Recall that $\textbf{w}(\widehat{\textbf{u}},\eta)=-\widetilde{\eta} \widehat{u}_1\textbf{e}_1+z\widetilde{\eta}_x \widehat{u}_1\textbf{e}_2$. Since $\widehat{u}_1=0$ on $\Sigma^{d}_{T}$ and $\widetilde{\eta}=\widetilde{\eta}_{x}=0$ on $\{0,L\}\times(0,T)$, we have  $\textbf{w}(\widehat{\textbf{u}},\eta)= 0$ on $\partial\Omega_{0}\times(0,T)$. System (\ref{s32e1}) becomes
\begin{equation}\label{s32e2}
\begin{aligned}
&\overline{\ug}_{t}-\text{div }\sigma(\overline{\ug},p)=\textbf{F}(\overline{\ug},p,\eta),\,\,\,\text{div }\overline{\ug}=0\,\text{ in }Q_{T},\\
&\overline{\ug}=\eta_{t}\textbf{e}_{2}\,\text{ on }\Sigma^{0}_{T},\\
&\overline{u}_2=0\,\text{ and }\,p+(1/2)\vert\overline{\ug} \vert^2=0\,\text{ on }\Sigma^{i,o}_{T},\\
&\overline{\ug} =0\,\text{ on }\Sigma^{b}_{T},\,\,\,\overline{\ug}(0)=\overline{\ug}^{0}\,\text{ in }\Omega_{0},\\
&\eta_{tt}-\beta\eta_{xx}-\gamma\eta_{txx}+\alpha\eta_{xxxx}=p+H(\overline{\ug},\eta)\,\text{ on }\Sigma^{0}_{T},\\
&\eta=0\,\text{ and }\,\eta_{x}=0\,\text{ on }\{0,L\}\times(0,T),\\
&\eta(0)=\eta_1^{0}\,\text{ and }\,\eta_{t}(0)=\eta_2^{0}\,\text{ in }\Gamma_{s},
\end{aligned}
\end{equation}
with $H(\overline{\ug},\eta)=\Psi(M(\overline{\ug},\eta),\eta)$ and 
$\overline{\ug}^{0}=\widehat{\textbf{u}}^{0} - \textbf{w}(\widehat{\textbf{u}},\eta)(0)$.

In order to solve the system (\ref{s32e2}) with a fixed point argument, consider the following linear system
\begin{equation}\label{s32e3}
\begin{aligned}
&\textbf{u}_{t}-\text{div }\sigma(\textbf{u},p)={\textbf{f}},\,\,\,\text{div }\textbf{u}=0\,\text{ in }Q_{T},\\
&\textbf{u}=\eta_{t}\textbf{e}_{2}\,\text{ on }\Sigma^{0}_{T},\\
& u_2=0\,\text{ and }\,p=\Theta\,\text{ on }\Sigma^{i,o}_{T},\\
&\ug=0\,\text{ on }\Sigma^{b}_{T},\,\,\,\textbf{u}(0)=\ug^{0}\,\text{ on }\Omega_{0},\\
&\big[\eta_{tt}-\beta\eta_{xx}-\gamma\eta_{txx}+\alpha\eta_{xxxx}\big](x,t)
\\
& \hspace{1.cm}
=p(x,1+\eta^0(x,t),t)+h\,\text{ for }(x,t)\in(0,L)\times(0,T),\\
&\eta=0\,\text{ and }\,\eta_{x}=0\,\text{ on }\{0,L\}\times(0,T),\\
&\eta(0)=\eta^{0}_{1}\,\text{ and }\,\eta_{t}(0)=\eta_2^{0}\,\text{ in }\Gamma_{s},
\end{aligned}
\end{equation}
with $\textbf{f}\in L^{2}(0,T;\textbf{L}^{2}(\Omega_{0}))$, $h\in L^{2}(0,T;L^{2}(\Gamma_s))$ and $\Theta\in L^{2}(0,T;H^{1/2}(\Gamma_{i,o}))$.
\section{Linear system}
Recall that $\Omega_{0}$ is given by \eqref{domain-ref} with a fixed $\eta^{0}\in H^{3}(\Gamma_{s})\cap H^{2}_{0}(\Gamma_{s})$ such that $1+\eta^{0}(x)>0$ for all $x\in (0,L)$.
\subsection{Function spaces}
To deal with the mixed boundary condition for the Stokes system
\begin{equation}\label{Stokes-def}
\begin{aligned}
&-\nu\Delta \ug +\nabla p=\textbf{f},\,\,\,\,\text{div }\ug=0\,\text{ in }\Omega_{0},\\
&\textbf{\ug}=0\,\text{ on }\Gamma_{d},\,\,\,u_{2}=0\,\text{ and }p=0\,\text{ on }\Gamma_{i,o},\\
\end{aligned}
\end{equation}
introduce the space
\[\textbf{V}^{0}_{n,\Gamma_{d}}(\Omega_{0})=\{\vg\in \textbf{L}^{2}(\Omega_{0})\mid \text{div }\vg=0\text{ in }\Omega_{0}, \vg\cdot\textbf{n}=0\text{ on }\Gamma_{d}\},\]
and the orthogonal decomposition of $\textbf{L}^{2}(\Omega_{0})=L^{2}(\Omega_{0},\mathbb{R}^{2})$
\[\textbf{L}^{2}(\Omega_{0})=\textbf{V}^{0}_{n,\Gamma_{d}}(\Omega_{0})\oplus\text{grad }H^{1}_{\Gamma_{i,o}}(\Omega_{0}),\]
where $H^{1}_{\Gamma_{i,o}}(\Omega_{0})=\{u\in H^{1}(\Omega_{0})\mid u=0 \text{ on }\Gamma_{i,o}\}$.  Let $\Pi :\textbf{L}^{2}(\Omega_{0})\rightarrow \textbf{V}^{0}_{n,\Gamma_{d}}(\Omega_{0})$ be the so-called Leray projector associated with this decomposition. If $\textbf{u}$ belongs to $\textbf{L}^{2}(\Omega_{0})$ then $\Pi\textbf{u}=\textbf{u}-\nabla p_{\textbf{u}}-\nabla q_{\textbf{u}}$ where $p_{\textbf{u}}$ and $q_{\textbf{u}}$ are solutions to the following elliptic equations
\begin{equation}\label{equation-Pif}
\begin{aligned}
&p_{\textbf{u}}\in H^{1}_{0}(\Omega_{0}),\,\,\,\Delta p_{\textbf{u}}=\text{div }\textbf{u}\in H^{-1}(\Omega_{0}),\\
&q_{\textbf{u}}\in H^{1}_{\Gamma_{i,o}}(\Omega_{0}),\,\,\,\Delta q_{\textbf{u}}=0,\,\,\,\frac{\partial q_{\textbf{u}}}{\partial \normal}=(\textbf{u}-\nabla p_{\textbf{u}})\cdot\normal\,\text{ on }\Gamma_{d},\,\,\,q_{\textbf{u}}=0\,\text{ on }\Gamma_{i,o}.
\end{aligned}
\end{equation}
Through this paper the functions with vector values are written with a bold typography. For example $\textbf{H}^{2}(\Omega_{0})=H^{2}(\Omega_{0},\mathbb{R}^{2})$. As the boundary $\Gamma_{0}$ is not $\mathcal{C}^{2,1}$ it is not clear that the operator $\Pi$ preserves the $\textbf{H}^{2}$-regularity. However, with extra conditions on $\textbf{u}$, this can be proved. Using the notations in \cite[Theorem 11.7]{MR0350177} we introduce the space $H^{3/2}_{00}(\Gamma_{0})=[H^{1}_{0}(\Gamma_{0}),H^{2}_{0}(\Gamma_{0})]_{1/2}$. The following lemma is proved in the appendix.
\begin{lem}\label{lemPif}Let $\ug$ be in $\textbf{H}^{2}(\Omega_{0})$ satisfying $\text{div }\ug=0$, $\ug=0$ on $\Gamma_{b}$ and $\ug=g\textbf{e}_{2}$ on $\Gamma_{0}$ with $g\in H^{3/2}_{00}(\Gamma_{0})$. Then $\Pi\ug$ belongs to $\textbf{H}^{2}(\Omega_{0})$.
\end{lem}
The energy space associated with (\ref{Stokes-def}) is
\[V=\{\ug\in \textbf{H}^{1}(\Omega_{0})\mid \text{div }\ug=0\text{ in }\Omega_{0}\text{, }\ug=0\text{ on }\Gamma_{d}\text{, }u_2=0\text{ on }\Gamma_{i,o}\}.\]
The regularity result for (\ref{Stokes-def}) stated in Theorem \ref{s61thm5} in the appendix allows us to introduce the Stokes operator $A$ defined in $\textbf{V}^{0}_{n,\Gamma_{d}}(\Omega_{0})$ by
\[\mathcal{D}(A)=\textbf{H}^{2}(\Omega_{0})\cap V,\]
and for all $\ug\in \mathcal{D}(A)$, $A\ug=\nu\Pi\Delta\ug$. We also use the notations
\[\textbf{V}^{s}(\Omega_{0})=\{\ug\in\textbf{H}^{s}(\Omega_{0})\mid \text{div }\ug=0\},\,\,\,\textbf{V}^{s}_{n,\Gamma_{d}}(\Omega_{0})=\textbf{V}^{0}_{n,\Gamma_{d}}(\Omega_{0})\cap \textbf{H}^{s}(\Omega_{0}),\]
for $s\geq 0$. For the Dirichlet boundary condition on $\Gamma_{0}$ set
\begin{align*}
\mathcal{L}^{2}(\Gamma_{0})&=\{0\}\times L^{2}(\Gamma_{0}),&\mathcal{H}^{3/2}_{00}(\Gamma_{0})&=\{0\}\times H^{3/2}_{00}(\Gamma_{0}),\\
\mathcal{H}^{s}(\Gamma_{0})&=\{0\}\times H^{s}(\Gamma_{0}),&\mathcal{H}^{s}_{0}(\Gamma_{0})&=\{0\}\times H^{s}_{0}(\Gamma_{0})
\,\text{ for }s\geq 0.
\end{align*}
For $s\geq 0$, the dual space of $\mathcal{H}^{s}(\Gamma_{0})$ with $\mathcal{L}^{2}(\Gamma_{0})$ as pivot space is denoted by $(\mathcal{H}^{s}(\Gamma_{0}))'$. Let $D\in \mathcal{L}(\mathcal{H}^{3/2}_{00}(\Gamma_{0}),\textbf{H}^{2}(\Omega_{0}))$ be the operator defined by $D\textbf{g}=\textbf{w}$ where $(\textbf{w},p)$ is the solution to
\[
\begin{aligned}
&-\nu\Delta \wg +\nabla p=0,\,\,\,\,\text{div }\wg=0\,\text{ in }\Omega_{0},\\
&\textbf{\wg}=\textbf{g}\,\text{ on }\Gamma_{0},\,\,\,\textbf{\wg}=0\,\text{ on }\Gamma_{b},\,\,\,w_{2}=0\,\text{ and }p=0\,\text{ on }\Gamma_{i,o}.\\
\end{aligned}
\]
given by Theorem \ref{s61thm5}. Using a weak regularity result (Theorem \ref{s61thm6}) and interpolation techniques (see \cite[Theorem 12.6]{MR0350177} and \cite[Remark 12.6]{MR0350177}), $D$ can be extended as a bounded linear operator in $\mathcal{L}(\mathcal{L}^{2}(\Gamma_{0}),\textbf{H}^{1/2}(\Omega_{0}))$.

For space-time dependent functions we use the notations introduced in \cite{MR0350178}:
\begin{align*}
\textbf{L}^{2}(Q_{T})&{}=L^{2}(0,T;\textbf{L}^{2}(\Omega_{0})),\,\,\,\textbf{H}^{p,q}(Q_{T})=L^{2}(0,T;\textbf{H}^{p}(\Omega_{0}))\cap H^{q}(0,T;\textbf{L}^{2}(\Omega_{0})),\,p,q\geq 0,\\
L^{2}(\Sigma_{T}^{s})&{}=L^{2}(0,T;L^{2}(\Gamma_{s})),\,\,\,\,H^{p,q}(\Sigma^{s}_{T})=L^{2}(0,T;H^{p}(\Gamma_{s}))\cap H^{q}(0,T;L^{2}(\Gamma_{s})),\,\,\,p,q\geq 0.
\end{align*}
\subsection{Semigroup formulation of the linear system}
We want to prove existence and regularity results for the coupled linear system (\ref{s32e3}). Let $\mathcal{R}\in \mathcal{L}(H^{1/2}(\Gamma_{i,o}),H^{1}(\Omega))$ be a lifting operator. Classically we transform (\ref{s32e3}) into a system with homogeneous boundary conditions (for the pressure) by looking for a solution to (\ref{s32e3}) under the form $(\ug,p,\eta)=(\ug,p_{1},\eta)+(0,\mathcal{R}(\Theta),0)$ with $(\ug,p_{1},\eta)$ solution to
\begin{equation}\label{s33e1}
\begin{aligned}
&\textbf{u}_{t}-\text{div }\sigma(\textbf{u},p_{1})=\textbf{f}-\nabla\mathcal{R}(\Theta),\,\,\,\text{div }\textbf{\ug}=0\,\text{ in }Q_{T},\\
&\textbf{u} =\eta_{t}\textbf{e}_{2}\,\text{ for }\Sigma^{0}_{T},\\
&u_2=0\,\text{ and }\,p_{1}=0\,\text{ on }\Sigma^{i,o}_{T},\\
&\textbf{\ug} =0\,\text{ on }\Sigma^{b}_{T},\,\,\,\textbf{u}(0)=\ug^{0}\,\text{ in }\Omega_{0},\\
&\big[\eta_{tt}-\beta\eta_{xx}-\gamma\eta_{txx}+\alpha\eta_{xxxx}\big](x,t)
\\
& \hspace{1.cm}
=\Big[p_{1}+\mathcal{R}(\Theta)\Big](x,1+\eta^0(x,t),t)+h\,\text{ for }(x,t)\in(0,L)\times(0,T),\\
&\eta=0\,\text{ and }\,\eta_{x}=0\,\text{ on }\{0,L\}\times(0,T),\\
&\eta(0)=\eta^{0}_{1}\,\text{ and }\,\eta_{t}(0)=\eta^{0}_{2}\,\text{ in }\Gamma_{s}.
\end{aligned}
\end{equation}
Set $F=\textbf{f}-\nabla\mathcal{R}(\Theta)$. As the boundary $\Gamma_{0}$ may not be flat and the beam equation is written on $\Gamma_{s}$ consider the transport operator $\mathcal{U}\in {\mathcal L}(L^{2}(\Gamma_{0}),L^{2}(\Gamma_{s}))$ defined by
$$
(\mathcal{U}g)(x,1)=g(x,1+\eta^{0}(x))\quad \mbox{for all}\ g\in L^{2}(\Gamma_{0}).
$$
We can easily verify that $\mathcal{U}$ is an isomorphism from $L^{2}(\Gamma_{0})$ to $L^{2}(\Gamma_{s})$, and that 
$$
(\mathcal{U}^{-1}\widetilde g)(x,1+\eta^{0}(x))=\widetilde g(x,1)\quad \mbox{for all}\ \widetilde g\in L^{2}(\Gamma_{s}).
$$
Moreover $\mathcal{U}^{-1}=\mathcal{U}^*$, if $L^{2}(\Gamma_{0})$ and $L^{2}(\Gamma_{s})$ are equipped with the inner products 
$$
\begin{array}{l}
\left(f,g \right)_{L^{2}(\Gamma_{0})}= \left(f(\cdot,1+\eta^{0}(\cdot))g(\cdot,1+\eta^{0}(\cdot)) \right)_{L^{2}(0,L)},
\\[2.mm] 
\mbox{and}
\\[2.mm] 
\left(\widetilde f,\widetilde g \right)_{L^{2}(\Gamma_{s})}=\left(\widetilde f(\cdot,1) \widetilde g(\cdot,1) \right)_{L^{2}(0,L)}. 
\end{array}
$$
In order to express the pressure $p_{1}$ in terms of $\Pi\ug$ and $\eta$ we introduce the Neumann-to-Dirichlet operator $N_{s}\in\mathcal{L}(L^{2}(\Gamma_{0}))$ defined by $N_{s,0}(g)=\pi|_{\Gamma_{0}}$ where $g\in L^{2}(\Gamma_{0})$ and $\pi$ is the solution to
\[
\begin{cases}
\begin{aligned}
&\Delta \pi=0\,\text{ in }\Omega_{0},\\
&\pi=0\,\text{ on }\Gamma_{i,o},\,\,\,\frac{\partial \pi}{\partial \normal}=g(1+(\eta^{0})^{2})^{-1/2}\,\text{ on }\Gamma_{0}\,\text{ and }\,\frac{\partial \pi}{\partial \normal}=0\,\text{ on }\Gamma_{b}.\\
\end{aligned}
\end{cases}
\]
As in \cite[Lemma 3.1]{MR2745779}, $N_{s,0}$ is a non-negative symmetric and compact operator in $L^{2}(\Gamma_{0})$. Hence, as $\mathcal{U}\in {\mathrm{isom}}(L^{2}(\Gamma_{0}),L^{2}(\Gamma_{s}))$ and $\mathcal{U}^{-1}=\mathcal{U}^*$, the operator 
$N_{s}=\mathcal{U}N_{s,0}\,\mathcal{U}^{-1}$ is a non-negative symmetric and compact operator in $L^{2}(\Gamma_{s})$. Consequently, 
the operator $(I_{L^{2}(\Gamma_{s})}+N_{s})$ is an automorphism in $L^{2}(\Gamma_{s})$. 

We also define the operator $N_{0}\in \mathcal{L}(H^{-1/2}(\Gamma_{d}),L^{2}(\Gamma_{s}))$ by $N_{0}(v)=\mathcal{U}\left(\rho|_{\Gamma_{0}}\right)$ for all $v\in H^{-1/2}(\Gamma_{d})$, where $\rho$ is the solution to
\[
\begin{cases}
\begin{aligned}
&\Delta \rho =0\,\text{ in }\Omega_0,\\
&\rho=0\,\text{ on }\Gamma_{i,o}\,\text{ and }\,\frac{\partial \rho}{\partial \normal}=v\,\text{ on }\Gamma_{d}.\\
\end{aligned}
\end{cases}
\]
Finally set $D_{s}(\eta_{t})=D(\mathcal{U}^{-1}\eta_{t}\textbf{e}_{2})$. The following lemma is similar to \cite[Lemma 3.1]{MR2745779} and is a direct application of Theorem \ref{s62thm5} in the appendix.
\begin{lem}\label{lemmaP}A pair $(\ug,p_{1})\in\textbf{H}^{2,1}(Q_{T})\times L^{2}(0,T;H^{1}(\Omega_{0}))$ obeys the fluid equations of (\ref{s33e1}) if and only if
\begin{align*}
&\Pi\ug'=A\Pi\ug + (-A)\Pi D_{s}(\eta_t)+\Pi F,\,\,\,\ug(0)=\ug^{0},\\
&(I-\Pi)\ug=(I-\Pi)D_{s}(\eta_t),\,\,\,p_{1}=\rho -q_{t}+p_{F},
\end{align*}
where 
\begin{itemize}
\item $q\in H^{1}(0,T;H^{1}(\Omega_{0}))$ is the solution to
\[\Delta q=0\,\text{ in }Q_{T},\,\,\,\rho=0\,\text{ on }\Sigma^{i,o}_{T},\,\,\,\frac{\partial q}{\partial \normal}=\mathcal{U}^{-1}\eta_{t}\textbf{e}_{2}\cdot\normal\,\text{ on }\Sigma^{0}_{T},\,\,\,\frac{\partial q}{\partial\normal}=0\,\text{ on }\Sigma^{b}_{T}.\]
\item $\rho\in L^{2}(0,T;H^{1}(\Omega_{0}))$ is the solution to
\[\Delta \rho =0\,\text{ in }Q_{T},\,\,\,\rho=0\,\text{ on }\Sigma^{i,o}_{T},\,\,\,\frac{\partial\rho}{\partial\normal}=\nu \Delta\Pi\ug\cdot\normal\,\text{ on }\Sigma^{d}_{T}.\]
\item $p_{F}\in L^{2}(0,T;H^{1}(\Omega_{0}))$ is given by the identity $(I-\Pi)F=\nabla p_{F}$.
\end{itemize}
\end{lem}
Using Lemma \ref{lemmaP} the pressure in the beam equation can be decomposed as follows $p_{1}=\nu N_{0}(\Delta\Pi\ug\cdot\normal)-\partial_{t}N_{s}(\eta_{t})+\mathcal{U}(p_{F}|_{\Gamma_{0}})$. Hence the beam equation becomes
\[(I_{L^{2}(\Gamma_{s})}+N_{s})\eta_{tt}-\beta\eta_{xx}-\gamma\eta_{txx}+\alpha\eta_{xxxx}=\nu N_{0}(\Delta\Pi\ug\cdot\normal)+\mathcal{U}\big[(p_{F}+\mathcal{R}(\Theta))|_{\Gamma_{0}}\big]+h.\]
The system (\ref{s33e1}) can be rewritten in terms of $(\Pi\ug,\eta,\eta_t)=(\Pi\ug,\eta_{1},\eta_{2})$ as
\begin{equation}\label{s33e2}
\begin{cases}
\begin{aligned}
&\frac{d}{dt}\begin{pmatrix}
\Pi\ug\\
\eta_{1}\\
\eta_{2}\\
\end{pmatrix}
=\mathcal{A}
\begin{pmatrix}
\Pi\ug\\
\eta_{1}\\
\eta_{2}\\
\end{pmatrix}
+\textbf{F},\\
&(I-\Pi)\ug=(I-\Pi)D_s(\eta_t),\\
\end{aligned}
\end{cases}
\end{equation}
where $\mathcal{A}$ is the unbounded operator in
\[\textbf{H}=\textbf{V}^{0}_{n,\Gamma_{d}}(\Omega_{0})\times H^{2}_{0}(\Gamma_s)\times L^{2}(\Gamma_s),\]
with domain
\[\mathcal{D}(\mathcal{A})=\{(\Pi\ug,\eta_1,\eta_2)\in \textbf{V}^{2}_{n,\Gamma_{d}}(\Omega_{0})\times (H^{4}(\Gamma_s)\cap H^{2}_{0}(\Gamma_s))\times H^{2}_{0}(\Gamma_s)\mid \Pi\ug -\Pi D_s(\eta_2)\in \mathcal{D}(A)\},\]
defined by
\begin{equation}\label{s5matriceA}
\mathcal{A}=
\begin{pmatrix}
I&0&0\\
0&I&0\\
0&0&(I+N_s)^{-1}\\
\end{pmatrix}
\begin{pmatrix}
A&0&(-A)\Pi D_s\\
0&0&I\\
\nu N_{0}(\Delta(\cdot)\cdot\normal)&\beta\Delta_s -\alpha\Delta^{2}_{s}&\delta\Delta_s\\
\end{pmatrix},
\end{equation}
with $\Delta_s=\partial_{xx}$ and 
\[\textbf{F}=
\begin{pmatrix}
\Pi F\\
0\\
(I+N_s)^{-1}(\mathcal{U}\big[(p_{F}+\mathcal{R}(\Theta))|_{\Gamma_{0}}\big]+h)\\
\end{pmatrix}.
\]
\subsection{Analyticity of $\mathcal{A}$}
Let $(A_{\alpha,\beta},\mathcal{D}(A_{\alpha,\beta}))$ be the unbounded operator in $L^{2}(\Gamma_s)$ defined by $\mathcal{D}(A_{\alpha,\beta})=H^{4}(\Gamma_s)\cap H^{2}_{0}(\Gamma_s)$ and for all $\eta\in\mathcal{D}(A_{\alpha,\beta})$, $A_{\alpha,\beta}\eta=\beta\eta_{xx}-\alpha\eta_{xxxx}$. The operator $A_{\alpha,\beta}$ is self-adjoint and is an isomorphism from $\mathcal{D}(A_{\alpha,\beta})$ to $L^{2}(\Gamma_s)$. The space $\textbf{H}$ will be equipped with the inner product
\[\langle (\textbf{u},\eta_1,\eta_2),(\vg,\zeta_1,\zeta_2)\rangle_{\textbf{H}}=\langle \ug,\vg\rangle_{\textbf{V}^{0}_{n,\Gamma_{d}}(\Omega_{0})}+\langle \eta_1,\zeta_1\rangle_{H^{2}_{0}(\Gamma_s)}+\langle \eta_2,\zeta_2 \rangle_{L^{2}(\Gamma_s)},\]
with $\textbf{V}^{0}_{n,\Gamma_d}(\Omega_{0})$ endowed with the natural scalar product of $L^{2}(\Omega_{0})$ and
\[\langle \eta_1,\zeta_1\rangle_{H^{2}_{0}(\Gamma_s)}=\int_{\Gamma_s}(-A_{\alpha,\beta})^{1/2}\eta_1(-A_{\alpha,\beta})^{1/2}\zeta_1=\int_{\Gamma_s}(\beta\eta_{1,x}\zeta_{1,x} + \alpha\eta_{1,xx}\zeta_{1,xx}).\]
This scalar product on $H^{2}_{0}(\Omega_{0})$ is used to simplify calculations involving the operator $A_{\alpha,\beta}$. The unbounded operator relative to the beam $(A_{s},\mathcal{D}(A_{s}))$ in
\[H_{s}=H^{2}_{0}(\Gamma_s)\times L^{2}(\Gamma_s),\]
is defined by $\mathcal{D}(A_{s})=(H^{4}(\Gamma_s)\cap H^{2}_{0}(\Gamma_s))\times H^{2}_{0}(\Gamma_s)$ and $A_{s}=
\begin{pmatrix}
0&I\\
A_{\alpha,\beta}&\delta\Delta_s\\
\end{pmatrix}
$.
\begin{thm}\label{s5thm1}The operator $(\mathcal{A},\mathcal{D}(\mathcal{A}))$ is the infinitesimal generator of an analytic semigroup on $\textbf{H}$.
\end{thm}
\begin{proof}
The idea of the proof is to split the operator $\mathcal{A}$ in two parts. The principal part of $\mathcal{A}$ will be the infinitesimal generator of an analytic semigroup on $\textbf{H}$ and the rest will be a perturbation bounded with respect to the principal part. Set $K_s=(I+N_s)^{-1}-I$. The operator $\mathcal{A}$ can be written
\[\mathcal{A}=\mathcal{A}_1 + \mathcal{A}_2,\]
with 
\[\mathcal{A}_1=
\begin{pmatrix}
A&0&(-A)\Pi D_s\\
0&0&I\\
0&A_{\alpha,\beta}&\delta\Delta_s\\
\end{pmatrix},
\]
and 
\[\mathcal{A}_2=
\begin{pmatrix}
0&0&0\\
0&0&0\\
\nu(I+N_s)^{-1}N_{0}(\Delta(\cdot)\cdot\normal)&K_s A_{\alpha,\beta}&\delta K_s \Delta_s\\
\end{pmatrix}.
\]
According to \cite[Section 3.2, Theorem 2.1]{MR0512912}, the result is a direct consequence of  Theorem \ref{s5thm2} and \ref{s5thm3}.
\end{proof}
\begin{thm}\label{s5thm2}The operator $(\mathcal{A}_1,\mathcal{D}(\mathcal{A}_1)=\mathcal{D}(\mathcal{A}))$ is the infinitesimal generator of an analytic semigroup on $\textbf{H}$.
\end{thm}
\begin{proof}
The proof follows the techniques used in \cite[Theorem 3.5]{MR2745779}. The first part is to prove that the semigroup $\mathcal{A}_{1}$ is strongly continuous. This property, established using regularity results on the unsteady Stokes equations, is proved in the appendix (Lemma \ref{A1continuous}).

The next step is to estimate the resolvent of $\mathcal{A}_1$. Using a perturbation argument to ensure the existence of the resolvent, we have, at least for $\text{Re}(\lambda)>0$,
\[
(\lambda I-\mathcal{A}_1)^{-1}=
\begin{pmatrix}
(\lambda I-A)^{-1}& (0\,\,(\lambda I-A)^{-1}(-A)\Pi D_s)(\lambda I-A_s)^{-1}\\
0&(\lambda I-A_s)^{-1}\\
\end{pmatrix},
\]
where $(\lambda I-A_s)^{-1}$ is given by
\[
(\lambda I-A_s)^{-1}=
\begin{pmatrix}
\mathcal{V}^{-1}(\lambda I-\delta\Delta_s)&\mathcal{V}^{-1}\\
\mathcal{V}^{-1}A_{\alpha,\beta}&\lambda\mathcal{V}^{-1}\\
\end{pmatrix},
\]
with $\mathcal{V}=\lambda^{2}I-\lambda\delta\Delta_{s}-A_{\alpha,\beta}$. From \cite{MR971932} we know that there exists $a\in\mathbb{R}$ and $\frac{\pi}{2}<\theta_0<\pi$ such that for all $\lambda$ in $S_{a,\theta_0}=\{\lambda\in \mathbb{C}\mid \lambda\ne a,\, \vert\text{arg}(\lambda-a)\vert<\theta_0\}$
\[\norme{(\lambda I-A_s)^{-1}}_{\mathcal{L}(H_s)}\leq \frac{C_s}{\vert \lambda-a\vert}.\]
The Stokes operator $A$ is the infinitesimal generator of an analytic semigroup on $\textbf{V}^{0}_{n,\Gamma_{d}}(\Omega_{0})$ and the proof of Theorem \ref{s62thm1} gives the existence of $\frac{\pi}{2}<\theta_1<\pi$ such that for all $\lambda\in S_{0,\theta_1}$
\[\norme{(\lambda I-A)^{-1}}_{\mathcal{L}(\textbf{V}^{0}_{n,\Gamma_{d}}(\Omega_{0}))}\leq \frac{C_{A}}{\vert \lambda\vert}.
\]
Choose $a'>a$ and $\theta'=\min(\theta_0,\theta_{1})$. For all $\lambda\in S_{a',\theta'}$ and all $(\textbf{f},\Phi)\in \textbf{V}^{0}_{n,\Gamma_{d}}(\Omega_{0})\times H_s$ we have
\[
(\lambda I-\mathcal{A}_1)^{-1}
\begin{pmatrix}
\textbf{f}\\
\Phi\\
\end{pmatrix}\\
=\begin{pmatrix}
(\lambda I-A)^{-1}\textbf{f}+(\lambda I-A)^{-1}(-A)\Pi D_s ((\lambda I-A_s)^{-1}\Phi)_{2}\\
(\lambda I-A_s)^{-1}\Phi\\
\end{pmatrix}.
\]
Remark that $(\lambda I-A)^{-1}(-A)\Pi D_s=\Pi D_s -\lambda(\lambda I-A)^{-1}\Pi D_s$. Using the previous estimates for the resolvent of $A$ and $A_s$ and the continuity of the operator $\Pi D_s$ we obtain
\begin{align*}
&\norme{(\lambda I-\mathcal{A}_{1})^{-1}\begin{pmatrix}
\textbf{f}\\
\Phi\\
\end{pmatrix}}_{\textbf{V}^{0}_{n,\Gamma_{d}}(\Omega_{0})}\\
&{}\leq \frac{C_{A}}{\vert \lambda-a'\vert}\norme{\textbf{f}}_{\textbf{V}^{0}_{n,\Gamma_{d}}(\Omega_{0})} + \frac{C_{\Pi D_s}C_s}{\vert \lambda-a'\vert}\norme{\Phi}_{H_s} + \frac{C_{A}C_{\Pi D_s}C_{s}}{\vert\lambda-a'\vert}\norme{\Phi}_{H_s} +\frac{C_s}{\vert \lambda-a'\vert}\norme{\Phi}_{H_s}.\\
\end{align*}
Hence there exists a constant $C>0$ such that
\[\norme{(\lambda I-\mathcal{A}_{1})^{-1}}_{\mathcal{L}(\textbf{H})}\leq \frac{C}{\vert \lambda-a'\vert},\]
for all $\lambda\in S_{a',\theta'}$ and $\mathcal{A}_{1}$ is the infinitesimal generator of an analytic semigroup on $\textbf{H}$.
\end{proof}
\begin{thm}\label{s5thm3}The operator $(\mathcal{A}_2,\mathcal{D}(\mathcal{A}_2)=\mathcal{D}(\mathcal{A}))$ is $\mathcal{A}_1$-bounded with relative bound zero.
\end{thm}
\begin{proof}
We proceed as in \cite{MR2745779}. Split the operator $\mathcal{A}_2$ in three parts $\mathcal{A}_{2}=\mathcal{A}_{2,1}+\mathcal{A}_{2,2}+\mathcal{A}_{2,3}$ with
\[\mathcal{A}_{2,1}=
\begin{pmatrix}
0&0&0\\
0&0&0\\
\nu(I+N_s)^{-1}N_{0}(\Delta(\cdot)\cdot\normal)&0&0\\
\end{pmatrix},
\]
\[\mathcal{A}_{2,2}=
\begin{pmatrix}
0&0&0\\
0&0&0\\
0&K_s A_{\alpha,\beta}&0\\
\end{pmatrix},\,\,
\mathcal{A}_{2,3}=
\begin{pmatrix}
0&0&0\\
0&0&0\\
0&0&\delta K_s \Delta_s\\
\end{pmatrix}.
\]
The following lemma is an adaptation of \cite[Proposition 3.3]{MR2745779}.
\begin{lem}\label{s5lem1}
The norm
\[(\Pi\ug,\eta_1,\eta_2)\mapsto \norme{(\Pi\ug,\eta_1,\eta_2)}_{\textbf{H}}+\norme{A\Pi\ug + (-A)\Pi D_s\eta_2}_{\textbf{V}^{0}_{n,\Gamma_{d}}(\Omega_{0})}+\norme{A_s(\eta_1,\eta_2)}_{H_s},\]
is a norm on $\mathcal{D}(\mathcal{A})$ which is equivalent to the norm
\[(\Pi\ug,\eta_1,\eta_2)\mapsto\norme{\Pi\ug}_{\textbf{V}^{2}_{n,\Gamma_{d}}(\Omega_{0})}+\norme{\eta_1}_{H^{4}(\Gamma_s)}+\norme{\eta_2}_{H^{2}_{0}(\Gamma_s)}.\]
\end{lem}
For $\mathcal{A}_{2,2}$ and $\mathcal{A}_{2,3}$ we can use \cite[Lemma 3.9]{MR2745779} and \cite[Lemma 3.10]{MR2745779} to prove that there exists $0<\theta_1<1$ and $0<\theta_2<1$ such that $\mathcal{A}_{2,2}$ (respectively $\mathcal{A}_{2,3}$) is bounded from $\mathcal{D}((-\mathcal{A}_{1})^{\theta_1})$ (respectively from $\mathcal{D}((-\mathcal{A}_{1})^{\theta_2})$) into $\textbf{H}$. Hence, according to \cite[Section 3.2, Corollary 2.4]{MR0512912}, the operators $\mathcal{A}_{2,2}$ and $\mathcal{A}_{2,3}$ are $\mathcal{A}_{1}$-bounded with relative bound zero. It remains to prove that $\mathcal{A}_{2,1}$ is $\mathcal{A}_{1}$-bounded with relative bound zero.
\begin{lem}\label{s5lem2}For all $\varepsilon>0$ there exists a constant $C_{\varepsilon}$ such that
\begin{equation}\label{s5est1}
\norme{N_{0}(\Delta\ug\cdot\normal)}_{L^{2}(\Gamma_s)}\leq \varepsilon\norme{\ug}_{\textbf{V}^{2}_{n,\Gamma_d}(\Omega_{0})}+C_{\varepsilon}\norme{\ug}_{\textbf{V}^{0}_{n,\Gamma_d}(\Omega_{0})},
\end{equation}
for all $\ug\in\textbf{V}^{2}_{n,\Gamma_{d}}(\Omega_{0})$.
\end{lem}
\begin{proof}
Using the transposition method, a density argument (as in Theorem \ref{s61lem1} and Theorem \ref{s61thm6}) and interpolation, the operator $N_{0}$ can be defined as a continuous operator from $H^{-1}(\Gamma_s)$ into $L^{2}(\Gamma_s)$. We prove the lemma by contradiction. Assume that there exists a sequence $\ug_{k}\in \textbf{V}^{2}_{n,\Gamma_{d}}(\Omega_{0})$ such that 
\[\norme{N_{0}(\Delta\ug_k\cdot\normal)}_{L^{2}(\Gamma_s)}=1,\,\,\norme{\ug_{k}}_{\textbf{V}^{0}_{n,\Gamma_{d}}(\Omega_{0})}\longrightarrow 0,\,\,\norme{\ug_k}_{\textbf{V}^{2}_{n,\Gamma_{d}}(\Omega_{0})}\leq M,\]
with $M>0$ a fixed constant. By reflexivity of the space $\textbf{V}^{2}_{n,\Gamma_{d}}(\Omega_{0})$, up to a subsequence, there exists $\ug\in \textbf{V}^{2}_{n,\Gamma_{d}}(\Omega_{0})$ such that $\ug_k\rightharpoonup \ug$ in $\textbf{V}^{2}_{n,\Gamma_{d}}(\Omega_{0})$. Since $\norme{\ug_k}_{\textbf{V}^{0}_{n,\Gamma_{d}}(\Omega_{0})}\longrightarrow 0$, we obtain $\ug=0$. Then $\Delta\ug_k\cdot\normal\rightharpoonup 0$ in $H^{-1/2}(\Gamma_s)$ and the compact embedding of $H^{-1/2}(\Gamma_s)$ into $H^{-1}(\Gamma_s)$ ensures that $\Delta\ug_k\cdot\normal\longrightarrow 0$ in $H^{-1}(\Gamma_s)$. Finally the continuity of $N_{0}$ implies that $N_{0}(\Delta\ug_k\cdot\normal)\longrightarrow 0$ in $L^{2}(\Gamma_s)$ which contradicts $\norme{N_{0}(\Delta\ug_k\cdot\normal)}_{L^{2}(\Gamma_s)}=1$.
\end{proof}
We come back to the proof that $\mathcal{A}_{2,1}$ is $\mathcal{A}_{1}$-bounded. Using the estimate (\ref{s5est1}) and the norm equivalence of Lemma \ref{s5lem1} it follows that for all $(\Pi\ug,\eta_1,\eta_2)\in \mathcal{D}(\mathcal{A}_{1})$
\begin{align*}
\norme{\mathcal{A}_{2,1}
\begin{pmatrix}
\Pi\ug\\
\eta_1\\
\eta_2\\
\end{pmatrix}}_{\textbf{H}}
&{}\leq \varepsilon\norme{\Pi\ug}_{\textbf{V}^{2}_{n,\Gamma_{d}}(\Omega_{0})}+C_{\varepsilon}\norme{\Pi\ug}_{\textbf{V}^{0}_{n,\Gamma_{d}}(\Omega_{0})}\\
&{}\leq \varepsilon(\norme{\Pi\ug}_{\textbf{V}^{2}_{n,\Gamma_{d}}(\Omega_{0})}+\norme{\eta_1}_{H^{4}(\Gamma_s)}+\norme{\eta_2}_{H^{2}_{0}(\Gamma_s)})+C_{\varepsilon}\norme{\Pi\ug}_{\textbf{V}^{0}_{n,\Gamma_{d}}(\Omega_{0})}\\
&{}\leq C_{1}\varepsilon\norme{\mathcal{A}_{1}\begin{pmatrix}
\Pi\ug\\
\eta_1\\
\eta_2\\
\end{pmatrix}}_{\textbf{H}}+C_{2,\varepsilon}\norme{\begin{pmatrix}
\Pi\ug\\
\eta_1\\
\eta_2\\
\end{pmatrix}}_{\textbf{H}}.\\
\end{align*}
This concludes the proof that $\mathcal{A}_{2}$ is $\mathcal{A}_{1}$-bounded with relative bound zero.
\end{proof}
\subsection{Regularity results}
We have seen that the system (\ref{s33e1}) can be rewritten
\begin{equation}\label{s34e1}
\begin{cases}
\begin{aligned}
&\frac{d}{dt}\begin{pmatrix}
\Pi\ug\\
\eta_{1}\\
\eta_{2}\\
\end{pmatrix}
=\mathcal{A}
\begin{pmatrix}
\Pi\ug\\
\eta_{1}\\
\eta_{2}\\
\end{pmatrix}
+\textbf{F},\,\,
\begin{pmatrix}
\Pi\ug(0)\\
\eta_{1}(0)\\
\eta_{2}(0)\\
\end{pmatrix}
=\begin{pmatrix}
\Pi\ug^{0}\\
\eta^{0}_{1}\\
\eta^{0}_{2}\\
\end{pmatrix},\\
&(I-\Pi)\ug=(I-\Pi)D_s(\eta_{2}).\\
\end{aligned}
\end{cases}
\end{equation}
We remark that there is no condition on $(I-\Pi)\ug^{0}$. As in \cite{MR2745779}, in order to satisfy the 
equality $(I-\Pi)\ug=(I-\Pi)D_s(\eta_{2})$ at time $t=0$, we introduce a subspace of initial conditions 
belonging to $\textbf{V}^{1}(\Omega_{0})\times H_s$ and satisfying a compatibility condition 
\[\textbf{H}_{\text{cc}}=\{(\ug^{0},\eta_1^0,\eta_2^0)\in \textbf{V}^{1}(\Omega_{0})\times H_s\mid (I-\Pi)\ug^{0}=(I-\Pi)D_s(\eta_2^0)\}.\] 
To obtain maximal regularity results, introduce the space $[\mathcal{D}(\mathcal{A}),\textbf{H}]_{1/2}$ given by
\begin{align*}
[\mathcal{D}(\mathcal{A}),\textbf{H}]_{1/2}=\{(\Pi\ug, {}&\eta_1,\eta_2)\in \textbf{V}^{1}_{n,\Gamma_{d}}(\Omega_{0})\times (H^{3}(\Gamma_s)\cap H^{2}_{0}(\Gamma_s))\times H^{1}_{0}(\Gamma_s)\mid\\
& \Pi\ug -\Pi D_s(\eta_2)\in V\}.
\end{align*}
It is equipped with the norm
\[(\Pi\ug,\eta_1,\eta_2)\longmapsto \left(\norme{\Pi \ug}^{2}_{\textbf{H}^{1}(\Omega_{0})}+\norme{\eta_1}^{2}_{H^{3}(\Gamma_s)}+\norme{\eta_2}^{2}_{H^{1}(\Gamma_s)}\right)^{1/2}.\]
If the initial condition $(\Pi \textbf{u}^{0},\eta_1^0,\eta_2^0)$ belongs to $[\mathcal{D}(\mathcal{A}),\textbf{H}]_{1/2}$, and if the compatibility condition $(I-\Pi)\ug^{0}=(I-\Pi)D_s(\eta_2^0)$ is satisfied, then $(\ug^{0},\eta_1^0,\eta_2^0)$ belongs to 
$$
\mathcal{X}(\Omega_{0})=\{(\ug^{0},\eta_1^0,\eta_2^0)\in  \textbf{H}_{\text{cc}} \mid (\Pi \textbf{u}^{0},\eta_1^0,\eta_2^0)\in [\mathcal{D}(\mathcal{A}),\textbf{H}]_{1/2}\}. 
$$
The space $\mathcal{X}(\Omega_{0})$ is equipped with the norm 
$$
\xg^0=(\textbf{u}^{0},\eta_1^0,\eta_2^0) \longmapsto \norme{\xg^0}_{\mathcal{X}(\Omega_{0})} =\left(\norme{\Pi \textbf{u}^{0}}^{2}_{\textbf{H}^{1}(\Omega_{0})}+\norme{\eta_1^0}^{2}_{H^{3}(\Gamma_s)}+\norme{\eta_2^0}^{2}_{H^{1}(\Gamma_s)}\right)^{1/2}. 
$$
We notice that the above mapping is indeed a norm since $(I-\Pi)\ug^{0}=(I-\Pi)D_s(\eta_2^0)$ if 
$\xg^0=(\textbf{u}^{0},\eta_1^0,\eta_2^0)\in \mathcal{X}(\Omega_{0})$. 
Defining $W_{T}$ by
\[W_{T}=\textbf{L}^{2}(Q_{T})\times L^{2}(0,T;H^{1/2}(\Gamma_{i,o}))\times L^{2}(0,T;L^{2}(\Gamma_s)),\]
we obtain the main theorem of this section.
\begin{thm}\label{s5thm5}For all $(\ug^{0},\eta^{0}_{1},\eta^{0}_{2})$ in $\mathcal{X}(\Omega_{0})$ and $(\textbf{f},\Theta,h)$ in $W_{T}$, system (\ref{s32e3}) admits a unique solution $(\ug,p,\eta)\in \textbf{H}^{2,1}(Q_{T})\times L^{2}(0,T;H^{1}(\Omega_{0}))\times H^{4,2}(\Sigma^{s}_{T})$. This solution satisfies
\begin{equation}\label{s5est3}
\begin{aligned}
&\norme{\ug}_{\textbf{H}^{2,1}(Q_T)}+\norme{\eta}_{H^{4,2}(\Sigma^{s}_{T})}+\norme{p}_{L^{2}(0,T;H^{1}(\Omega_{0}))}\\
&{}\leq C_{L}(\norme{(\ug^{0},\eta^{0}_{1},\eta^{0}_{2})}_{\mathcal{X}(\Omega_{0})}+\norme{(\textbf{f},\Theta,h)}_{W_{T}}).
\end{aligned}
\end{equation}
\end{thm}
\begin{proof}
According to \cite[Part II, Section 1.3, Theorem 3.1]{MR2273323} there exists a unique solution $(\ug,\eta_{1},\eta_{2})$ to (\ref{s34e1}) and the following estimate holds
\begin{align*}
&\norme{\Pi\ug}_{\textbf{H}^{2,1}(Q_T)}+\norme{(\eta_{1},\eta_{2})}_{L^{2}(0,T;\mathcal{D}(A_{s}))\cap H^{1}(0,T;H_{s})}\\
&\qquad\leq C(\norme{(\Pi\ug^{0},\eta^{0}_{1},\eta^{0}_{2})}_{[\mathcal{D}(\mathcal{A}),\textbf{H}]_{1/2}}+\norme{\textbf{F}}_{L^{2}(0,T;\textbf{H})})\\
&\norme{(I-\Pi)\ug}_{L^{2}(0,T;\textbf{H}^{2}(\Omega_{0}))}+\norme{(I-\Pi)\ug}_{H^{1}(0,T;\textbf{H}^{1/2}(\Omega_{0}))}\\
&\qquad\leq C\norme{(\eta_{1},\eta_{2})}_{L^{2}(0,T;\mathcal{D}(A_{s}))\cap H^{1}(0,T;H_{s})},
\end{align*}
with the estimate on $(I-\Pi)\ug$ coming from the properties of the operator $D_s$ and the identity $(I-\Pi)\ug=(I-\Pi)\eta_{2}$. Estimate (\ref{s5est3}) follows by writing $\textbf{F}$ explicitly.
\end{proof}
\begin{rmq}Let $T_{0}$ be a fixed time with $T<T_{0}$. The constant $C_{L}$ in the previous statement can be chosen independent of $T$ for all $T<T_{0}$. If we extend all the nonhomogeneous terms on $[T,T_{0}]$ by $0$ (still denoted by $(\textbf{f},\Theta,h)$) the previous result implies that there exists a unique solution $(\overset{\circ}{\ug},\overset{\circ}{p},\overset{\circ}{\eta})$ of (\ref{s32e3}) and the following estimate holds
\[
\begin{aligned}
&\norme{\overset{\circ}{\ug}}_{\textbf{H}^{2,1}(Q_{T_{0}})}+\norme{\overset{\circ}{\eta}}_{H^{4,2}(\Sigma^{s}_{T_{0}})}+\norme{\overset{\circ}{p}}_{L^{2}(0,T_{0};H^{1}(\Omega_{0}))}\\
&{}\leq C_{L}(T_{0})(\norme{(\ug^{0},\eta^{0}_{1},\eta^{0}_{2})}_{\mathcal{X}(\Omega_{0})}+\norme{(\textbf{f},\Theta,h)}_{W_{T_{0}}}).
\end{aligned}
\]
The uniqueness yields $(\ug,p,\eta)=(\overset{\circ}{\ug},\overset{\circ}{p},\overset{\circ}{\eta})$ on $[0,T]$ and the constant $C_{L}$ in Theorem \ref{s5thm5} can be taken as $C_{L}=C_{L}(T_{0})$.
\end{rmq}
\section{Nonlinear coupled system}
Throughout this section, excepted for Theorem \ref{s42thm3} which is stated in a rectangular domain, $\Omega_{0}$ is given by \eqref{domain-ref} for any fixed $\eta^{0}\in H^{3}(\Gamma_{s})\cap H^{2}_{0}(\Gamma_{s})$ such that $1+\eta^{0}(x)>0$ for all $x\in (0,L)$. We prove the existence of strong solutions for the complete nonlinear system (\ref{s32e2}). Let $T_{0}>0$ be a given time, fixed for this section. Let $\tilde{\mathcal{X}}(\Omega_{0})$ be the affine subspace of $\mathcal{X}(\Omega_{0})$ defined by
\[
\tilde{\mathcal{X}}(\Omega_{0}) = \{ (\ug^{0},\eta_1^0,\eta_2^0)\in \mathcal{X}(\Omega_{0})\mid \eta_1^0=\eta^0 \},
\]
that is, the space where the initial data of the beam $\eta^{1}_{0}$ and the geometric $\eta^{0}$ are equal. For $T>0$, set
\begin{align*}
\mathcal{Y}_{T}=\{&(\overline{\ug},p,\eta)\in \textbf{H}^{2,1}(Q_{T})\times L^{2}(0,T;H^{1}(\Omega_{0}))\times H^{4,2}(\Sigma^{s}_{T})\mid\\
& \overline{\ug}=0\text{ on }\Sigma^{b}_{T},\,\overline{\ug}=\eta_{t}\textbf{e}_2\text{ on }\Sigma^{0}_{T},\,\overline{u}_2=0\text{ on }\Sigma^{i,o}_{T},\ (\overline{\ug}(0),\eta(0),\eta_{t}(0))\in \tilde{\mathcal{X}}(\Omega_{0})\}.
\end{align*}
The usual norm on $\textbf{H}^{2,1}(Q_{T})\times L^{2}(0,T;H^{1}(\Omega_{0}))\times H^{4,2}(\Sigma^{s}_{T})$ is denoted by $\norme{\cdot}_{\mathcal{Y}_{T}}$. 

\subsection{Estimates} For every $\xg^{0}=(\ug^{0},\eta_1^0,\eta_2^0)\in \tilde{\mathcal{X}}(\Omega_{0})$, $R>0$, $\mu>0$ and $T>0$, define the ball
$$
\begin{array}{l}
\mathcal{B}(\xg^{0},R,\mu,T)
=\{(\overline{\ug},p,\eta)\in\mathcal{Y}_{T}\mid (\overline{\ug}(0),\eta(0),\eta_{t}(0))=\xg^{0},\  
\norme{(\overline{\ug},p,\eta)}_{\mathcal{Y}_{T}}\leq R,
\vspace{2.mm}\\
\hspace{5.4cm}\norme{(1+\eta)^{-1}}_{L^{\infty}(\Sigma^{s}_{T})}\leq 2\mu\} .
\end{array}
$$
For a given $\xg^{0}=(\ug^{0},\eta_1^0,\eta_2^0)\in \tilde{\mathcal{X}}(\Omega_{0})$, Theorem \ref{s5thm5} ensures the existence of $R>0$ and $\mu>0$ such that $\mathcal{B}(\xg^{0},R,\mu,T)$ is non empty for $T>0$ small enough (such a triple $(R,\mu,T)$ is explicitly chosen in the beginning of the proof of Theorem \ref{s42thm1}).

Throughout this section, $C(T_{0},R,\mu,\|\xg^{0}\|_{\mathcal{X}(\Omega_{0})})$ denotes a constant, depending on $T_{0}$, $R>0$, $\mu>0$, $\|\xg^{0}\|_{\mathcal{X}(\Omega_{0})} $ which may vary from one statement to another, but which is independent of $T$. 

The following lemmas are used to estimate the nonlinear terms (see Theorem \ref{s41thm1}).

\begin{lem}\label{s41lemma1} There exists a constant $C_{0}$ depending on $T_{0}$ such that, for all $0<T< T_0$ and all $u\in H^{2,1}(Q_{T})$ satisfying $u(0)=0$, the following estimate holds
\[\norme{u}_{L^{\infty}(0,T;H^{1}(\Omega_{0}))}+\norme{u}_{L^{4}(0,T;L^{\infty}(\Omega_{0}))}\leq C_{0}\norme{u}_{H^{2,1}(Q_{T})}.\]
Moreover for all $v\in H^{4,2}(\Sigma^{s}_{T})$ satisfying $v(0)=0$ the following estimate holds
$$
\norme{v}_{L^{\infty}(0,T;H^{3}(\Gamma_{s}))}\leq C_{0}\norme{v}_{H^{4,2}(\Sigma^{s}_{T})}. 
$$
If in addition $v_{t}(0)=0$, then
\[\norme{v_{t}}_{L^{\infty}(0,T;H^{1}(\Gamma_{s}))}+\norme{v_{t}}_{L^{2}(0,T;H^{2}(\Gamma_{s}))}\leq C_{0}\norme{v}_{H^{4,2}(\Sigma^{s}_{T})}.\]
\end{lem}
\begin{proof}
These estimates come from interpolation results (see \cite{MR0350178}). The only thing to prove is that the continuity constants can be made independent of $T$. Let $\overline{u}$ be the function defined by $\overline{u}=0$ on $[T-T_{0},0]$ and $\overline{u}=u$ on $[0,T]$. As $u(0)=0$, the function $\overline{u}$ is still in $H^{2,1}(Q_{T})$ and, using interpolation estimates, we have
\[\norme{\overline{u}}_{L^{\infty}(T-T_{0},T;H^{1}(\Omega_{0}))}\leq C(T_{0})\norme{\overline{u}}_{L^{2}(T-T_{0},T;H^{2}(\Omega_{0}))\cap H^{1}(T-T_{0},T;L^{2}(\Omega_{0}))}.\]
This implies that $\norme{u}_{L^{\infty}(0,T;H^{1}(\Omega_{0}))}\leq C_{0}\norme{u}_{H^{2,1}(Q_{T})}$ with $C_{0}=C(T_{0})$. The other estimates follow from the same argument.
\end{proof}

\begin{lem}\label{s41lemma2}Let $\xg^{0}$ belong to ${\tilde{\mathcal{X}}(\Omega_{0})}$, $R>0$, and $\mu>0$. There exists a constant $C(T_{0},R,\mu,\|\xg^{0}\|_{\mathcal{X}(\Omega_{0})})>0$ such that, for all $0<T< T_{0}$ and all $(\overline{\ug},p,\eta)$ in $\mathcal{B}(\xg^{0},R,\mu,T)$, the following estimates hold
\begin{align*}
&\norme{\overline{\ug}}_{L^{\infty}(0,T;\textbf{H}^{1}(\Omega_{0}))}+\norme{\overline{\ug}}_{L^{4}(0,T;\textbf{L}^{\infty}(\Omega_{0}))}+\norme{\eta_{t}}_{L^{\infty}(0,T;H^{1}(\Gamma_{s}))}\\
&+\norme{\eta_{t}}_{L^{2}(0,T;H^{2}(\Gamma_{s}))}+\norme{\eta_{xx}}_{L^{\infty}(0,T;H^{1}(\Gamma_{s}))}\leq C(T_{0},R,\mu,\|\xg^{0}\|_{\mathcal{X}(\Omega_{0})}).
\end{align*}
\end{lem}
\begin{proof}
Let $(\overset{\circ}{\ug},\overset{\circ}{\eta},\overset{\circ}{p})$ be the solution to (\ref{s32e3}) on the time interval $[0,T_{0}]$ with right-hand side $0$ and the initial condition $(\ug^{0},\eta_1^0,\eta_2^0)$. We have
\[
\norme{\overline{\ug}}_{L^{\infty}(0,T;\textbf{H}^{1}(\Omega_{0}))}\leq \norme{\overline{\ug}-\overset{\circ}{\ug}}_{L^{\infty}(0,T;\textbf{H}^{1}(\Omega_{0}))} +\norme{\overset{\circ}{\ug}}_{L^{\infty}(0,T;\textbf{H}^{1}(\Omega_{0}))},
\]
then using that $\overset{\circ}{\ug}(0)=\overline{\ug}(0)$ and Lemma \ref{s41lemma1} the following estimate holds with the constant $C_{0}=C_{0}(T_{0})$
\begin{align*}
\norme{\overline{\ug}-\overset{\circ}{\ug}}_{L^{\infty}(0,T;\textbf{H}^{1}(\Omega_{0}))}&{}\leq C_{0}\norme{\overline{\ug}-\overset{\circ}{\ug}}_{\textbf{H}^{2,1}(Q_{T})}\\
&{}\leq C_{0}\norme{\overline{\ug}}_{\textbf{H}^{2,1}(Q_{T})}+C_{0}\norme{\overset{\circ}{\ug}}_{\textbf{H}^{2,1}(Q_{T})}\\
&{}\leq C_{0}R+C_{0}\norme{\overset{\circ}{\ug}}_{\textbf{H}^{2,1}(Q_{T_{0}})}.
\end{align*}
The second part is estimated as follows
\[\norme{\overset{\circ}{\ug}}_{L^{\infty}(0,T;\textbf{H}^{1}(\Omega_{0}))}\leq \norme{\overset{\circ}{\ug}}_{L^{\infty}(0,T_0;\textbf{H}^{1}(\Omega_{0}))} \leq C_{0}'\norme{\overset{\circ}{\ug}}_{\textbf{H}^{2,1}(Q_{T_{0}})}.\]
Finally estimate (\ref{s5est3}) on $\overset{\circ}{\ug}$ implies
\[
\norme{\overline{\ug}}_{L^{\infty}(0,T;\textbf{H}^{1}(\Omega_{0}))}\leq C_{0}R+(C_{0}+C_{0}')C_{L}\norme{(\ug^{0},\eta_1^0,\eta_2^0)}_{\mathcal{X}(\Omega_{0})}.
\]
The estimates on $\eta_{t}$ and $\eta_{xx}$ follow similarly.
\end{proof}
\begin{lem}\label{s41lemma3}
Set $\mu_{0}=\norme{(1+\eta^{0})^{-1}}_{L^{\infty}(\Gamma_{s})}$. For $\eta\in H^{4,2}(\Sigma^{s}_{T})$ such that $\eta(0)=\eta^{0}$ the function $\widetilde{\eta}=\frac{\eta-\eta^{0}}{1+\eta^{0}}$ satisfies the following estimates
\begin{equation}\label{s41lemma3est1}
\norme{\widetilde{\eta}}_{L^{\infty}(\Sigma^{s}_{T})}\leq \mu_{0}T^{1/2}\norme{\eta_{t}}_{L^{2}(0,T;L^{\infty}(\Gamma_{s}))},
\end{equation}
\begin{equation}\label{s41lemma3est2}
\norme{\widetilde{\eta}_{x}}_{L^{\infty}(\Sigma^{s}_{T})}\leq \mu_{0}T^{1/2}\norme{\eta_{tx}}_{L^{2}(0,T,L^{\infty}(\Gamma_{s}))} + \mu_{0}^{2}\norme{\eta^{0}_{x}}_{L^{\infty}(\Gamma_{s})}T^{1/2}\norme{\eta_{t}}_{L^{2}(0,T;L^{\infty}(\Gamma_{s}))}.
\end{equation}
\end{lem}
\begin{proof}
The estimates come from the fundamental theorem of calculus and Cauchy-Schwarz inequality.
\end{proof}
\begin{lem}\label{s41lemma4}Let $u$ be in $L^{\infty}(0,T;H^{1}(\Omega_{0}))$ and $v$ be in $H^{2,1}(Q_{T})$. The following estimate holds
\begin{equation}\label{s41lemma4e1}
\norme{u\partial_{i}v}_{L^{2}(Q_{T})}\leq CT^{1/4}\norme{u}_{L^{\infty}(0,T;H^{1}(\Omega_{0}))}\norme{v}_{L^{\infty}(0,T;H^{1}(\Omega_{0}))}^{1/2}\norme{v}_{L^{2}(0,T,H^{2}(\Omega_{0}))}^{1/2},
\end{equation}
with $C$ independent of $T$, and $i=1,2$.
\end{lem}
\begin{proof}
We have
\[\int_{\Omega_{0}}\vert u\vert^{2}\vert\partial_{1}v\vert^{2}\leq \left(\int_{\Omega_{0}}\vert u\vert^{6}\right)^{1/3}\left(\int_{\Omega_{0}}\vert\partial_{i} v\vert^{3}\right)^{2/3},
\]
and, using Lebesgue interpolation, $\norme{\partial_{i}v}_{L^{3}(\Omega_{0})}\leq \norme{\partial_{i}v}_{L^{2}(\Omega_{0})}^{1/2}\norme{\partial_{i}v}_{L^{6}(\Omega_{0})}^{1/2}$. Sobolev embeddings then yield
\[\norme{\partial_{1}v }_{L^{4}(0,T;L^{3}(\Omega_{0}))}\leq C\norme{v}_{L^{\infty}(0,T;H^{1}(\Omega_{0}))}^{1/2}\norme{v}_{L^{2}(0,T;H^{2}(\Omega_{0}))}^{1/2}.
\]
Then we use the estimate $\norme{c}_{L^{2}(0,T)}\leq T^{1/4}\norme{c}_{L^{4}(0,T)}$ to obtain
\begin{align*}
\norme{u\partial_{i}v}_{L^{2}(0,T;L^{2}(\Omega_{0}))}&{}\leq \norme{u}_{L^{\infty}(0,T;H^{1}(\Omega_{0}))}\norme{\partial_{i}v}_{L^{2}(0,T;L^{3}(\Omega_{0}))}\\
&{}\leq T^{1/4}\norme{u}_{L^{\infty}(0,T;H^{1}(\Omega_{0}))}\norme{\partial_{i}v}_{L^{4}(0,T;L^{3}(\Omega_{0}))}\\
&{}\leq CT^{1/4}\norme{u}_{L^{\infty}(0,T;H^{1}(\Omega_{0}))}\norme{v}_{L^{\infty}(0,T;H^{1}(\Omega_{0}))}^{1/2}\norme{v}_{L^{2}(0,T;H^{2}(\Omega_{0}))}^{1/2}.
\end{align*}
\end{proof}
\begin{lem}\label{s41lemma5}Let $\xg^{0}$ belongs to ${\tilde{\mathcal{X}}(\Omega_{0})}$, $R>0$, and $\mu>0$. For all $0<T< T_{0}$ and $(\overline{\ug},p,\eta)$ in $\mathcal{B}(\xg^{0},R,\mu,T)$, the function $M(\overline{\ug},\eta)$ belongs to $\textbf{H}^{2,1}(Q_{T})$ and the following estimate holds
\begin{equation}\label{s41lemma5est1}
\norme{M(\overline{\ug},\eta)}_{\textbf{H}^{2,1}(Q_{T})}\leq C(T_{0},R,\mu,\|\xg^{0}\|_{\mathcal{X}(\Omega_{0})}).
\end{equation}
Furthermore, for all $(\overline{\ug}_{1},p_{1},\eta_{1})$ and $(\overline{\ug}_{2},p_{2},\eta_{2})$ belonging to $\mathcal{B}(\xg^{0},R,\mu,T)$, the following Lipschitz estimate holds
\begin{equation}\label{s41lemmaest2}\norme{M(\overline{\ug}_{1},\eta_{1})-M(\overline{\ug}_{2},\eta_{2})}_{\textbf{H}^{2,1}(Q_{T})}\leq C(T_{0},R,\mu,\|\xg^{0}\|_{\mathcal{X}(\Omega_{0})})\norme{(\overline{\ug}_{1},p_{1},\eta_{1})-(\overline{\ug}_{2},p_{2},\eta_{2})}_{\mathcal{Y}_{T}}.
\end{equation}
\end{lem}
\begin{proof}Through what follows we use the following basic estimate
\[\norme{(1+\widetilde{\eta})^{-1}}_{L^{\infty}(\Sigma^{s}_{T})}=\norme{\frac{1+\eta^{0}}{1+\eta}}_{L^{\infty}(\Sigma^{s}_{T})}\leq \mu\norme{1+\eta^{0}}_{L^{\infty}(\Gamma_{s})}.\]
Most of the estimates of $M(\overline{\ug},\eta)=\left(\displaystyle\frac{\overline{u}_{1}}{1+\widetilde{\eta}},\displaystyle\frac{z\widetilde{\eta}_{x}}{1+\widetilde{\eta}}\overline{u}_{1}+\overline{u}_{2}\right)^{T}$ and its derivatives are explicit $L^{\infty}\times L^{2}$ estimates using the previous lemmas and the regularity of $\widetilde{\eta}$. To estimate the second spatial derivative of $\frac{\overline{u}_{1}}{1+\widetilde{\eta}}$ we compute $\widetilde{\eta}_{xx}$:
\[\widetilde{\eta}_{xx}=\left(\frac{\eta_{xx}-\eta^{0}_{xx}}{1+\eta^{0}}\right)-\left(\frac{2(\eta_{x}-\eta^{0}_{x})\eta^{0}_{x}-(\eta-\eta^{0})\eta^{0}_{xx}}{(1+\eta^{0})^{2}}\right)+\frac{2(\eta^{0}_{x})^{2}(\eta-\eta^{0})}{(1+\eta^{0})^{4}}.
\]
Then
\begin{align*}
\norme{\widetilde{\eta}_{xx}}_{L^{\infty}(0,T;H^{1}(\Gamma_{0}))}&{}\leq C_{1}(\eta^{0})+C_{2}(\eta^{0})\norme{\eta_{xx}}_{L^{\infty}(0,T;H^{1}(\Gamma_{s}))}\\
&+C_{3}(\eta^{0})\norme{\eta_{x}}_{L^{\infty}(0,T;H^{1}(\Gamma_{s}))}+C_{4}(\eta^{0})\norme{\eta}_{L^{\infty}(0,T;H^{1}(\Gamma_{s}))}\\
&{}\leq C(T_{0},R,\mu,\|\xg^{0}\|_{\mathcal{X}(\Omega_{0})}).
\end{align*}
This estimate is more precise that the one needed here (it implies an estimate on $\norme{\widetilde{\eta}_{xx}}_{L^{\infty}(\Sigma^{s}_{T})}$ using spatial Sobolev embeddings); we stated it because it is used in the estimates of Theorem \ref{s41thm1}. For the time derivative we have
\[\norme{\frac{-\widetilde{\eta}_{t}}{(1+\widetilde{\eta})^{2}}}_{L^{\infty}(\Sigma^{s}_{T})}\leq \mu^{2}\norme{1+\eta^{0}}_{L^{\infty}(\Sigma^{s}_{T})}^{2}\mu \norme{\eta_{t}}_{L^{\infty}(\Sigma^{s}_{T})}\leq C(T_{0},R,\mu,\|\xg^{0}\|_{\mathcal{X}(\Omega_{0})}).\]
It follows from these estimates that
\[\norme{\frac{\overline{u}_{1}}{1+\widetilde{\eta}}}_{H^{2,1}(Q_{T})}\leq C(T_{0},R,\mu,\|\xg^{0}\|_{\mathcal{X}(\Omega_{0})}).
\]
The second component of $M(\overline{\ug},\eta)$ and its derivatives estimated similarly, except for the terms
\[\norme{\frac{z\widetilde{\eta}_{xxx}\overline{u}_{1}}{1+\widetilde{\eta}}}_{L^{2}(Q_{T})}\leq C(\mu,\eta^{0})\norme{\widetilde{\eta}_{xxx}\overline{u}_{1}}_{L^{2}(Q_{T})}.
\]
The term $\widetilde{\eta}_{xxx}$ is only in $L^{2}(\Gamma_{s})$ and we cannot use Lemma \ref{s41lemma4}. Let us write $\widetilde{\eta}=ND$ with $N=\eta-\eta^{0}$ and $D=(1+\eta^{0})^{-1}$. We have
\[\widetilde{\eta}_{xxx}=N_{xxx}D+3N_{xx}D_{x}+3N_{x}D_{xx}+ND_{xxx}.\]
When multiplied by $\displaystyle\frac{z\overline{u}_{1}}{1+\widetilde{\eta}}$, the terms involving up to two derivatives can be estimated directly. For
\[\frac{z(\eta_{xxx}-\eta^{0}_{xxx})\overline{u}_{1}}{(1+\eta^{0})(1+\widetilde{\eta})},\]
we have
\begin{align*}
\norme{\frac{z(\eta_{xxx}-\eta^{0}_{xxx})\overline{u}_{1}}{(1+\eta^{0})(1+\widetilde{\eta})}}_{L^{2}(Q_{T})}^{2}&{}\leq C(\mu,\eta^{0})\int_{0}^{T}\norme{(\eta_{xxx}-\eta^{0}_{xxx})(\cdot,t)}_{L^{2}(\Gamma_{s})}^{2}\norme{\overline{\ug}(\cdot,t)}_{\textbf{L}^{\infty}(\Omega_{0})}^{2}dt\\
&{}\leq C(\mu,\eta^{0})\norme{\eta_{xxx}-\eta^{0}_{xxx}}_{L^{4}(0,T;L^{2}(\Gamma_{s}))}^{2}\norme{\overline{\ug}}_{L^{4}(0,T;\textbf{L}^{\infty}(\Omega_{0}))}^{2}\\
&{}\leq C(\mu,\eta^{0})T^{1/4}\norme{\eta_{xxx}-\eta^{0}_{xxx}}_{L^{8}(0,T;L^{2}(\Gamma_{s}))}^{2}\norme{\overline{\ug}}_{L^{4}(0,T;\textbf{L}^{\infty}(\Omega_{0}))}^{2},
\end{align*}
and
\begin{align*}
\norme{\eta_{xxx}-\eta^{0}_{xxx}}_{L^{8}(0,T;L^{2}(\Gamma_{s}))}&{}\leq \norme{\eta_{xxx}}_{L^{8}(0,T;L^{2}(\Gamma_{s}))}+T^{1/8}\norme{\eta^{0}}_{H^{3}(\Gamma_{s})}\\
&{}\leq C(T_{0},R,\mu,\|\xg^{0}\|_{\mathcal{X}(\Omega_{0})})+T_{0}^{1/8}\norme{\eta^{0}}_{H^{3}(\Gamma_{s})}.
\end{align*}
This implies
\[
\norme{\frac{z\widetilde{\eta}_{xxx}\overline{u}_{1}}{1+\widetilde{\eta}}}_{L^{2}(Q_{T})}\leq C(T_{0},R,\mu,\|\xg^{0}\|_{\mathcal{X}(\Omega_{0})}).
\]
Thus \eqref{s41lemma5est1} is proved. For the Lipschitz estimate we use the same techniques. Let us make explicit the estimate on one of the terms, namely
\[
\frac{z\widetilde{\eta}_{1,x}}{1+\widetilde{\eta}_{1}}\overline{u}_{1,1}-\frac{z\widetilde{\eta}_{2,x}}{1+\widetilde{\eta}_{2}}\overline{u}_{2,1}
=z\left(\frac{\widetilde{\eta}_{1,x}}{1+\widetilde{\eta}_{1}}-\frac{\widetilde{\eta}_{2,x}}{1+\widetilde{\eta}_{2}}\right)\overline{u}_{1,1}+\frac{z\widetilde{\eta}_{2,x}}{1+\widetilde{\eta}_{2}}(\overline{u}_{1,1}-\overline{u}_{2,1}).
\]
Using the previous techniques and Lemma \ref{s41lemma1} we obtain
\[\norme{\frac{z\widetilde{\eta}_{2,x}}{1+\widetilde{\eta}_{2}}(\overline{u}_{1,1}-\overline{u}_{2,1})}_{H^{2,1}(Q_{T})}\leq C(T_{0},R,\mu,\|\xg^{0}\|_{\mathcal{X}(\Omega_{0})})\norme{\overline{\ug}_{1}-\overline{\ug}_{2}}_{\textbf{H}^{2,1}(Q_{T})}.\]
For the other term we write
\[
\frac{\widetilde{\eta}_{1,x}}{1+\widetilde{\eta}_{1}}-\frac{\widetilde{\eta}_{2,x}}{1+\widetilde{\eta}_{2}}=\frac{\widetilde{\eta}_{1,x}-\widetilde{\eta}_{2,x}}{(1+\widetilde{\eta}_{1})(1+\widetilde{\eta}_{2})}+\frac{\widetilde{\eta}_{1,x}(\widetilde{\eta}_{2}-\widetilde{\eta}_{1})}{(1+\widetilde{\eta}_{1})(1+\widetilde{\eta}_{2})}+\frac{\widetilde{\eta}_{1}(\widetilde{\eta}_{1,x}-\widetilde{\eta}_{2,x})}{(1+\widetilde{\eta}_{1})(1+\widetilde{\eta}_{2})},
\]
and 
\[\norme{z\left(\frac{\widetilde{\eta}_{1,x}}{1+\widetilde{\eta}_{1}}-\frac{\widetilde{\eta}_{2,x}}{1+\widetilde{\eta}_{2}}\right)\overline{u}_{1,1}}_{H^{2,1}(Q_{T})}\leq C(T_{0},R,\mu,\|\xg^{0}\|_{\mathcal{X}(\Omega_{0})})\norme{\eta_{1}-\eta_{2}}_{H^{4,2}(\Sigma_{T}^{s})}.\]
\end{proof}
The nonlinearities in (\ref{s31e1}) can now be estimated. 
\begin{thm}\label{s41thm1}Let $\xg^{0}$ belong to ${\tilde{\mathcal{X}}(\Omega_{0})}$, $R>0$, and $\mu>0$.   There exists a function $P_{\theta,n}(T)=\sum\nolimits_{k=0}^{n}T^{\theta_{k}}$ with $n\in\mathbb{N}^{*}$ and $\theta\in (\mathbb{R}^{*}_{+})^{n+1}$ such that, for all $0< T< T_{0}$ and all $(\overline{\ug},p,\eta)\in \mathcal{B}(\xg^{0},R,\mu,T)$, $(\textbf{F}(\overline{\ug},p,\eta),\Theta(\overline{\ug}),H(\overline{\ug},\eta))$ belongs to $W_{T}$ and the following estimate holds
\begin{equation}\label{s41thm1est1}
\norme{(\textbf{F}(\overline{\ug},p,\eta),\Theta(\overline{\ug}),H(\overline{\ug},\eta))}_{W_{T}}
\leq C(T_{0},R,\mu,\|\xg^{0}\|_{\mathcal{X}(\Omega_{0})})P_{\theta,n}(T).
\end{equation}
Moreover, for $(\overline{\ug}_{i},p_{i},\eta_{i})\in\mathcal{B}(\xg^{0},R,\mu,T)$ ($i=1,2$) the following estimate holds
\begin{equation}\label{s41thm1est2}
\begin{aligned}
&\norme{(\textbf{F}_{1},\Theta_{1},H_{1})-(\textbf{F}_{2},\Theta_{2},H_{2})}_{W_{T}}\\
&{}\leq C(T_{0},R,\mu,\|\xg^{0}\|_{\mathcal{X}(\Omega_{0})})P_{\theta,n}(T)\norme{(\overline{\ug}_{1},p_{1},\eta_{1})-(\overline{\ug}_{2},p_{2},\eta_{2})}_{\mathcal{Y}_{T}},
\end{aligned}
\end{equation}
with the notations $(\textbf{F}_{i},\Theta_{i},H_{i})=(\textbf{F}(\overline{\ug}_{i},p_{i},\eta_{i}),\Theta(\overline{\ug}_{i}),H(\overline{\ug}_{i},\eta_{i}))$.
\end{thm}
\begin{proof}
\textit{Step 1:} Estimate of $\textbf{F}(\overline{\ug},p,\eta)$.
We recall the form of $\textbf{F}(\overline{\ug},p,\eta)$:
\[\textbf{F}(\overline{\ug},p,\eta)=\textbf{G}(M(\overline{\ug},\eta),p,\eta)-\partial_tN(\overline{\ug},\eta) +\nu\Delta N(\overline{\ug},\eta).\]
Set $\ug=M(\overline{\ug},\eta)$. Thanks to Lemma \ref{s41lemma5} we can prove the estimates with $\ug$ and then obtain estimates in terms of $\overline{\ug}$. For
\begin{align*}
\textbf{G}(\ug,p,\eta)&=-\widetilde{\eta}\ug_t+\left[z\widetilde{\eta}_t+\nu z\left(\frac{\widetilde{\eta}_x^2}{1+\widetilde{\eta}}-\widetilde{\eta}_{xx}\right)\right]\ug_z\\
&+\nu\left[-2z\widetilde{\eta}_x\ug_{xz}+\widetilde{\eta}\ug_{xx}+\frac{z^2\widetilde{\eta}_x^2-\widetilde{\eta}}{1+\widetilde{\eta}}\ug_{zz}\right]\\
&+z(\widetilde{\eta}_x p_z-\widetilde{\eta} p_x)\textbf{e}_1-(1+\widetilde{\eta})u_1
\ug_x+(z\widetilde{\eta}_x
u_1-u_2)\ug_z,
\end{align*}
we use $L^{\infty}$ estimates on $\widetilde{\eta}$ to gain a factor $T$ for terms like the first one:
\begin{align*}\norme{-\widetilde{\eta} \ug_t}_{\textbf{L}^{2}(Q_{T})}&{}\leq \norme{\widetilde{\eta}}_{L^{\infty}(\Sigma^{s}_{T})}\norme{\ug_t}_{\textbf{L}^{2}(Q_{T})}\\
&{}\leq \mu T^{1/2}\norme{\eta_{t}}_{L^{2}(0,T;L^{\infty}(\Gamma_{s}))}\norme{\ug}_{\textbf{H}^{2,1}(Q_{T})}\leq C(T_{0},R,\mu,\|\xg^{0}\|_{\mathcal{X}(\Omega_{0})})T^{1/2}.
\end{align*}
For the product of functions in $L^{\infty}(0,T;H^{1}(\Omega_{0}))$ with derivatives of functions in $H^{2,1}(Q_{T})$ we use Lemma \ref{s41lemma4}, for example:
\begin{align*}
\norme{z\widetilde{\eta}_{t}\ug_{z}}_{\textbf{L}^{2}(Q_{T})}&{}\leq CT^{1/4}\norme{\widetilde{\eta}_{t}}_{L^{\infty}(0,T;H^{1}(\Gamma_{s}))}\norme{\ug}_{L^{\infty}(0,T;\textbf{H}^{1}(\Omega_{0}))}^{1/2}\norme{\ug}_{L^{2}(0,T;\textbf{H}^{2}(\Omega_{0}))}^{1/2}\\
&{}\leq C(T_{0},R,\mu,\|\xg^{0}\|_{\mathcal{X}(\Omega_{0})})T^{1/4},
\end{align*}
and
\begin{align*}
\norme{-\nu z\widetilde{\eta}_{xx}\ug_{z}}_{\textbf{L}^{2}(Q_{T})}&{}\leq CT^{1/4}\norme{\widetilde{\eta}_{xx}}_{L^{\infty}(0,T;H^{1}(\Omega_{0}))}\norme{\ug}_{L^{\infty}(0,T;\textbf{H}^{1}(\Omega_{0}))}^{1/2}\norme{\ug}_{L^{2}(0,T;\textbf{H}^{2}(\Omega_{0}))}^{1/2}\\
&{}\leq C(T_{0},R,\mu,\|\xg^{0}\|_{\mathcal{X}(\Omega_{0})})T^{1/4},
\end{align*}
where we have used the estimates on $\widetilde{\eta}_{t}$ and $\widetilde{\eta}_{xx}$ given in the proof of Lemma \ref{s41lemma5}.

The term $N(\overline{\ug},\eta)=\displaystyle\left(\frac{-\widetilde{\eta}\overline{u}_{1}}{1+\widetilde{\eta}},\frac{z\widetilde{\eta}_{x}\overline{u}_{1}}{1+\widetilde{\eta}}\right)^{T}$ has already been estimated in the proof of Lemma \ref{s41lemma5} and the factor $T$ is obtained with the previous techniques and Lebesgue interpolation for the terms
\[\frac{z(\eta_{xxx}-\eta^{0}_{xxx})\overline{u}_{1}}{(1+\eta^{0})(1+\widetilde{\eta})}.\]
\textit{Step 2:} Estimate of $\Theta(\overline{\ug})$.
In order to obtain an estimate in $L^{2}(0,T;H^{1/2}(\Gamma_{i,o}))$ we study $\Theta(\overline{\ug})=(1/2)\vert\overline{\ug}\vert^{2}$ on $\Omega_{0}$ and then look for the restriction to $\Gamma_{i,o}$. We have
\begin{align*}
\Theta(\overline{\ug})&=\frac{\overline{u}_{1}^{2}+\overline{u}_{2}^{2}}{2},\\
\Theta(\overline{\ug})_{x}&=\overline{u}_{1}\overline{u}_{1,x}+\overline{u}_{2}\overline{u}_{2,x},\\
\Theta(\overline{\ug})_{z}&=\overline{u}_{1}\overline{u}_{1,z}+\overline{u}_{2}\overline{u}_{2,z},\\
\end{align*}
and using Lemma \ref{s41lemma4}
\[
\norme{\overline{u}_{1}\overline{u}_{1,x}}_{L^{2}(Q_{T})}\leq CT^{1/4}\norme{\overline{\ug}}_{L^{\infty}(0,T;\textbf{H}^{1}(\Omega_{0}))}^{3/2}\norme{\overline{\ug}}_{L^{2}(0,T;\textbf{H}^{2}(\Omega_{0}))}
\leq C(T_{0},R,\mu,\|\xg^{0}\|_{\mathcal{X}(\Omega_{0})})T^{1/4}.
\]
which implies a $L^{2}(0,T;H^{1}(\Omega_{0}))$ estimate on $\Theta(\overline{\ug})$ and thus a $L^{2}(0,T;H^{1/2}(\Gamma_{i,o}))$ estimate for the trace.

\textit{Step 3:} Estimate of $H(\overline{\ug},\eta)$.
We recall that $H(\overline{\ug},\eta)=\Psi(M(\overline{\ug},\eta),\eta)$ with
\[\Psi(\textbf{u},\eta)=\nu\left(\frac{\eta_x}{1+\widetilde{\eta}}u_{1,z}+\eta_x u_{2,x}-\frac{\widetilde{\eta}_x\eta_{x}z-2}{1+\widetilde{\eta}}u_{2,z}\right).\]
For the terms without $\widetilde{\eta}$ we use directly the regularity of $\ug$ to gain a factor $T$. Using fractional Sobolev embeddings \cite[Theorem 7.58]{MR0450957} and trace theorems we know that $\ug_{x}(x,1+\eta^{0}(x),t),\ug_{z}(x,1+\eta^{0}(x),t)$ belong to $H^{1/2,1/4}(Q_{T})$ and $L^{4}(0,T;\textbf{L}^{2}(\Gamma_{s}))$. Hence
\[\norme{\nu\frac{\eta_x}{1+\widetilde{\eta}}u_{1,z}}_{L^{2}(0,T;L^{2}(\Gamma_{s}))}\leq 
C(\mu,\eta^{0})\norme{\eta_x u_{1,z}}_{L^{2}(0,T;L^{2}(\Gamma_{s}))},\]
and
\begin{align*}
\norme{\eta_x u_{1,z}}_{L^{2}(\Sigma^{s}_{T})}^{2}&{}\leq \norme{\eta_{x}}_{L^{\infty}(\Sigma^{s}_{T})}^{2}\int_{0}^{T}\norme{u_{1,z}}_{L^{2}(\Gamma_{s})}^{2}dt\\
&{}\leq \norme{\eta_{x}}_{L^{\infty}(\Sigma^{s}_{T})}^{2}
T^{1/2}\norme{u_{1,z}}_{L^{4}(0,T;L^{2}(\Gamma_{s}))}^{2}\\
&{}\leq C(T_{0},R,\mu,\|\xg^{0}\|_{\mathcal{X}(\Omega_{0})})T^{1/2}.
\end{align*}
The other terms are estimated with the previous techniques.

\textit{Step 4:} Lipschitz estimates.
The Lipschitz estimates are obtained with the same techniques. Let us make explicit some inequalities.
\begin{align*}
&\norme{\widetilde{\eta}_{1}\ug_{1,t}-\widetilde{\eta}_{2}\ug_{2,t}}_{\textbf{L}^{2}(Q_{T})}\\
&{}\leq \norme{(\widetilde{\eta}_{1}-\widetilde{\eta}_{2})\ug_{1,t}}_{\textbf{L}^{2}(Q_{T})}+\norme{\widetilde{\eta}_{2}(\ug_{1,t}-\ug_{2,t})}_{\textbf{L}^{2}(Q_{T})}\\
&{}\leq \norme{\widetilde{\eta}_{1}-\widetilde{\eta}_{2}}_{L^{\infty}(\Sigma^{s}_{T})}\norme{\ug_{1,t}}_{\textbf{L}^{2}(Q_{T})}+\norme{\widetilde{\eta}_{2}}_{L^{\infty}(\Sigma^{s}_{T})}\norme{\ug_{1,t}-\ug_{2,t}}_{\textbf{L}^{2}(Q_{T})}\\
&{}\leq C(T_{0},R,\mu,\|\xg^{0}\|_{\mathcal{X}(\Omega_{0})})T^{1/2}(\norme{\eta_{1}-\eta_{2}}_{H^{4,2}(\Sigma^{s}_{T})}+\norme{\overline{\ug}_{1}-\overline{\ug}_{2}}_{\textbf{H}^{2,1}(Q_{T})}).
\end{align*}
All the interest of working in the initial domain $\Omega_{0}$ instead of the rectangular $\Omega$ comes from the estimate (\ref{s41lemma3}) on $\widetilde{\eta}$. With the usual change of variables, the term $\eta_{2}(\ug_{1}-\ug_{2})$ cannot be estimated without smallness assumption on $\eta^{0}$. For $\nu z\widetilde{\eta}_{1,xx}\ug_{1,z}-\nu z\widetilde{\eta}_{2,xx}\ug_{2,z}$ we have
\begin{align*}
&\norme{\nu z(\widetilde{\eta}_{1,xx}-\widetilde{\eta}_{2,xx})\ug_{1,z}}_{\textbf{L}^{2}(Q_{T})}\\
&{}\leq CT^{1/4}\norme{\widetilde{\eta}_{1,xx}-\widetilde{\eta}_{2,xx}}_{L^{\infty}(0,T;H^{1}(\Omega_{0}))}\norme{\ug_{1}}_{L^{\infty}(0,T;\textbf{H}^{1}(\Omega_{0}))}^{1/2}\norme{\ug_{1}}_{L^{2}(0,T;\textbf{H}^{2}(\Omega_{0}))}^{1/2}\\
&{}\leq C(T_{0},R,\mu,\|\xg^{0}\|_{\mathcal{X}(\Omega_{0})})T^{1/4}\norme{\eta_{1}-\eta_{2}}_{H^{4,2}(\Sigma^{s}_{T})},
\end{align*}
and
\begin{align*}
&\norme{\nu z\widetilde{\eta}_{2,xx}(\ug_{1,z}-\ug_{2,z})}_{\textbf{L}^{2}(Q_{T})}\\
&{}\leq CT^{1/4}\norme{\widetilde{\eta}_{2,xx}}_{L^{\infty}(0,T;H^{1}(\Omega_{0}))}\norme{\ug_{1}-\ug_{2}}_{L^{\infty}(0,T;\textbf{H}^{1}(\Omega_{0}))}^{1/2}\norme{\ug_{1}-\ug_{2}}_{L^{2}(0,T;\textbf{H}^{2}(\Omega_{0}))}^{1/2}\\
&{}\leq C(T_{0},R,\mu,\|\xg^{0}\|_{\mathcal{X}(\Omega_{0})})T^{1/4}\norme{\overline{\ug}_{1}-\overline{\ug}_{2}}_{\textbf{H}^{2,1}(Q_{T})}.
\end{align*}
The Lipschitz estimates with lower regularity terms like $\widetilde{\eta}_{xxx}$ are obtained as in the proof of Lemma \ref{s41lemma5}.
\end{proof}
\subsection{Fixed point procedure}
For all $\xg^{0}=(\ug^{0},\eta_1^0,\eta_2^0)\in \tilde{\mathcal{X}}(\Omega_{0})$, $R>0$, $\mu>0$ and $T>0$ such that $\mathcal{B}(\xg^{0},R,\mu,T)\neq \emptyset$, consider the map
\begin{equation}\label{s42map1}
\mathcal{F}:
\begin{cases}
\begin{array}{lr}
\mathcal{B}(\xg^{0},R,\mu,T)\longrightarrow \mathcal{Y}_{T},\\
(\overline{\ug},p,\eta)\mapsto (\overline{\ug}^{*},p^{*},\eta^{*}),
\end{array}
\end{cases}
\end{equation}
where $(\overline{\ug}^{*},p^{*},\eta^{*})$ is the solution to \eqref{s32e3}  with right-hand side $(\textbf{F}(\overline{\ug},p,\eta),\Theta(\overline{\ug}),H(\overline{\ug},\eta))$. In order to solve (\ref{s32e2}) we look for a fixed point of the map $\mathcal{F}$.
\begin{thm}\label{s42thm1}
For all $\xg^{0}=(\ug^{0},\eta_1^0,\eta_2^0)$ in $\tilde{\mathcal{X}}(\Omega_{0})$, there exist $R>0$, $\mu>0$ and $T>0$ such that $\mathcal{F}$ is a contraction from $\mathcal{B}(\xg^{0},R,\mu,T)$ into $\mathcal{B}(\xg^{0},R,\mu,T)$. Hence $\mathcal{F}$ has a fixed point.
\end{thm}

\begin{proof}
Set $\mu=\norme{(1+\eta^{0})^{-1}}_{L^{\infty}(\Gamma_{s})}$. Let $(\overset{\circ}{\ug},\overset{\circ}{\eta},\overset{\circ}{p})$ be the solution on $[0,T_{0}]$ to (\ref{s32e3}) with right-hand side $0$. We choose $R_1$ such that $\norme{(\overset{\circ}{\ug},\overset{\circ}{\eta},\overset{\circ}{p})}_{\mathcal{X}_{T_{0}}}\leq R_1$. 
Writing 
\[1+\overset{\circ}{\eta}(t)=1+\eta^{0}+\int_{0}^{t}\overset{\circ}{\eta}_{t}(s)ds,\]
we can choose $0<T_{1}< T_{0}$ such that
\begin{equation}\label{continuityeta}
\norme{(1+\overset{\circ}{\eta})^{-1}}_{L^{\infty}(\Sigma^{s}_{T_{1}})}\leq \frac{1}{\norme{1+\eta^{0}}_{L^{\infty}(\Gamma_{s})}-T_{1}C(\mu,T_{0},R_1,\|\xg^{0}\|_{\mathcal{X}(\Omega_{0})})}\leq 2\mu.
\end{equation}
Thus the ball $\mathcal{B}(\xg^{0},R_1,\mu,T_{1})$ is non-empty. From 
Theorems \ref{s41thm1} and \ref{s5thm5}, it follows that 
\begin{align*}
\norme{\mathcal{F}(\overline{\ug},p,\eta)}_{\mathcal{Y}_{T_{1}}}\leq C_{L}\left(\norme{(\ug^{0},\eta_1^0,\eta_2^0)}_{\mathcal{X}(\Omega_{0})}+C(T_{0},R_1,\mu,\|\xg^{0}\|_{\mathcal{X}(\Omega_{0})})P_{\theta,n}(T_{1})\right), 
\end{align*}
for all $(\overline{\ug},p,\eta)\in\mathcal{B}(\xg^{0},R_1,\mu,T_{1})$. 

Then we choose $R_2\geq R_1$ such that $R_2 \geq 2C_{L}\norme{(\ug^{0},\eta_1^0,\eta_2^0)}_{\mathcal{X}(\Omega_{0})}$, and $T_{2}\leq T_{1}$ such that 
\[
C_{L}C(T_{0},R_2,\mu,\|\xg^{0}\|_{\mathcal{X}(\Omega_{0})})P_{\theta,n}(T_{2})\leq C_{L}\norme{(\ug^{0},\eta_1^0,\eta_2^0)}_{\mathcal{X}(\Omega_{0})}\text{ and }\norme{(1+\eta^{*})^{-1}}_{L^{\infty}(\Sigma^{s}_{T_{2}})}\leq 2\mu.
\]
Therefore $\mathcal{F}$ is well defined from $\mathcal{B}(\xg^{0},R_2,\mu,T_2)$ into $\mathcal{B}(\xg^{0},R_2,\mu,T_{2})$. Still with Theorems \ref{s41thm1} and \ref{s5thm5}, it follows that 
\begin{align*}
&\norme{\mathcal{F}(\overline{\ug}_{1},p_{1},\eta_{1})-\mathcal{F}(\overline{\ug}_{2},p_{2},\eta_{2})}_{\mathcal{Y}_{T_{2}}}\\
&{}\leq C_{L}C(T_{0},R_2,\mu,\|\xg^{0}\|_{\mathcal{X}(\Omega_{0})})P_{\theta,n}(T_{2})\norme{(\overline{\ug}_{1},p_{1},\eta_{1})-(\overline{\ug}_{2},p_{2},\eta_{2})}_{\mathcal{Y}_{T_{2}}}, 
\end{align*}
for all $(\overline{\ug}_{1},p_{1},\eta_{1})$ and $(\overline{\ug}_{2},p_{2},\eta_{2})$ belonging to $\mathcal{B}(\xg^{0},R_2,\mu,T_2)$. We choose $0<T_{3}\leq T_{2}$ such that $C_{L}C(T_{0},R_2,\mu,\|\xg^{0}\|_{\mathcal{X}(\Omega_{0})})P_{\theta,n}(T_{3})\leq 1/2$. 
The mapping $\mathcal{F}$ is a contraction from the complete metric space $\mathcal{B}(\xg^{0},R_2,\mu,T_3)$ into itself, and the Banach fixed point theorem concludes the proof.
\end{proof}
\begin{thm}\label{s42thm2} For all $(\ug^{0},\eta_1^0,\eta_2^0)\in\tilde{\mathcal{X}}(\Omega_{0})$ there exists $T>0$ such that the system (\ref{s32e2}) has a unique strong solution $(\overline{\ug},p,\eta)$ in $\textbf{H}^{2,1}(Q_T)\times L^{2}(0,T;H^{1}(\Omega_{0}))\times H^{4,2}(\Sigma^{s}_{T})$.
\end{thm}
\begin{proof}
The existence is already proved. Set $\xg^{0}=(\ug^{0},\eta_1^0,\eta_2^0)$ and let $(\overline{\ug},p,\eta)$ be the unique solution to (\ref{s32e2}) in $\mathcal{B}(\xg^{0},R,\mu,T)$ with
\[\mu=\norme{(1+\eta^{0})^{-1}}_{L^{\infty}(\Gamma_{s})},\,\,\,R = 2C_{L}\norme{(\ug^{0},\eta_1^0,\eta_2^0)}_{\mathcal{X}(\Omega_{0})},\]
and $T>0$, constructed by the fix point method in the previous Theorem. Let $(\overline{\ug}',p',\eta')$ be another solution to (\ref{s32e2}) defined on $[0,T]$ with the same initial data. Define the constants $R_{0}=\norme{(\overline{\ug}',p',\eta')}_{\mathcal{Y}_{T}}$ and $\mu_{0}=\norme{(1+\eta')^{-1}}_{L^{\infty}(\Sigma_{T}^{s})}$. Assume that 
$T>0$ is small enough such that $\norme{(1+\eta')^{-1}}_{L^{\infty}(\Sigma_{T}^{s})}\leq 2\mu$. 
From Theorems \ref{s41thm1} and \ref{s5thm5}, it follows that 
\[
\norme{(\overline{\ug}',p',\eta')}_{\mathcal{Y}_{T_{1}}}\leq C_{L}\left(\norme{(\ug^{0},\eta_1^0,\eta_2^0)}_{\mathcal{X}(\Omega_{0})}+C(T_{0},R_{0},\mu_{0},\|\xg^{0}\|_{\mathcal{X}(\Omega_{0})})P_{\theta,n}(T_{1})\right), 
\]
for all $0<T_1\leq T$. Let us choose $0< T_{1}\leq T$ such that
$$
C_{L}C(T_{0},R_{0},\mu_{0},\|\xg^{0}\|_{\mathcal{X}(\Omega_{0})})P_{\theta,n}(T_{1})\leq C_{L}\norme{(\ug^{0},\eta_1^0,\eta_2^0)}_{\mathcal{X}(\Omega_{0})}. 
$$
Hence $(\overline{\ug}',p',\eta')$ belongs to $\mathcal{B}(\xg^{0},R,\mu,T_{1})$ and $(\overline{\ug},p,\eta)=(\overline{\ug}',p',\eta')$ on $[0,T_{1}]$. 

Let $0<T^{*}\leq T$ be the greatest time such that the two solutions are equal. We then consider the system (\ref{s32e2}) starting at the time $T^{*}$, rewritten in $\Omega_{\eta(T^{*})}$, with the initial conditions $(\overline{\ug}(T^{*}),\eta(T^{*}),\eta_{t}(T^{*}))=(\overline{\ug}'(T^{*}),\eta'(T^{*}),\eta_{t}'(T^{*}))$. If $T^{*}<T$, using the fixed point procedure we prove the existence of a solution $(\overline{\ug}'',\eta'',\eta_{t}'')$ on $[T^{*},T_{2}]$ with $T_{2}>T^{*}$. The previous argument shows that there exists $T_{3}>0$ such that the three solutions are equal (after a change of variable in order to consider functions in the domain $\Omega_{\eta(T^{*})}$) on $[T^{*},T_{3}]$ which is a contradiction with the definition of $T^{*}$. Hence $T^{*}=T$ and the solution to (\ref{s32e2}) is unique.
\end{proof}
The previous ideas and techniques can be applied on system \eqref{s1main} with the Dirichlet boundary conditions $\ug=0$ on $\Gamma_{i,o}$ and thus fix the gap in the proof of local existence in \cite{MR2765696}.

To conclude this section, we state the existence and uniqueness of a solution for (\ref{s32e2}) on $[0,T]$, with $T>0$ a fixed time and smallness assumptions on the initial data. This result is proved on the rectangular domain since the estimates of the nonlinear terms are done through the radius of the ball in the fixed point argument. The existence technique is similar to the one in \cite[Theorem 10.1]{MR2745779} and the uniqueness comes from the local existence and uniqueness result. Let us notice that with this approach the nonlinear term in the beam equation in (\ref{s31e1}) with $\eta^{0}=0$ writes
\begin{align*}
\Psi(\widehat{\textbf{u}},\eta)&{}=\nu\left(\frac{\eta_x}{1+\eta}\widehat{u}_{1,z}+\eta_x \widehat{u}_{2,x}-\frac{\eta_x^{2}z-2}{1+\eta}\widehat{u}_{2,z}\right)\\
&{}=-2\nu \widehat{u}_{2,z}+\nu\left(\frac{\eta_x}{1+\eta}\widehat{u}_{1,z}+\eta_x u_{2,x}-\frac{\eta_x^{2}z-2\eta}{1+\eta}\widehat{u}_{2,z}\right)\\
&{}=-2\nu \widehat{u}_{2,z} + \overline{\Psi}(\widehat{\ug},\eta).
\end{align*}
After writing $\widehat{\ug}=M(\overline{\ug},\eta)=\overline{\ug}+N(\overline{\ug},\eta)$ this nonlinear term becomes
\[H(\overline{\ug},\eta)=-2\nu \overline{u}_{2,z} -2\nu N(\overline{\ug},\eta)_{2,z} + \overline{\Psi}(M(\overline{\ug},\eta),\eta),\]
and as $\text{div }\overline{\ug}=\overline{u}_{1,x}+\overline{u}_{2,z}=0$ in $Q_{T}$ and $\overline{u}_{1,x}=0$ on $\Sigma^{s}_{T}=\Sigma^{0}_{T}$ we obtain $\overline{u}_{2,z}=0$ on $\Sigma^{s}_{T}$. Hence all the nonlinear terms in the beam equation are at least quadratic.
\begin{thm}\label{s42thm3}Let $T>0$ be a fixed time and recall that $\Omega=(0,L)\times(0,1)$. There exists $r>0$ such that for all $(\ug^{0},\eta_1^0,\eta_2^0)$ in $\mathcal{X}(\Omega)$ satisfying $\norme{(\ug^{0},\eta_1^0,\eta_2^0)}_{\mathcal{X}(\Omega)}\leq r$, the system (\ref{s32e2}) admits a unique solution in $\textbf{H}^{2,1}(Q_{T})\times L^{2}(0,T;H^{1}(\Omega))\times H^{4,2}(\Sigma^{s}_{T})$.
\end{thm}
\begin{rmq}Note that the initial condition is taken in $\mathcal{X}(\Omega)$, not $\tilde{\mathcal{X}}(\Omega)$, which means that $\eta^{1}_{0}$ can be different from $0$.
\end{rmq}
\section{Appendix} \label{sec6}
\subsection{Steady Stokes equations}
Consider the steady Stokes equations
\begin{equation}\label{s61e3}
\begin{aligned}
&-\nu\Delta\ug + \nabla p={\textbf{f}},\,\,\,\text{div }\textbf{u}=0\,\text{ in }\Omega_{0},\\
&\textbf{u}=\textbf{g}\,\text{ on }\Gamma_{0},\,\ug=0\,\text{ on }\Gamma_{b},\,u_2=0\,\text{ and }\,p=h\,\text{ on }\Gamma_{i,o},\\
\end{aligned}
\end{equation}
with $\textbf{f}\in \textbf{L}^{2}(\Omega_{0})$, $\textbf{g}=(0,g)^{T}\in\mathcal{H}^{3/2}_{00}(\Gamma_{0})$ and $h\in H^{1/2}(\Gamma_{i,o})$. We prove in Theorem \ref{s61thm5} the existence and uniqueness of a pair $(\ug,p)\in \textbf{H}^{2}(\Omega_{0})\times H^{1}(\Omega_{0})$ solution to (\ref{s61e3}). An existence and uniqueness result for (\ref{s61e3}) with weaker data is given in Theorem \ref{s61thm6}. The nonhomogeneous boundary condition on the pressure is handled directly with a lifting operator $\mathcal{R}\in \mathcal{L}(H^{1/2}(\Gamma_{i,o}),H^{1}(\Omega_{0}))$. For the nonhomogeneous Dirichlet boundary condition we use the following theorem.
\begin{thm}\label{s61thm1}For all $\textbf{g}=(0,g)^{T}\in\mathcal{H}^{3/2}_{00}(\Gamma_{0})$ the system
\begin{equation}\label{s61e2}
\begin{array}{l}
\text{div }\textbf{w}=0\,\text{ in }\Omega_{0},\\
\textbf{w}=\textbf{g}\,\text{ on }\Gamma_{0},\,\textbf{w}=0\,\text{ on }\Gamma_{b},\,w_2=0\,\text{ on }\Gamma_{i,o},\\
\end{array}
\end{equation}
admits a solution $\textbf{w}\in \textbf{H}^{2}(\Omega_{0})$ satisfying the estimate
\[\norme{\textbf{w}}_{\textbf{H}^{2}(\Omega_{0})}\leq C\norme{\textbf{g}}_{\mathcal{H}^{3/2}_{00}(\Gamma_{0})}.\]
\end{thm}
\begin{proof}
We look for $\textbf{w}$ under the form $\textbf{w}=(-\partial_{2}\phi,\partial_{1}\phi)^{T}$, which ensures the property $\text{div }\textbf{w}=0$. The boundary conditions on $\textbf{w}$ imply the following conditions on $\phi$
\begin{equation}\label{s61thm1e1}
\begin{array}{lr}
\partial_{2}\phi=0\text{ and }\partial_{1}\phi=g\text{ on }\Gamma_{0},\\
\frac{\partial\phi}{\partial\normal}=\partial_{1}\phi=0\text{ on }\Gamma_{i,o},\,\partial_{2}\phi=\partial_{1}\phi=0\text{ on }\Gamma_{b}.
\end{array}
\end{equation}
Let $\eta^{0}_{e}$ be an $H^{3}(\mathbb{R})$ extension of $\eta^{0}$. We consider the change of variables
\[\psi^{\pm} :
\begin{cases}
\begin{array}{lr}
\mathbb{R}^{2}\longrightarrow \mathbb{R}^{2}\\
(x,y)\mapsto (x,y\pm\eta^{0}_{e}(x)).\\
\end{array}
\end{cases}
\]
Let $\widehat{v}$ be a function in $H^{3}(\mathbb{R}^{2})$. Thanks to the $H^{3}$-regularity of $\eta^{0}_{e}$, the function $\widehat{v}\circ\psi^{\pm}$ is still in $H^{3}(\mathbb{R}^{2})$. We search for $\phi$ solution to (\ref{s61thm1e1}) under the form $\phi=\widehat{\phi}\circ\psi^{-}$ with $\widehat{\phi}\in H^{3}((0,L)\times(-\infty,1))$ satisfying
\begin{equation}\label{s61thm1e2}
\begin{array}{lr}
\partial_{2}\widehat{\phi}=0\text{ and }\partial_{1}\widehat{\phi}=\widehat{g}\text{ on }\Gamma_{s},\\
\partial_{1}\widehat{\phi}=0\text{ on }\Gamma_{i,o},\,\widehat{\phi}=0\text{ on }(0,L)\times(-\infty,1-\delta),\\
\end{array}
\end{equation}
with $\widehat{g}=g\circ\psi^{+}$ and
\begin{equation}\label{delta}
\delta=
\begin{cases}
\begin{array}{lr}
\displaystyle\min_{x\in(0,L)}(1+\eta^{0}(x))\text{ if }\min_{x\in(0,L)}(1+\eta^{0}(x))<1,\\
\displaystyle \alpha \text{ if }\min_{x\in(0,L)}(1+\eta^{0}(x))\geq 1,\\
\end{array}
\end{cases}
\end{equation}
for a fixed $\alpha\in(0,1)$. This condition is used to ensure that the function $\phi=\widehat{\phi}\circ\psi^{-}$ is equal to zero near $\Gamma_{b}$, in order to fulfil the boundary conditions $\partial_{1}\phi=\partial_{2}\phi=0$ on $\Gamma_{b}$. To build $\widehat{\phi}$ we first search for $\widehat{\phi}_{o}$ such that
\begin{equation}\label{s61thm1e3}
\begin{array}{lr}
\frac{\partial\widehat{\phi}_{o}}{\partial\normal}=0\text{ on }\Gamma_{s}\cup\Gamma_{o},\\
\widehat{\phi}_{o}(x,y)=G(x,y)=\int_{0}^{x}\widehat{g}(s)ds\text{ for }(x,y)\in \Gamma_{s},\\
\widehat{\phi}_{o}=0\text{ on }(0,L)\times(-\infty,1-\delta),
\end{array}
\end{equation}
The boundary conditions on $\Gamma_{s}$ are handled directly thanks to a lifting and a symmetry argument is used to obtain the homogeneous Neumann boundary condition on $\Gamma_{o}$. We set
\[G^{*}:
\begin{cases}
\begin{array}{lr}
G^{*}(x,y)=G(x,y)\text{ for }(x,y)\in \Gamma_{s},\\
G^{*}(x,y)=G(2L-x,y)\text{ for }(x,y)\in (L,2L)\times\{1\}.\\
\end{array}
\end{cases}
\]
Denote by $\widehat{g}_{s}$ the odd extension of $\widehat{g}$ on $\Gamma_{s,s}=(0,2L)\times\{1\}$. As $\widehat{g}\in H^{3/2}_{00}(\Gamma_{s})$, the function $\widehat{g}_{s}$ belongs to $H^{3/2}(\Gamma_{s,s})$. Indeed odd and even symmetries preserve the $H^1$-regularity (resp. $H^2$-regularity) for functions in $H^{1}_{0}(\Gamma_{0})$ (resp. in $H^{2}_{0}(\Gamma_{0})$), thus, by interpolation, the $H^{3/2}$-regularity is also preserved for functions in $H^{3/2}_{00}(\Gamma_{0})=[H^{1}_{0}(\Gamma_{0}),H^{2}_{0}(\Gamma_{0})]_{1/2}$.

As $\partial_{1}G^{*}(\cdot,1)=\widehat{g}_{s}(\cdot)$ we have $G^{*}\in H^{5/2}(\Gamma_{s,s})$. We still denote by $G^{*}$ a regular extension of $G^{*}$ on $\mathbb{R}\times\{1\}$. The lifting results in \cite{MR0350178} in the case of the half-plan give a function $\widehat{\phi}_{1}\in H^{3}(\mathbb{R}\times(-\infty,1))$ such that $\widehat{\phi}_{1}=G^{*}$ and $\frac{\partial\widehat{\phi}_{1}}{\partial\normal}=0$ on $\mathbb{R}\times\{1\}$. We then use cut-off functions to ensure that $\widehat{\phi}_{1}=0$ on $(0,2L)\times(-\infty,1-\delta)$.

Introduce the symmetric function $\widehat{\phi}_{2}$ to $\widehat{\phi}$ with respect to the axis $x=L$ defined by $\widehat{\phi}_{2}(x,y)=\widehat{\phi}_{2}(2L-x,y)$ for $(x,y)\in (0,2L)\times(-\infty,1)$. As the Dirichlet boundary condition $G^{*}$ is symmetric, $\widehat{\phi}_{2}$ satisfies the same boundary conditions as $\widehat{\phi}_{1}$ on $\Gamma_{s,s}$. We finally set $\widehat{\phi}_{o}=\frac{\widehat{\phi}_{1}+\widehat{\phi}_{1,s}}{2}$. The function $\widehat{\phi}_{o}$ belongs to $H^{3}((0,2L)\times(-\infty,1))$ and admits $x=L$ as an axis of symmetry. Hence we have $\frac{\partial\widehat{\phi}_{o}}{\partial\normal}=0$ on $\Gamma_{o}$ and the restriction on $(0,L)\times(-\infty,1)$ is a solution to (\ref{s61thm1e3}).

Using the same tools we obtain a function $\widehat{\phi}_{i}\in H^{3}((0,L)\times(-\infty,1))$ such that
\begin{equation}\label{s61thm1e4}
\begin{array}{lr}
\frac{\partial\widehat{\phi}_{i}}{\partial\normal}=0\text{ on }\Gamma_{s}\cup\Gamma_{i},\\
\widehat{\phi}_{i}(x,y)=G(x,y)=\int_{0}^{x}\widehat{g}(s)ds\text{ for }(x,y)\in \Gamma_{s},\\
\widehat{\phi}_{i}=0\text{ on }(0,L)\times(-\infty,1-\delta).
\end{array}
\end{equation}
Then we combine $\widehat{\phi}_{o}$ and $\widehat{\phi}_{i}$. Let $\alpha$ be a function defined on $[0,L]$ such that $\alpha=1$ near $\Gamma_{i}$, $\alpha=0$ near $\Gamma_{o}$ and $\alpha\in \mathcal{C}^{\infty}([0,L])$. The function $\widehat{\phi}$ defined by
\[\widehat{\phi}(x,y)=\alpha(x) \widehat{\phi}_{i}(x,y) + (1-\alpha(x))\widehat{\phi}_{o}(x,y)\text{ for all }(x,y)\in (0,L)\times(-\infty,1),\]
is a solution to (\ref{s61thm1e2}). Finally the restriction to $\Omega_{0}$ of the function $\phi=\widehat{\phi}\circ\psi^{-}$ is a solution to (\ref{s61thm1e1}). Indeed,
\[
\begin{array}{lr}
\partial_{2}\phi=\partial_{2}\widehat{\phi}\circ\psi^{-}=0\,\text{ on }\Gamma_{0},\\
\partial_{1}\phi=\partial_{1}\widehat{\phi}\circ\psi^{-} -\eta^{0}_{x}\partial_{2}\widehat{\phi}\circ\psi^{-}=\partial_{1}\widehat{\phi}\circ\psi^{-}=\widehat{g}\circ\psi^{-}=g\text{ on }\Gamma_{0},\\
\frac{\partial\phi}{\partial\normal}=\partial_{1}\phi=0\text{ on }\Gamma_{i,o},\,\partial_{2}\phi=\partial_{1}\phi=0\text{ on }\Gamma_{b}.
\end{array}
\]
and $\textbf{w}=(-\partial_2 \phi,\partial_1 \phi)^{T}$ is a solution of (\ref{s61e2}). We have $\textbf{w}\in \textbf{H}^{2}(\Omega_{0})$ and the estimate follows from the continuity of the lifting operator in \cite{MR0350178}.
\end{proof}
Let $\textbf{w}\in \textbf{H}^{2}(\Omega_{0})$ be the lifting of $\textbf{g}$ given by Theorem \ref{s61thm1} and $H=\mathcal{R}(h)$. By setting $(\vg,q)=(\ug,p)-(\textbf{w},H)$ the Stokes system (\ref{s61e3}) is equivalent to
\begin{equation}\label{s61thm2e}
\begin{aligned}
&-\nu\Delta \vg +\nabla q=\overline{\textbf{f}},\,\,\,\,\text{div }\vg=0\,\text{ in }\Omega_{0},\\
&\vg=0\,\text{ on }\Gamma_{d},\,\,\,v_{2}=0\,\text{ and }q=0\,\text{ on }\Gamma_{i,o},\\
\end{aligned}
\end{equation}
with $\overline{\textbf{f}}=\textbf{f}+\nu\Delta \textbf{w}-\nabla H$. Using Green formula one can derive the following variational formulation for (\ref{s61thm2e}).
\begin{thm}\label{s61thm3}Let $(\vg,q)\in\textbf{H}^{2}(\Omega_{0})\times H^{1}(\Omega_{0})$ be a solution to (\ref{s61thm2e}). Then $\vg$ satisfies the variational formulation :

Find $\vg\in V$ such that $\quad\displaystyle\nu\int_{\Omega_{0}}\nabla \textbf{v}:\nabla \boldsymbol{\varphi}=\int_{\Omega_{0}}\overline{\textbf{f}}\cdot \boldsymbol{\varphi}\text{ for all }\boldsymbol{\varphi}\in V.\quad (\star)$
\end{thm}
\begin{thm}\label{s61thm4}
The variational formulation $(\star)$ admits a unique solution $\vg\in V$. Moreover there exists a pressure $\mathcal{Q}\in L^{2}(\Omega_{0})$, unique up to an additive constant, such that $\displaystyle -\nu\Delta \vg +\nabla \mathcal{Q} = \overline{\textbf{f}}\,$ in $\textbf{H}^{-1}$. 
\end{thm}
The pressure $\mathcal{Q}$ is mentioned as a pressure associated with $\vg$.
\begin{proof}As the only constant in $V$ is the null function we can use a Poincar\'e inequality to prove that the bilinear form
\[a(\textbf{v},\boldsymbol{\varphi})=\nu\int_{\Omega_{0}}\nabla \textbf{v}:\nabla \boldsymbol{\varphi},\]
is coercive on $V$. Hence the Lax-Milgram lemma gives us the existence of a unique solution $\vg\in V$ to the variational formulation $(\star)$. For the pressure, we use the equality
\[\left\langle -\nu\Delta\vg - \overline{\textbf{f}},\boldsymbol{\varphi} \right\rangle_{\textbf{H}^{-1},\textbf{H}^{1}_{0}}=0 \text{, for all }\boldsymbol{\varphi}\in ({H}^{1}_{0}(\Omega_{0}))^{2}\text{ such that }\text{div }\boldsymbol{\varphi}=0, \]
and \cite[Chap 4, Theorem 2.3]{MR2986590} to prove the existence of $\mathcal{Q}\in L^{2}(\Omega_{0})$, unique up to an additive constant and such that $-\nu\Delta \vg +\nabla \mathcal{Q} = \overline{\textbf{f}}$ in $\textbf{H}^{-1}$.
\end{proof}
We now state the main theorem of this section.
\begin{thm}\label{s61thm5}For all $(\textbf{f},\textbf{g},h)\in\textbf{L}^{2}(\Omega_{0})\times\mathcal{H}^{3/2}_{00}(\Gamma_0)\times H^{1/2}(\Gamma_{i,o})$ the equation (\ref{s61e3}) admits a unique solution $(\ug,p)\in \textbf{H}^{2}(\Omega_{0})\times H^{1}(\Omega_{0})$. This solution satisfies the estimate
\[\norme{\ug}_{\textbf{H}^{2}(\Omega_{0})}+\norme{p}_{H^{1}(\Omega_{0})}\leq C(\norme{\textbf{f}}_{\textbf{L}^{2}(\Omega_{0})}+\norme{\textbf{g}}_{\mathcal{H}^{3/2}_{00}(\Gamma_0)}+\norme{h}_{H^{1/2}(\Gamma_{i,o})}).\]
\end{thm}
\begin{proof}
Let us work directly on the homogeneous system (\ref{s61thm2e}). We prove the existence of a unique pair $(\vg,q)\in \textbf{H}^{2}(\Omega_{0})\times H^{1}(\Omega_{0})$ solution to this system. According to Theorems \ref{s61thm3} and \ref{s61thm4}, $\vg$ has to solve the variational formulation $(\star)$. Hence we start with the solution of the variational formulation $(\star)$ and we prove that it is the solution to (\ref{s61thm2e}). The plan is the following:
\begin{itemize}
\item \textit{Step 1:} We extend the variational formulation $(\star)$ on a larger domain $\Omega_{0,e}$ with a solution denoted by $\vg_{e}$.
\item \textit{Step 2:} We prove that the solution $\vg_{e}$ to this new variational formulation is in $\textbf{H}^{2}$ in a neighbourhood of $\Gamma_{i}$.
\item \textit{Step 3:} We prove that the restriction of $\vg_{e}$ to the initial domain $\Omega_{0}$ is the solution $\vg$ to $(\star)$ which implies that $\vg$ is $\textbf{H}^{2}$ in a neighbourhood of $\Gamma_{i}$, and finally that $\vg\in \textbf{H}^{2}(\Omega_{0})$.
\item \textit{Step 4:} We prove that all the pressures associated with $\vg$ are in $H^{1}(\Omega_{0})$ and are constant on $\Gamma_{i,o}$.
\item \textit{Step 5:} We conclude by taking the pressure satisfying $q=0$ on $\Gamma_{i,o}$, so that the pair $(\vg,q)$ is the unique solution to (\ref{s61thm2e}).
\end{itemize}

\textit{Step 1:} Let $\eta^{0}_{e}$ be the function defined by
\[\eta^{0}_{e}:
\begin{cases}
\begin{array}{lr}
\eta^{0}(x)\,\text{ for all }x\in(0,L),\\
\eta^{0}(-x)\,\text{ for all }x\in (-L,0).
\end{array}
\end{cases}
\]
We recall that $\eta^{0}$ is in $H^{3}(0,L)$ and that $\eta^{0}(0)=\eta^{0}_{x}(0)=0$. Due to the even symmetry we have $\eta^{0}_{e}(0^{-})=\eta^{0}_{e}(0^{+})=0$, $\eta^{0}_{e,x}(0^{-})=\eta^{0}_{e,x}(0^{+})=0$, $\eta^{0}_{e,xx}(0^{-})=\eta^{0}_{e,xx}(0^{+})$ and thus we obtain $\eta^{0}_{e}\in H^{3}(-L,L)$ and the curve $\Gamma_{0,e}=\{(x,y)\in\mathbb{R}^{2}\mid x\in(-L,L),\,y=1+\eta^{0}_{e}(x)\}$ is $\mathcal{C}^{2}$. We set $\Omega_{0,e}=\{(x,y)\in\mathbb{R}^{2}\mid x\in(-L,L),\,0<y<1+\eta^{0}_{e}(x)\}$.

Let $\vg_e$ be the solution to
\[\nu\int_{\Omega_{0,e}}\nabla \vg_e:\nabla \psi=\int_{\Omega_{0,e}}\overline{\textbf{f}}_e\cdot \psi\,\text{ for all }\psi\in V_{e},
\]
where 
\begin{align*}
&V_{e}=\{\vg\in \textbf{H}^{1}(\Omega_{0,e})\mid \text{div }\vg=0\text{ in }\Omega_{0,e}\text{, }\vg=0\text{ on }\Gamma_{d,e}\text{, }v_2=0\text{ on }\Gamma_{i,o,e}\},\\
&\Gamma_{d,e}=(-L,L)\times \{0\}\cup\Gamma_{0,e},\,\,\,\Gamma_{i,e}=\{-L\} \times (0,1),\,\,\,\Gamma_{i,o,e}=\Gamma_{i,e}\cup\Gamma_{o},\\
\end{align*}
and $\overline{\textbf{f}}_e$ is the function defined by
\[\overline{\textbf{f}}_e:\begin{cases}
\begin{aligned}
\overline{\textbf{f}}_e&=\overline{\textbf{f}}\,\text{ in }\Omega_{0},\\
\overline{\textbf{f}}_e(x,y)&=\begin{pmatrix}
1&0\\
0&-1\\
\end{pmatrix}\overline{\textbf{f}}(-x,y)\,\text{ for all }(x,y)\in\Omega_{0,s},\\
\end{aligned}
\end{cases}
\]
with $\Omega_{0,s}=\{(x,y)\in\mathbb{R}^{2}\mid x\in(-L,0),\,0<y<1+\eta^{0}_{e}(x)\}$ .

\textit{Step 2:} We use cutoff functions to prove the $\textbf{H}^{2}$ regularity result near $\Gamma_{i}$. Let $\varphi$ be a function in $\mathcal{C}^{\infty}_{0}(\mathbb{R}^{2})$ such that $\varphi=1$ on $\Omega_{\varphi,1}$ and $\text{support}(\varphi)\subset \Omega_{\varphi,2}$, with $\Omega_{\varphi,1}$ and $\Omega_{\varphi,2}$ two open sets with smooth boundaries such that $\overline{\Omega_{\varphi,1}}\subset\overline{\Omega_{\varphi,2}}\subset \Omega_{0,e}$ and $\Omega_{\varphi,1}$ containing a neighbourhood of $\Gamma_i$. 

Let $\mathcal{Q}_{e}$ be a pressure associated to $\vg_{e}$. The pair $(\vg_{c},q_{c})=(\varphi\vg_{e},\varphi\mathcal{Q}_{e})$ satisfies, in $\textbf{H}^{-1}(\Omega_{\varphi,2})$,
\[-\nu\Delta \vg_{c} +\nabla q_{c}=-\nu\Delta\varphi\vg_{e} - 2\nu\nabla \vg_{e}\text{ }\nabla \varphi +\mathcal{Q}_{e}\nabla \varphi+\varphi\overline{\textbf{f}}_e.\]
Since $(\vg_{c},q_{c})$ belongs to $\textbf{H}^{1}_{0}(\Omega_{\varphi,2})\times L^{2}(\Omega_{\varphi,2})$, the previous equality implies that $(\vg_{c},q_{c})$ is a solution to the following Stokes equations (in the usual variational sense)
\begin{equation}\label{s61cutoff1}
\begin{aligned}
-\nu\Delta \vg_{c} +\nabla q_{c}&{}=-\nu\Delta\varphi\vg_{e} - 2\nu\nabla \vg_{e}\text{ }\nabla \varphi +\mathcal{Q}_{e}\nabla \varphi+\varphi\overline{\textbf{f}}_e\,\text{ in }\Omega_{\varphi,2},\\
\text{div }\vg_{c}&=\vg_{e}\cdot\nabla\varphi\,\text{ in }\Omega_{\varphi,2},\,\vg_c=0\,\text{ on }\partial\Omega_{\varphi,2}. \\
\end{aligned}
\end{equation}
We then use known results for Stokes equations with Dirichlet boundary conditions (see for example \cite[Chap IV, Theorem 5.8]{MR2986590}) to obtain $(\vg_{c},q_{c})\in \textbf{H}^{2}(\Omega_{\varphi,2})\times H^{1}(\Omega_{\varphi,2})$. As $(\vg_{c},q_{c})$ is equal to $(\vg_{e},\mathcal{Q}_{e})$ on $\Omega_{\varphi,1}$ we obtain the regularity result for $(\vg_{e},\mathcal{Q}_{e})$ in a neighbourhood of $\Gamma_{i}$.

\textit{Step 3:} We want to prove that the restriction to $\Omega_{0}$ of $\vg_{e}$ is the solution $\vg$ to the variational formulation $(\star)$. Using the Lax-Milgram lemma we know that $\vg_{e}$ satisfies
\begin{equation}\label{s61min}
\frac{1}{2}\nu\int_{\Omega_{0,e}}\vert\nabla\vg_{e}\vert^{2} - \int_{\Omega_{0,e}}\overline{\textbf{f}}_e\cdot\vg_{e}=\min_{\boldsymbol{\varphi}\in V_{e}}\left(\frac{1}{2}\nu\int_{\Omega_{0,e}}\vert\nabla\boldsymbol{\varphi}\vert^{2}-\int_{\Omega_{0,e}}\overline{\textbf{f}}_e\cdot\boldsymbol{\varphi}\right).\end{equation}
Hence, using the symmetry properties of $\overline{\textbf{f}}_e$ we can prove that the function $\vg_s$ defined by
\[\vg_s(x,y)=\begin{pmatrix}
1&0\\
0&-1\\
\end{pmatrix}\textbf{v}_{e}(-x,y)\,\text{ for all }(x,y)\in \Omega_{0,e},\]
is also a solution to the minimization problem (\ref{s61min}). As (\ref{s61min}) admits a unique solution we obtain that $\vg_s=\vg_{e}$. The symmetry properties and the regularity of $\vg_{e}$ imply that $v_{e,2}=0$ on $\Gamma_{i}$. We can now prove that the restriction to $\Omega_{0}$ of $\vg_{e}$ is the solution $\vg$ to $(\star)$. Let $\boldsymbol{\varphi}$ be a test function in $V$ and denote by $\boldsymbol{\varphi}_{e}$ the function defined by 
\[\boldsymbol{\varphi}_{e}:\begin{cases}
\begin{aligned}
\boldsymbol{\varphi}_{e}&=\boldsymbol{\varphi}\text{ on }\Omega_{0},\\
\boldsymbol{\varphi}_{e}(x,y)&=\begin{pmatrix}
1&0\\
0&-1\\
\end{pmatrix}\boldsymbol{\varphi}(-x,y)\,\text{ for all }(x,y)\in\Omega_{0,s}.
\end{aligned}
\end{cases}
\]
Thanks to the condition $\varphi_{2}=0$ on $\Gamma_{i,o}$ we notice that $\boldsymbol{\varphi}_{e}$ is in $\textbf{H}^{1}(\Omega_{0,e})$, and more precisely in $V_{e}$. Hence we can use $\boldsymbol{\varphi}_{e}$ as a test function in the variational formulation satisfied by $\vg_{e}$, we obtain
\[\nu\int_{\Omega_{0,e}}\nabla \vg_{e}:\nabla \boldsymbol{\varphi}_{e}=\int_{\Omega_{0,e}}\overline{\textbf{f}}_e\cdot\boldsymbol{\varphi}_{e}.\]
Using the symmetry properties of $\vg_{e}$, $\boldsymbol{\varphi}_{e}$ and $\overline{\textbf{f}}_e$ we have
\[
\int_{\Omega_{0,s}}\nabla \vg_{e}:\nabla \boldsymbol{\varphi}_{e}=\int_{\Omega_{0}}\nabla \vg_{e}:\nabla \boldsymbol{\varphi}_{e},
\]
and 
\[\int_{\Omega_{0,s}}\overline{\textbf{f}}_e\cdot\boldsymbol{\varphi}_{e}=\int_{\Omega_{0}}
\overline{\textbf{f}}_e\cdot\boldsymbol{\varphi}_{e}.\]
Hence,
\[\nu\int_{\Omega_{0}}\nabla \vg_{e}:\nabla \boldsymbol{\varphi}=\int_{\Omega_{0}}\overline{\textbf{f}}\cdot\boldsymbol{\varphi},\]
for all $\boldsymbol{\varphi}$ in $V$, which proves that the restriction to $\Omega_{0}$ of $\vg_{e}$ is the solution $\vg$ to the variational formulation $(\star)$. Hence $\vg$ is $\textbf{H}^{2}$ in a neighbourhood of $\Gamma_{i}$. The same technique works for the boundary $\Gamma_{o}$ which implies the regularity result on the whole domain $\Omega_{0}$.

\textit{Step 4:} Let $\mathcal{Q}$ be a pressure associated with $\vg$. The regularity of $\vg$ and the equality (in the sense of the distributions)
\[-\nu\Delta \vg + \nabla \mathcal{Q}=\overline{\textbf{f}},\]
imply that $\mathcal{Q}$ belongs to $H^{1}(\Omega_{0})$. We now have to prove that $\mathcal{Q}$ is equal to a constant on $\Gamma_{i,o}$. Thanks to the regularity of $(\vg,\mathcal{Q})$, the equality $-\Delta\vg + \nabla\mathcal{Q}=\overline{\textbf{f}}$ holds in $\textbf{L}^{2}(\Omega_{0})$. For all $\boldsymbol{\psi}$ in $V$ we have
\[
\int_{\Omega_{0}}\textbf{f}\cdot\boldsymbol{\psi} =\int_{\Omega_{0}}(-\nu\Delta\vg + \nabla \mathcal{Q})\cdot\boldsymbol{\psi}=\int_{\Omega_{0}}\nu\nabla \vg : \nabla \boldsymbol{\psi} + \int_{\Gamma_{i,o}}\mathcal{Q}(\boldsymbol{\psi}\cdot\normal),
\]
and, using the definition of $\vg$,
\[\int_{\Gamma_{i,o}}\mathcal{Q}(\boldsymbol{\psi}\cdot\normal)=0.\]
This implies that $\mathcal{Q}$ is constant on $\Gamma_{i,o}$. To see this, it is sufficient to prove that for all $\phi\in\mathcal{C}^{\infty}_{c}(\Gamma_{i,o})$ satisfying
\[\int_{\Gamma_{i,o}}\phi=0,\]
there exists $\boldsymbol{\psi}\in V$ such that $\boldsymbol{\psi}\cdot\normal=\phi$ on $\Gamma_{i,o}$. Let $\boldsymbol{\phi}$ be the function defined by 
\[\boldsymbol{\phi}:\begin{cases}
\begin{aligned}
\boldsymbol{\phi}&=0 \text{ on }\Gamma_{d},\\
\boldsymbol{\phi}&=\begin{pmatrix}
\phi\\
0\\
\end{pmatrix} \text{ on }\Gamma_{i,o}.\\
\end{aligned}
\end{cases}
\]
Using \cite[Lemma 2.2]{MR851383} the equations
\[
\begin{cases}
\begin{array}{rr}
\begin{aligned}
\text{div }\boldsymbol{\psi}&=0 \\
\boldsymbol{\psi}&=\boldsymbol{\phi}\\
\end{aligned}
&
\begin{aligned}
&\Omega_{0},\\
&\Gamma_{0},\\
\end{aligned}
\end{array}
\end{cases}
\]
admit a solution $\boldsymbol{\psi}$ in $\textbf{H}^{1}(\Omega_{0})$. Such a $\boldsymbol{\psi}$ belongs to $V$ and satisfies $\boldsymbol{\psi}\cdot\normal=\phi$ on $\Gamma_{i,o}$. Hence $\mathcal{Q}$ is constant on $\Gamma_{i,o}$.

\textit{Step 5:} Among the pressures $\mathcal{Q}$ associated with $\vg$ there exists a unique $q$ in $H^{1}(\Omega_{0})$ satisfying $q=0$ in $\Gamma_{i,o}$ in the sense of the trace for Sobolev functions. The pair $(\vg,q)$ in $\textbf{H}^{2}(\Omega_{0})\times H^{1}(\Omega_{0})$ is the unique solution to (\ref{s61thm2e}) and $(\ug,p)=(\vg,q)+(\textbf{w},H)$ is the unique solution to (\ref{s61e3}). The estimate on $(\ug,p)$ follows from classical estimate for the Stokes equations (\ref{s61cutoff1}) and Theorem \ref{s61thm1} to estimate $\textbf{w}$.
\end{proof}
According to Theorem \ref{s61thm5} the Stokes operator $A$ associated to (\ref{s61e3}) with homogeneous boundary condition is defined by
\[\mathcal{D}(A)=\textbf{H}^{2}(\Omega_{0})\cap V,\]
and for all $\ug\in \mathcal{D}(A)$, $A\ug=\nu\Pi\Delta\ug$.
\begin{thm}\label{s62thm1}The operator $(A,\mathcal{D}(A))$ is the infinitesimal generator of an analytic semigroup on $\textbf{V}^{0}_{n,\Gamma_{d}}(\Omega_{0})$. Moreover we have $\mathcal{D}(A^{1/2})=V$.
\end{thm}
\begin{proof}
The bilinear form associated with the operator $A$ defined by
\[\forall (\vg,\boldsymbol{\varphi})\in V\times V,\,\,\, a(\textbf{v},\boldsymbol{\varphi})=\nu\int_{\Omega_{0}}\nabla \textbf{v}:\nabla \boldsymbol{\varphi},\]
is continuous and coercive, hence \cite[Part 2, Theorem 2.2]{MR2273323} proves that the operator $A$ is the infinitesimal generator of an analytic semigroup. For the second part of the theorem we have, for all $\ug\in\mathcal{D}(A)$,
\[\norme{\ug}_{V}=\langle -A\ug,\ug\rangle=\norme{(-A)^{1/2}\ug}_{\textbf{V}^{0}_{n,\Gamma_{d}}(\Omega_{0})}.\]
By density, the previous equality is still true for $\ug\in V$ which concludes the proof.
\end{proof}
We now want to study (\ref{s61e3}) for weaker data using transposition method. The following lemma, used to solved non-zero divergence Stokes equations, is needed to obtain weak estimates on the pressure in Theorem \ref{s61lem1}.
\begin{lem}\label{non-zero-div}For all $\Phi\in H^{1}_{0}(\Omega_{0})$ the system
\begin{equation}\label{non-zero-div-eq}
\begin{array}{l}
\text{div }\wg=\Phi\,\text{ in }\Omega_{0},\\
\wg=0\,\text{ on }\Gamma_{d},\,w_2=0\,\text{ on }\Gamma_{i,o},\\
\end{array}
\end{equation}
admits a solution $\wg\in \textbf{H}^{2}(\Omega_{0})$ satisfying the estimate
\[\norme{\wg}_{\textbf{H}^{2}(\Omega_{0})}\leq C\norme{\Phi}_{H^{1}_{0}(\Omega_{0})}.\]
\end{lem}
\begin{proof}
If $\Phi$ has a zero average the result comes directly from \cite[Chap II.2, Lemma 2.3.1]{MR3013225}. This lemma gives the existence of a function $\wg\in \textbf{H}^{2}_{0}(\Omega_{0})$ such that $\text{div }\wg=\Phi$. In the general case, the idea is to find a pair $(\wg_{0},\Phi_{0})$ solution to (\ref{non-zero-div-eq}), where $\Phi_{0}$ has a non zero average, and to use it to come back to the previous framework.

Let $\delta>0$ be the constant defined by (\ref{delta}) in Theorem \ref{s61thm1} and $\rho\in\mathcal{C}^{\infty}(\mathbb{R})$ be a non zero non negative function compactly supported in $(0,\delta)$. Let $\theta\in\mathcal{C}^{\infty}(0,L)$ be such that $\theta=0$ near $0$ and $\theta=1$ near $L$. Define $\wg_{0}(x,y)=(\rho(y)\theta(x),0)^{T}$ for all $(x,y)\in\Omega_{0}$. The function $\wg_{0}$ is smooth and satisfies the boundary conditions in (\ref{non-zero-div-eq}). Finally, set $\Phi_{0}(x,y)=\text{div }\wg_{0}(x,y)=\rho(y)\theta'(x)$ for all $(x,y)\in\Omega_{0}$ and remark that $\Phi_{0}\in H^{1}_{0}(\Omega_{0})$ and
\[\int_{\Omega_{0}}\rho(y)\theta'(x)dxdy=\int_{0}^{\delta}\rho(y)dy>0.\]
We look for a solution to (\ref{non-zero-div-eq}) under the form $\textbf{w}=\widetilde{\wg}+c\wg_{0}$ with $c=\int_{\Omega_{0}}\Phi/\int_{\Omega_{0}}\Phi_{0}$. The function $\widetilde{\wg}$ needs to satisfy
\[
\begin{array}{l}
\text{div }\widetilde{\wg}=\Phi-c\Phi_{0}\,\text{ in }\Omega_{0},\\
\widetilde{\wg}=0\,\text{ on }\Gamma_{d},\,\widetilde{w}_2=0\,\text{ on }\Gamma_{i,o}.\\
\end{array}
\]
The function $\widetilde{\Phi}=\Phi-c\Phi_{0}$ is in $H^{1}_{0}(\Omega_{0})$ and has a zero average. The existence of $\widetilde{\wg}$ follows from \cite[Chap II.2, Lemma 2.3.1]{MR3013225}. To prove the estimate on $\wg$ remark that
\[c\leq \frac{\sqrt{\mu(\Omega_{0})}}{\int_{\Omega_{0}}\Phi_{0}}\norme{\Phi}_{H^{1}_{0}(\Omega_{0})}.
\]
\end{proof}
\begin{thm}\label{s61lem1}For all $(\textbf{f},\textbf{g},h)\in \textbf{L}^{2}(\Omega_{0})\times\mathcal{H}^{3/2}_{00}(\Gamma_{0})\times H^{1/2}(\Gamma_{i,o})$ the solution $(\ug,p)$ of the equation (\ref{s61e3}) satisfies the estimate
\begin{equation}\label{s61lem1est}
\norme{\ug}_{\textbf{L}^{2}(\Omega_{0})}+\norme{p}_{H^{-1}(\Omega_{0})}\leq C(\norme{\textbf{f}}_{(\textbf{H}^{2}(\Omega_{0}))'}+\norme{\textbf{g}}_{(\mathcal{H}^{1/2}(\Gamma_0))'}+\norme{h}_{(H^{3/2}(\Gamma_{i,o}))'}).
\end{equation}
\end{thm}
\begin{proof}
The fluid part estimate is similar to \cite[Lemma A.3]{MR2371113} using as test function the solution $(\boldsymbol{\Psi},\pi)$, given by Theorem \ref{s61thm5}, to
\begin{equation}\label{s61lem1e1}
\begin{aligned}
&-\nu\Delta \boldsymbol{\Psi} +\nabla \pi=\boldsymbol{\varphi},\,\,\,\text{div }\boldsymbol{\Psi}=0\,\text{ in }\Omega_{0},\\
&\boldsymbol{\Psi}=0\,\text{ on }\Gamma_{d},\,\,\,\Psi_2=0\,\text{ and }\,\pi=0\,\text{ on }\Gamma_{i,o},
\end{aligned}
\end{equation}
with $\boldsymbol{\varphi}\in \textbf{L}^{2}(\Omega_{0})$. Let us prove the pressure estimate. For all $\Phi\in H^{1}_{0}(\Omega_{0})$ consider the system
\begin{equation}\label{s61lem1e2}
\begin{aligned}
&-\nu\Delta \vg +\nabla q=0,\,\,\,\text{div }\vg=\Phi\,\text{ in }\Omega_{0},\\
&\vg=0\,\text{ on }\Gamma_{d},\,\,\,v_2=0\,\text{ and }\,q=0\,\text{ on }\Gamma_{i,o}.
\end{aligned}
\end{equation}
Using Lemma \ref{non-zero-div} and Theorem \ref{s61thm5} this system admits a unique solution $(\vg,q)$ in $\textbf{H}^{2}(\Omega_{0})\times H^{1}(\Omega_{0})$ which satisfies
\[\norme{\vg}_{\textbf{H}^{2}(\Omega_{0})}+\norme{q}_{H^{1}(\Omega_{0})}\leq C\norme{\Phi}_{H^{1}_{0}(\Omega_{0})}.\]
Using Green's formula the following computations hold
\begin{align*}
0&=\int_{\Omega_{0}}(-\nu\Delta\vg + \nabla q)\cdot\ug\\
&=-\nu\int_{\Omega_{0}}\Delta\ug\cdot\vg -\nu \int_{\partial\Omega_{0}}\ug\cdot(\nabla\vg\,\normal) +\nu \int_{\partial\Omega_{0}}\vg\cdot(\nabla\ug\,\normal) + \int_{\partial\Omega_{0}}q(\ug\cdot\normal)\\
&=\int_{\Omega_{0}}\textbf{f}\cdot\vg - \int_{\Omega_{0}}\nabla p\cdot\vg +\nu\int_{\partial\Omega_{0}}\ug\cdot(\nabla\vg\,\normal) + \int_{\Gamma_{0}}q(\ug\cdot\normal)\\
&=\int_{\Omega_{0}}\textbf{f}\cdot\vg + \int_{\Omega_{0}}p\,\Phi - \int_{\Gamma_{i,o}}h(\vg\cdot\normal) +\nu\int_{\partial\Omega_{0}}\textbf{u}\cdot(\nabla\vg\,\normal)+\int_{\Gamma_{0}}q(\textbf{g}\cdot\normal), \\
\end{align*}
and 
\begin{align*}
\int_{\partial\Omega_{0}}\textbf{u}\cdot(\nabla\vg\,\normal)&=\int_{\Gamma_{0}}\textbf{g}\cdot(\nabla\vg\,\normal) + \int_{\Gamma_{i,o}}\ug\cdot(\nabla\vg\,\normal)\\
&=\int_{\Gamma_{0}}g \,P_{2}(\nabla\vg\,\normal) + \int_{\Gamma_{i,o}}u_1\partial_1 v_1,\\
\end{align*}
where $\textbf{g}=(0,g)^{T}$ and $P_{2}$ is the vectorial projection on the second component. As $\partial_1 v_1+\partial_2 v_2=\Phi$ and $v_2=0$ on $\Gamma_{i,o}$ we notice that $\partial_1 v_1=\Phi$ on $\Gamma_{i,o}$ and as $\Phi\in H^{1}_{0}(\Omega_{0})$ we obtain $\partial_1 v_1=0$ on $\Gamma_{i,o}$. Finally
\begin{align*}
\left\vert\int_{\Omega_{0}}p\,\Phi\right\vert \leq{} &C(\norme{\textbf{f}}_{(\textbf{H}^{2}(\Omega_{0}))'}\norme{\vg}_{\textbf{H}^{2}(\Omega_{0})}+\norme{h}_{(H^{3/2}(\Gamma_{i,o}))'}\norme{\vg\cdot\normal}_{H^{3/2}(\Gamma_{i,o})}\\
&+\norme{\textbf{g}}_{(\mathcal{H}^{1/2}(\Gamma_0))'}\norme{P_{2}(\nabla\vg\,\normal)}_{\mathcal{H}^{1/2}(\Gamma_{0})}+\norme{\textbf{g}}_{(\mathcal{H}^{1/2}(\Gamma_{0}))'}\norme{P_{2}(q\normal)}_{\mathcal{H}^{1/2}(\Gamma_{0})}),\\
&{}\leq C(\norme{\textbf{f}}_{(\textbf{H}^{2}(\Omega_{0}))'}+\norme{\textbf{g}}_{(\mathcal{H}^{1/2}(\Gamma_0))'}+\norme{h}_{(H^{3/2}(\Gamma_{i,o}))'})\norme{\Phi}_{H^{1}_{0}(\Omega_{0})},
\end{align*}
which implies the pressure estimate.
\end{proof}
As for \cite[Theorem A.1]{MR2371113} we now define a notion of weak solutions for (\ref{s61e3}). For $(\textbf{f},\textbf{g},h)$ in $(\textbf{H}^{2}(\Omega_{0}))'\times (\mathcal{H}^{1/2}(\Gamma_0))'\times (H^{3/2}(\Gamma_{i,o}))'$ consider the following variational formulation:

Find $(\ug,p)\in \textbf{L}^{2}(\Omega_{0})\times H^{-1}(\Omega_{0})$ such that
\begin{equation}\label{s61var2}
\begin{aligned}
\int_{\Omega_{0}}\ug\cdot\boldsymbol{\varphi}={}&\langle \textbf{f},\boldsymbol{\Psi}\rangle_{(\textbf{H}^{2}(\Omega_{0}))',\textbf{H}^{2}(\Omega_{0})} + \langle h,\boldsymbol{\Psi}\cdot\normal\rangle_{(H^{3/2}(\Gamma_{i,o}))',H^{3/2}(\Gamma_{i,o})}\\
 - &\langle\textbf{g},P_{2}(\nabla\boldsymbol{\Psi}\,\normal)\rangle_{(\mathcal{H}^{1/2}(\Gamma_{0}))',\mathcal{H}^{1/2}(\Gamma_{0})} +\langle \textbf{g},P_{2}(\pi\normal)\rangle_{(\mathcal{H}^{1/2}(\Gamma_{0}))',\mathcal{H}^{1/2}(\Gamma_{0})},\\
\end{aligned}
\end{equation}
for all $\boldsymbol{\varphi}\in \textbf{L}^{2}(\Omega_{0})$ and $(\boldsymbol{\Psi},\pi)$ solution of (\ref{s61lem1e1}), and
\begin{equation}\label{s61var3}
\begin{aligned}
\langle p,\Phi\rangle_{H^{-1}(\Omega_{0}),H^{1}_{0}(\Omega_{0})}={}&-\langle \textbf{f},\vg\rangle_{(\textbf{H}^{2}(\Omega_{0}))',\textbf{H}^{2}(\Omega_{0})} + \langle h,\vg\cdot\normal\rangle_{(H^{3/2}(\Gamma_{i,o}))',H^{3/2}(\Gamma_{i,o})}\\
 -&\langle \textbf{g},P_{2}(\nabla\vg\,\normal)\rangle_{(\mathcal{H}^{1/2}(\Gamma_{0}))',\mathcal{H}^{1/2}(\Gamma_{0})}+\langle \textbf{g},P_{2}(q\normal)\rangle_{(\mathcal{H}^{1/2}(\Gamma_{0}))',\mathcal{H}^{1/2}(\Gamma_{0})},
\end{aligned}
\end{equation}
for all $\Phi\in H^{1}_{0}(\Omega_{0})$ and $(\vg,q)$ solution of (\ref{s61lem1e2}).
\begin{thm}\label{s61thm6}For all $(\textbf{f},\textbf{g},h)\in (\textbf{H}^{2}(\Omega_{0}))'\times (\mathcal{H}^{1/2}(\Gamma_0))'\times (H^{3/2}(\Gamma_{i,o}))'$ there exists a unique solution $(\ug,p)\in \textbf{L}^{2}(\Omega_{0})\times H^{-1}(\Omega_{0})$ of (\ref{s61e3}) in the sense of the variational formulation (\ref{s61var2})-(\ref{s61var3}). This solution satisfies the following estimate
\begin{equation}\label{s61thm6est}
\norme{\ug}_{\textbf{L}^{2}(\Omega_{0})}+\norme{p}_{H^{-1}(\Omega_{0})}\leq C(\norme{\textbf{f}}_{(\textbf{H}^{2}(\Omega_{0}))'}+\norme{\textbf{g}}_{(\mathcal{H}^{1/2}(\Gamma_0))'}+\norme{h}_{(H^{3/2}(\Gamma_{i,o}))'}).
\end{equation}
\end{thm}
\begin{proof}See \cite[Theorem A.1]{MR2371113}.
\end{proof}
\subsection{Unsteady Stokes equations}
Consider the unsteady Stokes equations
\begin{equation}\label{s62e1}
\begin{aligned}
&\textbf{u}_{t}-\nu\Delta\ug+\nabla p={\textbf{f}},\,\,\,\text{div }\textbf{u}=0\,\text{ in }Q_{T},\\
&\textbf{u}=\textbf{g}\,\text{ on }\Sigma^{0}_{T},\,\ug=0\,\text{ on }\Sigma^{b}_{T},\\
& u_2=0\,\text{ and }\,p=0\,\text{ on }\Sigma^{i,o}_{T},\\
&\textbf{u}(0)=\ug^{0}\,\text{ on }\Omega_{0}.\\
\end{aligned}
\end{equation}
As for the steady Stokes equations, a nonhomogeneous boundary condition on the pressure $p=h$ in (\ref{s62e1}) can be handled directly with a lifting, hence through this section we assume that $h=0$. We prove the existence and uniqueness of a solution to (\ref{s62e1}) in Theorem \ref{s62thm2}. Then we transform (\ref{s62e1}) to prove existence uniqueness and regularity result when the Dirichlet boundary condition $\textbf{g}$ is less regular (see Theorem \ref{s62thm3}). We use this result to prove Lemma \ref{lemmaP}. Finally we specify the regularity result used in the study of the fluid structure system in Theorem \ref{s62thm4} and we apply this result in Lemma \ref{A1continuous}.

Writing the equations satisfied by $\ug-D\textbf{g}$ and using standard semigroup techniques we obtain the following theorem. Remark that the assumption $\ug^0-D\textbf{g}(0)\in V$ is equivalent to $\ug^{0}\in \textbf{V}^{1}(\Omega_{0})$, $\ug^{0}=\textbf{g}$ on $\Gamma_{0}$ and $u^{0}_{2}=0$ on $\Gamma_{i,o}$.
\begin{thm}\label{s62thm2}For all $\textbf{g}\in L^{2}(0,T;\mathcal{H}^{3/2}_{00}(\Gamma_0))\cap H^{1}(0,T;(\mathcal{H}^{1/2}(\Gamma_0))')$, $\textbf{f}\in\textbf{L}^{2}(Q_{T})$ and $\ug^0\in \textbf{H}^{1}(\Omega_{0})$ satisfying the compatibility condition $\ug^0-D\textbf{g}(0)$ belongs to $V$, the equation (\ref{s62e1}) admits a unique solution $(\ug,p)\in \textbf{H}^{2,1}(Q_{T})\times L^{2}(0,T;H^{1}(\Omega_{0}))$. This solution satisfies the following estimate
\begin{align*}
&\norme{\ug}_{\textbf{H}^{2,1}(Q_{T})}+\norme{p}_{L^2(0,T;H^{1}(\Omega_{0}))}\\
&\leq C(\|\ug^0\|_{\textbf{H}^{1}(\Omega_{0})}+\norme{\textbf{g}}_{L^2(0,T;\mathcal{H}^{3/2}_{00}(\Gamma_0))}+\norme{\textbf{g}'}_{L^{2}(0,T;(\mathcal{H}^{1/2}(\Gamma_0)))'}+\norme{\textbf{f}}_{\textbf{L}^{2}(Q_{T})}).\\
\end{align*}
\end{thm}
We now want to study (\ref{s62e1}) for $\textbf{g}\in L^{2}(0,T;\mathcal{L}^{2}(\Gamma_0))$. We follow the approach of \cite{MR2371113}. The operator $A$, using extrapolation method, can be extended to an unbounded operator $\overset{\sim}{A}$ defined on $(\mathcal{D}(A^{*}))'$ with domain $\mathcal{D}(\overset{\sim}{A})=\textbf{V}^{0}_{n,\Gamma_{d}}(\Omega_{0})$.
\begin{defi}\label{s62d1}A function $\ug\in \textbf{L}^{2}(Q_{T})$ is called a weak solution to (\ref{s62e1}) if $\Pi\ug$ is a weak solution to the evolution equation
\begin{equation}\label{s62evol}\Pi\ug'=\overset{\sim}{A}\Pi\ug +(-\overset{\sim}{A})\Pi D\textbf{g}+\Pi\textbf{f},\,\,\Pi\ug(0)=\Pi\ug^0,
\end{equation}
and $(\mathbb{I}-\Pi)\ug$ is given by
\begin{equation}\label{s62Dg}(\mathbb{I}-\Pi)\ug=(\mathbb{I}-\Pi)D\textbf{g}\,\text{ in }\textbf{L}^{2}(Q_{T}).
\end{equation}
\end{defi}
Remark that $A=A^{*}$ (the operator $A$ is symmetric and onto from $\mathcal{D}(A)$ into $\textbf{V}^{0}_{n,\Gamma_{d}}(\Omega_{0})$). By definition to  a weak solution for (\ref{s62evol}) (see \cite{MR2273323}), $\Pi\ug\in L^{2}(0,T,\textbf{V}^{0}_{n,\Gamma_{d}}(\Omega_{0}))$ is solution to (\ref{s62evol}) if and only if for all $\Phi\in\mathcal{D}(A^{*})=\mathcal{D}(A)$ the map $t\mapsto \int_{\Omega_{0}}\Pi\ug\cdot\Phi$ belongs to $H^{1}(0,T)$ and
\begin{equation}\label{s62var1}
\frac{d}{dt}\int_{\Omega_{0}}\Pi\ug\cdot\Phi=\langle \overset{\sim}{A}\Pi\ug,\Phi \rangle_{\mathcal{D}(A)',\mathcal{D}(A)} + \langle -\overset{\sim}{A}\Pi D\textbf{g}, \Phi \rangle_{\mathcal{D}(A)',\mathcal{D}(A)}+\langle \Pi\textbf{f}, \Phi \rangle_{\mathcal{D}(A)',\mathcal{D}(A)}.
\end{equation}
Using Green formula we compute the adjoint of the operator $D$.
\begin{lem}\label{s62lem1}For all $\textbf{f}\in \textbf{L}^{2}(\Omega_{0})$ the adjoint operator $D^{*}$ of $D$ is defined by
\[D^{*}\textbf{f}=(-\nu\nabla \vg+q)\,\normal,\]
where $(\vg,q)\in \textbf{H}^{2}(\Omega_{0})\times H^{1}(\Omega_{0})$ is the solution to
\[
\begin{aligned}
&-\nu\Delta \vg +\nabla q=\textbf{f},\,\,\,\text{div }\vg=0\,\text{ in }\Omega_{0},\\
&\vg=0\,\text{ on }\Gamma_{d},\,\,\,v_{2}=0\,\text{ and }q=0\,\text{ on }\Gamma_{i,o}.
\end{aligned}
\]
\end{lem}
Using that $\overset{\sim}{A}^{*}=A$ on $\mathcal{D}(A)$ the variational formulation (\ref{s62var1}) becomes
\begin{align*}
\frac{d}{dt}\int_{\Omega_{0}}\Pi\ug\cdot\Phi&=\int_{\Omega_{0}}\Pi\ug\cdot A\Phi+\int_{\Gamma_0}\textbf{g}\cdot D^{*}(-A)\Phi+\int_{\Omega_{0}}\Pi\textbf{f}\cdot\Phi\\
&=\int_{\Omega_{0}}\Pi\ug\cdot A\Phi+\int_{\Gamma_0}\textbf{g}\cdot(-\nu\nabla \Phi +q)\,\normal,+\int_{\Omega_{0}}\Pi\textbf{f}\cdot\Phi,
\end{align*}
with $\nabla q=\nu(\mathbb{I}-\Pi)\Delta\Phi$. The previous equality follows from the uniqueness of the stationary Stokes system and the identity $-\nu\Delta \Phi + \nu(\mathbb{I}-\Pi)\Delta\Phi=-A\Phi$. Finally $\Pi\ug$ is a weak solution to \ref{s62evol} if and only if 
\begin{equation}\label{s62var2}
\frac{d}{dt}\int_{\Omega_{0}}\Pi\ug\cdot\Phi=\int_{\Omega_{0}}\Pi\ug\cdot A\Phi+\int_{\Gamma_0}\textbf{g}\cdot(-\nu\nabla \Phi +q)\,\normal+\int_{\Omega_{0}}\Pi\textbf{f}\cdot\Phi\,\text{ for all }\Phi\in\mathcal{D}(A).
\end{equation}
We can now state a theorem analogue to \cite[Theorem 2.3]{MR2371113}.
\begin{thm}\label{s62thm3}For all $\Pi\ug^0\in \textbf{V}^{0}_{n,\Gamma_{d}}(\Omega_{0})$, $\textbf{g}\in L^{2}(0,T;\mathcal{L}^{2}(\Gamma_0))$ and $\textbf{f}\in \textbf{L}^{2}(Q_{T})$ the equation (\ref{s62e1}) admits a unique weak solution $\ug$ in the sense Definition \ref{s62d1}. This solution satisfies the following estimate
\begin{equation}\label{s62thm3est1}
\begin{aligned}
&\norme{\Pi\ug}_{L^{2}(0,T;\textbf{V}^{1/2-\varepsilon}_{n,\Gamma_{d}}(\Omega_{0}))}+\norme{\Pi\ug}_{H^{1/4-\varepsilon/2}(0,T;\textbf{V}^{0}(\Omega_{0}))}+\norme{(\mathbb{I}-\Pi)\ug}_{L^{2}(0,T;\textbf{V}^{1/2}(\Omega_{0}))}\\
&{}\leq C\left(\|\Pi\ug^0\|_{\textbf{V}^{0}_{n,\Gamma_{d}}(\Omega_{0})}+\norme{\textbf{g}}_{L^{2}(0,T;\mathcal{L}^{2}(\Gamma_0))}+\norme{\Pi\textbf{f}}_{L^{2}(0,T;\textbf{V}^{0}_{n,\Gamma_{d}}(\Omega_{0}))}\right),\text{ for all }\varepsilon>0.\\
\end{aligned}
\end{equation}
\end{thm}
\begin{proof}See \cite[Theorem 2.3]{MR2371113}.
\end{proof}
As in \cite{MR2371113} we can prove that for $\textbf{g}\in L^{2}(0,T;\mathcal{H}^{3/2}_{00}(\Gamma_0))\cap H^{1}(0,T;\mathcal{H}^{-1/2}(\Gamma_0))$ a function $\ug$ is solution to (\ref{s62e1}) in the sense of Theorem \ref{s62thm2} if and only if $\ug$ is a weak solution to (\ref{s62e1})(in the sense of Definition \ref{s62d1}). The following theorem characterize the pressure.
\begin{thm}\label{s62thm5}For all $\textbf{g}\in L^{2}(0,T;\mathcal{H}^{3/2}_{00}(\Gamma_0))\cap H^{1}(0,T;(\mathcal{H}^{1/2}(\Gamma_0))')$, $\textbf{f}\in\textbf{L}^{2}(Q_{T})$ and $\ug^0\in \textbf{H}^{1}(\Omega_{0})$ satisfying the compatibility condition $\ug^0-D\textbf{g}(0)$ belongs to $V$, a pair $(\ug,p)\in \textbf{H}^{2,1}(Q_{T})\times L^{2}(0,T;H^{1}(\Omega_{0}))$ is solution of (\ref{s62e1}) if and only if 
\begin{align*}
&\Pi\ug'=A\Pi\ug + (-A)\Pi D\textbf{g}+\Pi \textbf{f},\,\,\,\ug(0)=\ug^{0},\\
&(I-\Pi)\ug=(I-\Pi)D\textbf{g},\,\,\,p=\rho -q_{t}+p_{\textbf{f}},
\end{align*}
where
\begin{itemize}
\item $q\in H^{1}(0,T;H^{1}(\Omega_{0}))$ is the solution to
\begin{equation}\label{eqq}\Delta q=0\,\text{ in }Q_{T},\,\,\,\rho=0\,\text{ on }\Sigma^{i,o}_{T},\,\,\,\frac{\partial q}{\partial \normal}=\textbf{g}\cdot\normal\,\text{ on }\Sigma^{0}_{T},\,\,\,\frac{\partial q}{\partial\normal}=0\,\text{ on }\Sigma^{b}_{T}.
\end{equation}
\item $\rho\in L^{2}(0,T;H^{1}(\Omega_{0}))$ is the solution to
\begin{equation}\label{eqrho}\Delta \rho =0\,\text{ in }Q_{T},\,\,\,\rho=0\,\text{ on }\Sigma^{i,o}_{T},\,\,\,\frac{\partial\rho}{\partial\normal}=\nu \Delta\Pi\ug\cdot\normal\,\text{ on }\Sigma^{d}_{T},
\end{equation}
where $\nu\Delta\Pi\ug\cdot\normal$ is in $L^{2}(0,T;H^{-1/2}(\Gamma_{d}))$ thanks to the divergence theorem.
\item $p_{\textbf{f}}\in L^{2}(0,T;H^{1}(\Omega_{0}))$ is given by the identity $(I-\Pi)\textbf{f}=\nabla p_{\textbf{f}}$.
\end{itemize}
\end{thm}
\begin{proof}
Writing $\ug=\Pi\ug + (\mathbb{I}-\Pi)\ug$ in Equation (\ref{s62e1}), we have
\[\ug_t -\nu\Delta \ug +\nabla p =\Pi\ug_{t}+(\mathbb{I}-\Pi)\ug_t -\nu\Delta \Pi\ug -\nu\Delta(\mathbb{I}-\Pi)\ug +\nabla p=0.\]
By definition of $(\mathbb{I}-\Pi)$ there exists $q\in H^{1}_{\Gamma_{i,o}}(\Omega_{0})$ such that $\nabla q=(\mathbb{I}-\Pi)\ug$. Using the condition $\text{div }\ug=0$ and $(\mathbb{I}-\Pi)\ug=(\mathbb{I}-\Pi)D\textbf{g}$ we obtain that $q$ is solution to (\ref{eqq}). As $\textbf{g}\in H^{1}(0,T;\mathcal{H}^{-1/2}(\Gamma_0))$ the function $q$ belongs to $H^{1}(0,T;H^{1}(\Omega_{0}))$.

The function $\Pi\ug$ is solution to the equation
\[\Pi\ug_{t}-\nu\Delta \Pi\ug +\nabla \rho=0,\]
with $\rho = p-\nu\Delta q + q_t=p+q_t$. Taking the divergence of the previous equation and the normal trace on $\Gamma_{d}$ (which is well defined as $\Delta\Pi\ug$ is in $L^{2}(0,T;\textbf{L}^{2}(\Omega_{0}))$ with a divergence equal to zero) we obtain (\ref{eqrho}) and $\rho\in L^{2}(0,T;H^{1}(\Omega_{0}))$.
\end{proof}

We conclude this section with a regularity result, coming from the interpolation of the regularity results stated in Theorem \ref{s62thm2} and Theorem \ref{s62thm3}, and an application to the operator $\mathcal{A}_{1}$ defined in Section 3.3.
\begin{thm}\label{s62thm4}For all $\textbf{g}\in L^{2}(0,T;\mathcal{H}^{1}_{0}(\Gamma_0))\cap H^{1/2}(0,T;\mathcal{L}^{2}(\Gamma_0))$, $\textbf{f}=0$ and $\Pi\ug^{0}=0$, the solution $\ug$ to (\ref{s62evol})-(\ref{s62Dg}) satisfies the estimate
\[\norme{\Pi\ug}_{\textbf{H}^{3/2-\varepsilon,3/4-\varepsilon/2}(Q_{T})}\leq C(\norme{\textbf{g}}_{L^{2}(0,T;\mathcal{H}^{1}_{0}(\Gamma_0))} + \norme{\textbf{g}}_{H^{1/2}(0,T;\mathcal{L}^{2}(\Gamma_0))}),\text{ for all }\varepsilon>0.
\]
\end{thm}
\begin{lem}\label{A1continuous}
The operator $(\mathcal{A}_{1},\mathcal{D}(\mathcal{A}_{1}))$ is the infinitesimal generator of a strongly continuous semigroup on $\textbf{H}$.
\end{lem}
\begin{proof}
The first part is to prove that the unbounded operator $(\overset{\sim}{\mathcal{A}_1},\mathcal{D}(\overset{\sim}{\mathcal{A}_1}))$, defined by 
\[\mathcal{D}(\overset{\sim}{\mathcal{A}_1})=\{(\Pi\ug,\eta_1,\eta_2)\in \textbf{V}^{1}_{n,\Gamma_{d}}(\Omega_{0})\times (H^{4}(\Gamma_s)\cap H^{2}_{0}(\Gamma_s))\times H^{2}_{0}(\Gamma_s)\mid \Pi\ug -\Pi D_s(\eta_2)\in V\}\]
and 
\[
\overset{\sim}{\mathcal{A}_1}=
\begin{pmatrix}
A&0&(-A)\Pi D_s\\
0&0&I\\
0&A_{\alpha,\beta}&\delta\Delta_s\\
\end{pmatrix},
\]
is the infinitesimal generator of a strongly continuous semigroup on $V^{-1}\times H_s$. Here, $V^{-1}$ is the dual of $V$ endowed with the norm
\[\textbf{v}\mapsto \left(\left\langle (-A)^{-1}\textbf{v},\textbf{v}\right\rangle_{V,V^{-1}}\right)^{1/2}.\]
This proof is similar to \cite[Theorem 3.5]{MR2745779}. Then we consider the evolution equation
\begin{equation}\label{s5e4}
\frac{d}{dt}\begin{pmatrix}
\Pi\ug\\
\eta_{1}\\
\eta_{2}\\
\end{pmatrix}
=\overset{\sim}{\mathcal{A}_1}
\begin{pmatrix}
\Pi\ug\\
\eta_{1}\\
\eta_{2}\\
\end{pmatrix},\,\,
\begin{pmatrix}
\Pi\ug(0)\\
\eta_{1}(0)\\
\eta_{2}(0)\\
\end{pmatrix}
=\begin{pmatrix}
\Pi\ug^0\\
\eta^{0}_{1}\\
\eta^{0}_{2}\\
\end{pmatrix}.
\end{equation}
The solution to (\ref{s5e4}) can be found in two steps. First we determine $(\eta_{1},\eta_{2})$ and then $\Pi\ug$. We recall that $(A_s,\mathcal{D}(A_s))$ is the infinitesimal generator of an analytic semigroup on $H_s$ (see \cite{MR971932}). Let $(\Pi\ug^{0},\eta^{0}_{1},\eta^{0}_{2})$ be in $V^{-1}\times H_s$. Using \cite[Chap 3, Theorem 2.2]{MR2273323} we obtain $\eta_{1}\in H^{3,3/2}(\Sigma^{s}_{T})$ and $\eta_{2}\in H^{1,1/2}(\Sigma^{s}_{T})$. Now let us assume that $(\Pi\ug^{0},\eta^{0}_{1},\eta^{0}_{2})\in \textbf{H}$. We have to solve
\[(\Pi\ug)'=A\Pi\ug + (-A)\Pi D_s(\eta_2),\,\, \Pi \ug(0)=\Pi\ug^{0}.\]
We split this equation in two parts $\Pi\ug=\Pi\ug_{1} + \Pi\ug_{2}$ with 
\[(\Pi\ug_{1})'=A\Pi\ug_{1} + (-A)\Pi D_s(\eta_2),\,\, \Pi \ug_{1}(0)=0,\]
and 
\[(\Pi\ug_2)'=A\Pi\ug_{2},\,\, \Pi \ug(0)=\Pi\ug^{0}.\]
Using Theorem \ref{s62thm4} we remark that $\Pi\ug_{1}\in \textbf{H}^{3/2-\varepsilon,3/4-\varepsilon/2}(Q_{T})$. For $\Pi\ug_{2}$, \cite[Chap 3, Theorem 2.2]{MR2273323} shows that $\Pi\ug_{2}\in L^{2}(0,T;V)\cap H^{1}(0,T;V^{-1})$. Interpolation result \cite[Theorem 3.1]{MR0350178} ensures that $\Pi\ug_{2}\in \mathcal{C}([0,T];\textbf{V}^{0}_{n,\Gamma_{d}}(\Omega_{0}))$. 

Hence $(\Pi\ug,\eta_{1},\eta_2)\in \mathcal{C}([0,T];\textbf{H})$ and the restriction to the semigroup $(e^{t\overset{\sim}{\mathcal{A}_1}})_{t\in\mathbb{R}^{+}}$ to $\textbf{H}$ is a strongly continuous semigroup on $\textbf{H}$. Finally we can verify that the infinitesimal generator associated with this restriction is exactly the operator $(\mathcal{A}_1,\mathcal{D}(\mathcal{A}_{1}))$.
\end{proof}
\subsection{Elliptic equations for the projector $\Pi$}
In this section we prove higher regularity result for an elliptic equation, which implies the regularity result on the projector $\Pi$ given in Lemma \ref{lemPif}.
\begin{lem}\label{lemmaPiu}Let $f$ be in $H^{1}(\Omega_{0})$ such that $f=0$ on $\Gamma_{i,o}$ and $g$ be in $H^{3/2}_{00}(\Gamma_{0})$. Then the elliptic equation
\begin{equation}\label{lemmaH3}
\begin{cases}
\begin{aligned}
&\Delta \rho =f\,\text{ in }\Omega_{0},\\
&\frac{\partial\rho}{\partial \normal}=g(1+(\eta^{0})^{2})^{-1/2}\,\text{ on }\Gamma_{0}\text{ and }\frac{\partial\rho}{\partial \normal}=0\,\text{ on }\Gamma_{b},\\
&\rho=0\,\text{ on }\Gamma_{i,o},\\
\end{aligned}
\end{cases}
\end{equation}
admits a unique solution $\rho\in H^{3}(\Omega_{0})$.
\end{lem}
\begin{proof}$H^{3}$ regularity far from the corners of $\Omega_{0}$ is obtained through classical arguments. To prove the $H^{3}$ regularity at the corners, say along $x=0$, we first perform a symmetry with respect to $x=0$ (step 1) and then a change of variables to transport the PDE on $(-L,L)\times (0,1)$ (step 2).

\textit{Step 1:} Using the notations of step 1 in the proof of Theorem \ref{s61thm5} for $\eta^{0}_{e}$, $\Gamma_{0,e}$, $\Omega_{0,s}$ and $\Omega_{0,e}$ we define $f_{e}$ and $g_{e}$ by
\[
\begin{array}{cc}
f_{e}:\begin{cases}
\begin{aligned}
&f_{e}=f\text{ in }\Omega_{0},\\
&f_{e}(x,y)=-f(-x,y)\,\text{ in }\Omega_{0,s},
\end{aligned}
\end{cases}& g_{e}:\begin{cases}
\begin{aligned}
&g_{e}=g\text{ in }\Gamma_{0},\\
&g_{e}(x,y)=g(-x,y)\,\text{ in }\Gamma_{0,e}\setminus\Gamma_{0}.
\end{aligned}
\end{cases}
\end{array}
\]
Assumptions on $f$ and $g$ ensure that $(f_{e},g_{e})$ is in $H^{1}(\Omega_{0,e})\times \textbf{H}^{3/2}(\Gamma_{0,e})$. Define $\rho_{e}$ by 
\[\rho_{e}:\begin{cases}
\begin{aligned}
&\rho_{e}=\rho\text{ in }\Omega_{0},\\
&\rho_{e}=-\rho(-x,y)\,\text{ for all }(x,y)\in\Omega_{0,s}.\\
\end{aligned}
\end{cases}
\]
Then $\rho_{e}\in H^{2}(\Omega_{0,e})$ and satisfies
\[
\begin{cases}
\begin{aligned}
&\Delta \rho_{e} =f_{e}\,\text{ in }\Omega_{0},\\
&\frac{\partial\rho_{e}}{\partial \normal}=g_{e}(1+(\eta^{0}_{e})^{2})^{-1/2}\,\text{ on }\Gamma_{0,e}\text{ and }\frac{\partial\rho_{e}}{\partial \normal}=0\,\text{ on }(-L,L)\times\{0\},\\
&\rho_e=0\,\text{ on }\left(\{-L\}\times(0,1)\right)\cup\Gamma_{o}.\\
\end{aligned}
\end{cases}
\]
\textit{Step 2:} Let $\Omega_{e}=(-L,L)\times(0,1)$ and $\varphi$ be the change of variables
\[\varphi:\begin{cases}
\begin{aligned}
&\Omega_{0,e}\longrightarrow \Omega_{e},\\
&(x,y)\mapsto (x,z)=\left(x,\frac{y}{1+\eta^{0}_{e}(x)}\right).\\
\end{aligned}
\end{cases}
\]
As in Theorem \ref{s61thm1} the function $\varphi$ transports $H^{3}(\Omega_{0,e})$ to $H^{3}(\Omega_{e})$. Hence it is sufficient to prove the $H^{3}$ regularity after transport. Let $\textbf{J}_{\varphi}$ be the Jacobian matrix of $\varphi$. Setting $\widetilde{\rho_{e}}=\rho\circ\varphi^{-1}$, $\widetilde{f_{e}}=\vert\textbf{J}_{\varphi}\vert^{-1} f_{e}\circ\varphi^{-1}$ and $\displaystyle\widetilde{g_{e}}(x,1)=g_{e}(x,1+\eta^{0}_{e
}(x))$ the function $\widetilde{\rho_{e}}$ is solution to
\begin{equation}\label{lemmaH3e1}\begin{cases}
\begin{aligned}
&\text{div}(A\nabla \widetilde{\rho_{e}})=\widetilde{f_{e}}\,\text{ in }\Omega_{e},\\
&A\nabla \widetilde{\rho}_{e}\cdot\normal=\widetilde{g_{e}}\,\text{ on }(-L,L)\times\{1\}\text{ and }A(x,z)\nabla \widetilde{\rho}_{e}\cdot\normal=0\,\text{ on }(-L,L)\times\{0\},\\
&\widetilde{\rho_e}=0\,\text{ on }\left(\{-L\}\times(0,1)\right)\cup\Gamma_{o},\\
\end{aligned}
\end{cases}
\end{equation}
where the matrix $A=(A_{i,j})_{1\leq i,j\leq 2}=\vert\text{det}(\textbf{J}_{\varphi})\vert^{-1}\textbf{J}_{\varphi}\textbf{J}_{\varphi}^{T}$ is uniformly positive definite symmetric with coefficients in $W^{1,\infty}\cap H^{2}$.

\textit{Step 3:} Deriving (\ref{lemmaH3e1}) with respect to $x$ shows that $\partial_{x}\widetilde{\rho_{e}}$ satisfies (with $\partial_{1}=\partial_{x}$ and $\partial_{2}=\partial_{z}$)
\begin{equation}\label{lemmaH3e2}
\text{div}(A\nabla(\partial_{x}\widetilde{\rho_e}))=\partial_{x}\widetilde{f_e}-F(A,\widetilde{\rho_e}),
\end{equation}
with
\begin{align*}
F(A,\widetilde{\rho_e})={}&\left(\partial_{11}A_{11}\right)\partial_{1}\widetilde{\rho_e}+\left(\partial_{1}A_{11}\right)\partial_{11}\widetilde{\rho_e}+\left(\partial_{11}A_{12}\right)\partial_{2}\widetilde{\rho_e}+\left(\partial_{1}A_{12}\right)\partial_{12}\widetilde{\rho_e}\\
&+\left(\partial_{12}A_{21}\right)\partial_{1}\widetilde{\rho_e}+\left(\partial_{1}A_{21}\right)\partial_{21}\widetilde{\rho_e}+\left(\partial_{12}A_{22}\right)\partial_{2}\widetilde{\rho_e}+\left(\partial_{1}A_{22}\right)\partial_{22}\widetilde{\rho_e},
\end{align*}
in the sense of the distributions on $\Omega_{e}$. From here on, we localize near $(0,1)$.

\textit{Step 4:} We use a bootstrap argument. The first step is to find an $L^{\infty}$ estimate on $\nabla\widetilde{\rho_e}$. In the right hand-side of (\ref{lemmaH3e2}) the least regular terms are under the form $\left(\partial_{11}A_{11}\right)\partial_{x}\widetilde{\rho_e}$ or $\left(\partial_{12}A_{22}\right)\partial_{z}\widetilde{\rho_e}$. Sobolev embeddings show that these terms are in $L^{r}$ for all $1<r<2$. Moreover the Neumann boundary condition involves $\partial_{x}\widetilde{g_e}-\left(\partial_{1}A_{21}\right)\partial_{x}\widetilde{\rho_e}-\left(\partial_{1}A_{22}\right)\partial_{z}\widetilde{\rho_e}$ where the least regular terms are traces of $W^{1,r}$ functions. Using the results of \cite{MR0125307} and \cite{MR0162050} we obtain that $\partial_{x}\widetilde{\rho_e}$ is in $W^{2,r}$. Then the embeddings $W^{2,r}\subset W^{1,r^{*}}\subset L^{\infty}$ with $r^{*}=\frac{2r}{2-r}>2$ show that the terms under the form $\left(\partial_{11}A_{11}\right)\partial_{x}\widetilde{\rho_e}$ are in $L^{2}$ and $\left(\partial_{1}A_{21}\right)\partial_{x}\widetilde{\rho_e}$ is in $H^{1/2}$ (on the boundary). Moreover using the equation (\ref{lemmaH3e1}) we obtain that $\partial_{zz}\widetilde{\rho_e}$ is in $L^{r^{*}}$ and thus $\partial_{z}\widetilde{\rho_e}\in W^{1,r^{*}}\subset L^{\infty}$. Finally the right hand-side is in $L^{2}$ and the Neumann boundary condition in $H^{1/2}$ and thus $\partial_{x}\widetilde{\rho_e}$ is $H^{2}$ near $(0,1)$. For the regularity with respect to $z$ we can use the equation (\ref{lemmaH3e1}) and $\widetilde{\rho_e}$ is $H^{3}$ in a neighbourhood of $(0,0)$.

\textit{Step 5:} The strategy applies for $(0,0)$. If we come back to the initial equation on the domain $\Omega_{0}$ we have proved that $\rho$ is $H^{3}$ near $\Gamma_{i}$. The same proof can be used for the regularity near $\Gamma_{o}$ and finally $\rho\in H^{3}(\Omega_{0})$.
\end{proof}

\bibliographystyle{plain}
\bibliography{biblio_1}

\begin{thebibliography}{10}

\bibitem{MR0450957}
R.~A. Adams.
\newblock {\em Sobolev spaces}.
\newblock Academic Press [A subsidiary of Harcourt Brace Jovanovich,
  Publishers], New York-London, 1975.
\newblock Pure and Applied Mathematics, Vol. 65.

\bibitem{MR0125307}
S.~Agmon, A.~Douglis, and L.~Nirenberg.
\newblock Estimates near the boundary for solutions of elliptic partial
  differential equations satisfying general boundary conditions. {I}.
\newblock {\em Comm. Pure Appl. Math.}, 12:623--727, 1959.

\bibitem{MR0162050}
S.~Agmon, A.~Douglis, and L.~Nirenberg.
\newblock Estimates near the boundary for solutions of elliptic partial
  differential equations satisfying general boundary conditions. {II}.
\newblock {\em Comm. Pure Appl. Math.}, 17:35--92, 1964.

\bibitem{MR2359448}
G.~Avalos and R.~Triggiani.
\newblock The coupled {PDE} system arising in fluid/structure interaction. {I}.
  {E}xplicit semigroup generator and its spectral properties.
\newblock In {\em Fluids and waves}, volume 440 of {\em Contemp. Math.}, pages
  15--54. Amer. Math. Soc., Providence, RI, 2007.

\bibitem{MR2379673}
G.~Avalos and R.~Triggiani.
\newblock Mathematical analysis of {PDE} systems which govern fluid-structure
  interactive phenomena.
\newblock {\em Bol. Soc. Parana. Mat. (3)}, 25(1-2):17--36, 2007.

\bibitem{MR2027753}
H.~Beir\~ao~da Veiga.
\newblock On the existence of strong solutions to a coupled fluid-structure
  evolution problem.
\newblock {\em J. Math. Fluid Mech.}, 6(1):21--52, 2004.

\bibitem{MR2273323}
A.~Bensoussan, G.~Da~Prato, M.~C. Delfour, and S.~K. Mitter.
\newblock {\em Representation and control of infinite dimensional systems}.
\newblock Systems \& Control: Foundations \& Applications. Birkh\"auser Boston,
  Inc., Boston, MA, second edition, 2007.

\bibitem{MR1900648}
J.~M. Bernard.
\newblock Non-standard {S}tokes and {N}avier-{S}tokes problems: existence and
  regularity in stationary case.
\newblock {\em Math. Methods Appl. Sci.}, 25(8):627--661, 2002.

\bibitem{MR1978563}
J.~M. Bernard.
\newblock Time-dependent {S}tokes and {N}avier-{S}tokes problems with boundary
  conditions involving pressure, existence and regularity.
\newblock {\em Nonlinear Anal. Real World Appl.}, 4(5):805--839, 2003.

\bibitem{MR2986590}
F.~Boyer and P.~Fabrie.
\newblock {\em Mathematical tools for the study of the incompressible
  {N}avier-{S}tokes equations and related models}, volume 183 of {\em Applied
  Mathematical Sciences}.
\newblock Springer, New York, 2013.

\bibitem{MR971932}
S.~P. Chen and R.~Triggiani.
\newblock Proof of extensions of two conjectures on structural damping for
  elastic systems.
\newblock {\em Pacific J. Math.}, 136(1):15--55, 1989.

\bibitem{MR1308419}
C.~Conca, F.~Murat, and O.~Pironneau.
\newblock The {S}tokes and {N}avier-{S}tokes equations with boundary conditions
  involving the pressure.
\newblock {\em Japan. J. Math. (N.S.)}, 20(2):279--318, 1994.

\bibitem{MR851383}
V.~Girault and P.-A. Raviart.
\newblock {\em Finite element methods for {N}avier-{S}tokes equations},
  volume~5 of {\em Springer Series in Computational Mathematics}.
\newblock Springer-Verlag, Berlin, 1986.
\newblock Theory and algorithms.

\bibitem{MR3466847}
C.~Grandmont and M.~Hillairet.
\newblock Existence of global strong solutions to a beam-fluid interaction
  system.
\newblock {\em Arch. Ration. Mech. Anal.}, 220(3):1283--1333, 2016.

\bibitem{MR2765696}
J.~Lequeurre.
\newblock Existence of strong solutions to a fluid-structure system.
\newblock {\em SIAM J. Math. Anal.}, 43(1):389--410, 2011.

\bibitem{MR0350177}
J.-L. Lions and E.~Magenes.
\newblock {\em Non-homogeneous boundary value problems and applications. {V}ol.
  {I}}.
\newblock Springer-Verlag, New York-Heidelberg, 1972.
\newblock Translated from the French by P. Kenneth, Die Grundlehren der
  mathematischen Wissenschaften, Band 181.

\bibitem{MR0350178}
J.-L. Lions and E.~Magenes.
\newblock {\em Non-homogeneous boundary value problems and applications. {V}ol.
  {II}}.
\newblock Springer-Verlag, New York-Heidelberg, 1972.
\newblock Translated from the French by P. Kenneth, Die Grundlehren der
  mathematischen Wissenschaften, Band 182.

\bibitem{MR3017292}
B.~Muha and S.~Cani\'c.
\newblock Existence of a weak solution to a nonlinear fluid-structure
  interaction problem modeling the flow of an incompressible, viscous fluid in
  a cylinder with deformable walls.
\newblock {\em Arch. Ration. Mech. Anal.}, 207(3):919--968, 2013.

\bibitem{MR0512912}
A.~Pazy.
\newblock {\em Semi-groups of linear operators and applications to partial
  differential equations}.
\newblock Department of Mathematics, University of Maryland, College Park, Md.,
  1974.
\newblock Department of Mathematics, University of Maryland, Lecture Note, No.
  10.

\bibitem{MR2371113}
J.-P. Raymond.
\newblock Stokes and {N}avier-{S}tokes equations with nonhomogeneous boundary
  conditions.
\newblock {\em Ann. Inst. H. Poincar\'e Anal. Non Lin\'eaire}, 24(6):921--951,
  2007.

\bibitem{MR2745779}
J.-P. Raymond.
\newblock Feedback stabilization of a fluid-structure model.
\newblock {\em SIAM J. Control Optim.}, 48(8):5398--5443, 2010.

\bibitem{MR3013225}
H.~Sohr.
\newblock {\em The {N}avier-{S}tokes equations}.
\newblock Modern Birkh\"auser Classics. Birkh\"auser/Springer Basel AG, Basel,
  2001.
\newblock An elementary functional analytic approach.

\end{thebibliography}

\end{document}